\documentclass[11pt,a4paper]{amsart}

\usepackage[utf8]{inputenc} 
\usepackage{enumitem,hyperref}       
\usepackage{url}            
\usepackage{amsfonts,amsmath,amsthm,amssymb}       
\usepackage{amsaddr}
\usepackage{tikz}
\usetikzlibrary{calc,shapes,snakes,patterns,matrix,backgrounds}
\usetikzlibrary{arrows.meta}
\usepackage{pgfplots}
\usepackage{capt-of}
\usepackage{mathtools}
\usepackage{algorithm}
\usepackage{algpseudocode}
\usepackage{xparse}
\usepackage{listings}  
\usepackage[numbered,framed]{mcode}
\usepackage{xcolor}
\usepackage{textcomp}
\usepackage[margin=0.75in]{geometry}
\usepackage{tabularx}

\DeclareMathOperator{\COARSEN}{COARSEN}
\DeclareMathOperator{\refe}{ref}

\definecolor{pastelred}{rgb}{1.0, 0.41, 0.38}
\definecolor{lightgreen}{rgb}{0.56, 0.93, 0.56}
\definecolor{lightblue}{rgb}{0.53, 0.81, 0.98}
\definecolor{markred}{rgb}{0.8, 0.0, 0.0}
\definecolor{lemon}{rgb}{1.0, 1.0, 0.4}
\definecolor{rednodes}{rgb}{0.7, 0.11, 0.11}
\def\revision#1{{\color{black}#1}} 
\def\matlab#1{{\small\mcode{#1}}}

\newenvironment{axioms}
 {\enumerate[label=\textbf{R\arabic*.}, ref=\textbf{R\arabic*.},align=left, labelwidth=1ex]}
 {\endenumerate}
\newtheorem*{definition}{Definition}
\newtheorem{theorem}{Theorem}

\newtheorem{remark}{Remark}

\tikzset{
  schraffiert/.style={pattern=north west lines,pattern color=#1},
  schraffiert/.default=black
}
\tikzset{cross/.style={cross out, draw=black, minimum size=2*(#1-\pgflinewidth), inner sep=0pt, outer sep=0pt},
cross/.default={3pt}}

\title[A coarsening algorithm on adaptive red-green-blue refined meshes]{A coarsening algorithm on\\ adaptive red-green-blue refined meshes}
\author{Stefan A.~Funken \and Anja Schmidt}
\address{Ulm University, Institute for Numerical Mathematics}
\email{stefan.funken@uni-ulm.de, anja.schmidt@uni-ulm.de}
\thanks{\emph{Acknowledgement.} The authors would like to thank Mazen Ali, Stefan Ehard, and Marcus Heitel for comments that greatly improved this work.}

\begin{document}
\maketitle

\begin{abstract}
Adaptive meshing is a fundamental component of adaptive finite element methods. This includes refining and coarsening meshes locally. In this work, we are concerned with the red-green-blue refinement strategy \revision{in two dimensions} and its counterpart -- coarsening. In general, coarsening algorithms are mostly based on an explicitly given refinement history. In this work, we present a coarsening algorithm on adaptive red-green-blue meshes \revision{in two dimensions} without explicitly knowing the refinement history. To this end, we examine the local structure of these meshes, find an easy-to-verify criterion to \revision{adaptively} coarsen red-green-blue meshes and prove that this criterion generates meshes with the desired properties. We present a MATLAB implementation built on the red-green-blue refinement routine of the \texttt{ameshref}-package \cite{web,funkenschmidt,ameshref}. 
\end{abstract}

\section{Introduction}
Adaptive meshing is a popular tool to efficiently solve partial differential equations where solutions exhibit local singularities \cite{nochetto2009theory}. In time-dependent problems, singularities, interfaces and forces may move or change in time. This requires coarsening meshes locally. Otherwise the algorithm's efficiency would decrease with time since degrees of freedom needed for an earlier time step are not released as the singularity or interface progresses. To this end, it is common to deploy coarsening algorithms to maintain the adaptive efficiency \cite{bartels,alberta}. Furthermore, coarsening routines are used in multigrid techniques where a sequence of coarse and fine meshes is needed \cite{multigrid,gooch}.

Local geometric refinement is a major part of adaptive meshing. The goal is to reduce the element size by adding further nodes to a given mesh. Several refinement strategies are known which have desired properties and are therefore well suited for adaptive meshing.  An overview and a list of public code is provided by Schneiders in \cite{schneiders}. Local coarsening is the counterpart of local refinement and is thus also an important part of adaptive meshing. There are different approaches to coarsening. Local coarsening refers to deleting nodes from a given mesh to increase the element size. Possible approaches are based on edge collapsing \cite{bankxu,gooch}, centroidal Voronoi tessellations \cite{Shu2009} or the refinement history \cite{bartels,chenzhang,kossa,alberta}. The latter approach aims to invert the refinement based on the refinement history. Desired properties such as the inscribed ball condition \cite{ciarlet} are automatically fulfilled during coarsening. The first two approaches, in contrast, do not use the refinement history. Desired properties are thus not automatically preserved within the coarsening process.

Early works on coarsening based on the refinement history refer to the hierarchical structure of the refinement and use this information to coarsen elements to their corresponding father element \cite{kossa,alberta}. Chen and Zhang proposed a new concept to identify admissible-to-coarsen nodes without explicitly knowing the hierarchical structure for the newest vertex bisection (NVB) \cite{chenzhang}. Bartels and Schreier extended this result to three dimensions \cite{bartels}. To the best of our knowledge, this has not been done for other refinement strategies \revision{of triangular meshes}. To this end, we bridge the gap and present a new criterion to adaptively coarsen meshes generated by the red-green-blue (RGB) refinement strategy \revision{in two dimensions} introduced in \cite{carstensen2004} and implemented in the \texttt{ameshref}-package \cite{web,funkenschmidt}. The only information we use to describe a mesh is the element-connectivity and the coordinates of the vertices. No information about neighbours or father-child connections is stored. A key observation within this paper is that this minimal data structure can also be kept for coarsening, i.e., no additional information is needed to coarsen the meshes. However, as hierarchical data is non-present, the determination of nodes that can be eliminated -- while preserving desired properties -- is more difficult. We present an algorithm that determines those ``admissible'' nodes. 

This paper is organized as follows. In Section~\ref{sect:prelim}, we introduce some notations and definitions and shortly present the red-green-blue refinement. \revision{We highlight the requirements for the data structure of an RGB refinement implementation such that the proposed RGB coarsening algorithm can be realized based on this implementation. We further present the RGB} implementation in the \texttt{ameshref}-package as we build our coarsening routine on this code. In Section~\ref{sect:requirements}, we focus on coarsening requirements \revision{and compare the RGB refinement to the newest vertex bisection. For newest vertex bisection, a coarsening strategy is already known. Thus, we show the limitations of this approach for RGB and, in Section~\ref{sect:algorithm}, adapt it in a way that some ideas can be carried over and the limitations motivate the RGB coarsening algorithm presented and examined}. In Section~\ref{sect:Matlab}, we focus on the efficient implementation in MATLAB by use of vectorization and conclude with numerical experiments presented in Section~\ref{sect:experiments}. 
\section{Preliminaries}\label{sect:prelim}

Let $\Omega$ be a polygonal domain in $\mathbb{R}^2$. An element $T \subset \mathbb{R}^2$ is a triangle including edges. We call $\mathcal{T}$ a \emph{triangulation} of $\Omega$ if 
\begin{itemize}
\item $\mathcal{T}$ is a finite set of elements $T$ with positive area $|T|>0$ ,
\item the union of all elements in $\mathcal{T}$ covers the closure $\overline{\Omega}$,
\item for $T_i,T_j \in \mathcal{T}$ with $T_i\neq T_j$ for $i\neq j$ holds $\mathring{T_i}\cap \mathring{T_j}=\emptyset$, where $\mathring{T}$ denotes the interior of $T$.
\end{itemize}
We denote the set of all vertices of a triangulation $\mathcal{T}$ with $\mathcal{N}$, and the set of all edges with $\mathcal{E}$. With this, $\mathcal{N}(T) \coloneqq \left\{v \in \mathcal{N}~|~v\in T\right\}$ is the set of nodes of an element $T \in \mathcal{T}$. Analogously, $\mathcal{E}(T) \coloneqq \left\{e \in \mathcal{E}~|~e \subset \partial T\right\}$ is the set of edges of an element $T \in \mathcal{T}$. \revision{We index $\mathcal{T}$ and $\mathcal{N}$ with a zero when we reference to the initial triangulation $\mathcal{T}_0$ and the nodes $\mathcal{N}_0$ in the initial triangulation.}
We call $\mathcal{T}$ a \emph{\revision{conforming} triangulation} of $\Omega$ if additionally
\begin{itemize}
\item for all $T_i,T_j$ with $T_i\neq T_j$ for $i\neq j$ holds that $T_i\cap T_j$ is the empty set, a common node or a common edge.
\end{itemize}
The aforementioned definition prevents a triangulation from having hanging nodes. A node $v \in \mathcal{N}$ is called \emph{hanging node} if for some element $K \in \mathcal{T}$ it satisfies $v \in \partial K \setminus \mathcal{N}(K)$. We define an \emph{extended \revision{conforming} triangulation} $(\mathcal{T},\refe_\mathcal{T})$ where $\mathcal{T}$ is a \revision{conforming} triangulation and $\refe_\mathcal{T}$ is a mapping $\refe_{\mathcal{T}}: \mathcal{T} \rightarrow \mathcal{E}(\mathcal{T})$ that assigns a \emph{reference edge} to each triangle $T \in \mathcal{T}$ such that for $T \in \mathcal{T}$ holds: $\refe_\mathcal{T}(T) \in \mathcal{E}(T)$. 
For a triangle $T\in \mathcal{T}$ with reference edge $\refe_{T}$, a \emph{refinement} $(r(T),\refe_{r(T)})$ is a finite set of triangles such that
\begin{itemize}
\item for all $\tilde{T} \in r(T)$ holds $\tilde{T} \subset T$,
\item $\bigcup\limits_{\tilde{T} \in r(T)}\tilde{T}=T$,
\item for all $\tilde{T},\hat{T} \in r(T)$ with $\tilde{T}\neq \hat{T}$ holds that $\tilde{T}\cup\hat{T}$ is the empty set, a common node or a common edge, and
\item for all $\tilde{T} \in r(T)$ a new reference edge $\refe_{r(T)}: r(T) \rightarrow \mathcal{E}(r(T))$ is assigned such that for $\tilde{T} \in r(T)$ holds $\refe_{r(T)}(\tilde{T}) \in \mathcal{E}(\tilde{T})$.
\end{itemize} 
We call $(\tilde{\mathcal{T}},\refe_{\tilde{\mathcal{T}}})$ a \emph{refinement of a triangulation} $\mathcal{T}$ if 
\begin{itemize}
\item each $(T,\refe_T) \in (\mathcal{T},\refe_\mathcal{T})$ is refined to $(r(T),\refe_{r(T)})$, and
\item the resulting triangulation $(\tilde{\mathcal{T}},\refe_{\tilde{\mathcal{T}}})$ is an extended \revision{conforming} triangulation.
\end{itemize}
The last point in particular ensures that the resulting triangulation does not have any hanging nodes. Eliminating hanging nodes by refining further elements is called CLOSURE. For further details we refer to Section~\ref{sect:RGBrefine} and \cite{funkenschmidt}.

In this work, we are concerned with the red-green-blue refinement \revision{in two dimensions}.
\begin{definition}[red-green-blue refinement (RGB), cf. \cite{carstensen2004}] We call a refinement $(r(T),\refe_{r(T)})$ of a triangle $T$ with reference edge $\refe(T)$ a 
\revision{
\begin{itemize}
\item \emph{red} refinement if triangle $T$ is divided into four subtriangles by joining the midpoints of its edges;
\item \emph{green} refinement if triangle $T$ is divided into two subtriangles by joining the midpoint of the reference edge $r(T)$ to the vertex opposite to this edge;
\item \emph{blue} refinement if triangle $T$ is divided into three subtriangles by joining the midpoint of the reference edge $r(T)$ to the vertex opposite to this edge and to the midpoint of one of the other edges;
\end{itemize} and for each subtriangle a new reference edge is assigned according to }Figure~\ref{fig:TrefineRGBpattern}.\end{definition}
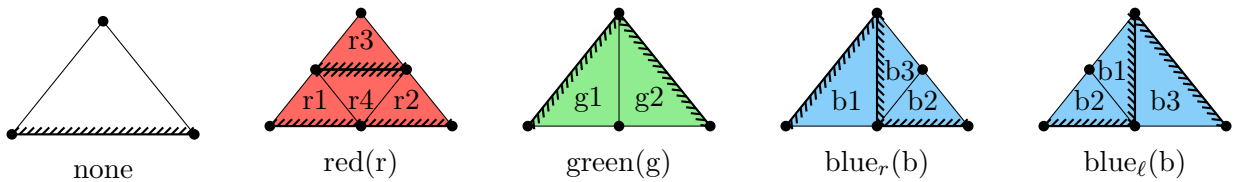
\begin{figure}
\begin{minipage}[t]{0.19\linewidth}
\hspace*{6mm}\begin{tikzpicture}
[ interface/.style={
        postaction={draw,decorate,decoration={border,angle=45,
                    amplitude=0.13cm,segment length=1mm}}},
    ]
\coordinate (1) at (-1.2,0);
\coordinate (2) at (1.2,0);
\coordinate (3) at (0,1.5);
\draw(1)--(2) -- (3)--(1);
\draw [interface, thick](1) --(2);
\fill (1) circle (2pt); 
\fill (2) circle (2pt); 
\fill (3) circle (2pt); 
\node at (0,-0.5) {none};
\end{tikzpicture}
\end{minipage}
\begin{minipage}[t]{0.19\linewidth}
\hspace*{6mm}\begin{tikzpicture}[ interface/.style={
        postaction={draw,decorate,decoration={border,angle=45,
                    amplitude=0.13cm,segment length=1mm}}},
    ]
\coordinate (1) at (-1.2,0);
\coordinate (2) at (1.2,0);
\coordinate (3) at (0,1.5);
\draw[fill=pastelred] (1)--(2) -- (3)--(1);
\draw[interface,thick] (1) --(2);
\fill (1) circle (2pt); 
\fill (2) circle (2pt); 
\fill (3) circle (2pt); 
\fill ($(1)!0.5!(2)$) circle (2pt);
\fill ($(2)!0.5!(3)$) circle (2pt);
\fill ($(1)!0.5!(3)$) circle (2pt);
\draw ($(1)!0.5!(2)$) -- ($(2)!0.5!(3)$)--($(1)!0.5!(3)$)--($(1)!0.5!(2)$);
\draw[interface,thick] ($(2)!0.5!(3)$)--($(1)!0.5!(3)$);
\draw[interface,thick] ($(1)!0.5!(3)$)--($(2)!0.5!(3)$);
\node at (-0.6,0.35) {r1};
\node at (0.6,0.35) {r2};
\node at (0,1.15) {r3};
\node at (0,0.35) {r4};
\node at (0,-0.5) {red(r)};
\end{tikzpicture}
\end{minipage}
\begin{minipage}[t]{0.19\linewidth}
\hspace*{6mm}\begin{tikzpicture}[ interface/.style={
        postaction={draw,decorate,decoration={border,angle=45,
                    amplitude=0.13cm,segment length=1mm}}},
    ]
\coordinate (1) at (-1.2,0);
\coordinate (2) at (1.2,0);
\coordinate (3) at (0,1.5);
\draw[fill=lightgreen] (1) --(2)--(3)--(1);
\draw[interface,thick] (2) -- (3);
\draw[interface,thick] (3) -- (1);
\fill (1) circle (2pt); 
\fill (2) circle (2pt); 
\fill (3) circle (2pt); 
\fill ($(1)!0.5!(2)$) circle (2pt);
\draw ($(1)!0.5!(2)$) --(3);
\node at (-0.4,0.35) {g1};
\node at (0.4,0.35) {g2};
\node at (0,-0.5) {green(g)};
\end{tikzpicture}
\end{minipage}
\begin{minipage}[t]{0.19\linewidth}
\hspace*{6mm}\begin{tikzpicture}[ interface/.style={
        postaction={draw,decorate,decoration={border,angle=45,
                    amplitude=0.13cm,segment length=1mm}}},
    ]
\coordinate (1) at (-1.2,0);
\coordinate (2) at (1.2,0);
\coordinate (3) at (0,1.5);
\draw[fill=lightblue] (1)--(2) -- (3)--(1);
\draw[interface,thick] (3) -- (1);
\fill (1) circle (2pt); 
\fill (2) circle (2pt); 
\fill (3) circle (2pt); 
\fill ($(1)!0.5!(2)$) circle (2pt);
\fill ($(2)!0.5!(3)$) circle (2pt);
\draw ($(2)!0.5!(3)$) -- ($(1)!0.5!(2)$);
\draw[interface,thick] (3)--($(1)!0.5!(2)$);
\draw[interface,thick] ($(1)!0.5!(2)$) --(2);
\node at (0,-0.5) {blue$_r$(b)};
\node at (-0.4,0.35) {b1};
\node at (0.6,0.35) {b2};
\node at (0.3,0.75) {b3};
\end{tikzpicture}
\end{minipage}
\begin{minipage}[t]{0.19\linewidth}
\hspace*{6mm}\begin{tikzpicture}[ interface/.style={
        postaction={draw,decorate,decoration={border,angle=45,
                    amplitude=0.13cm,segment length=1mm}}},
    ]
\coordinate (1) at (-1.2,0);
\coordinate (2) at (1.2,0);
\coordinate (3) at (0,1.5);
\draw[fill=lightblue] (1)--(2) -- (3)--(1);
\draw[interface,thick] (2) -- (3);
\fill (1) circle (2pt); 
\fill (2) circle (2pt); 
\fill (3) circle (2pt); 
\fill ($(1)!0.5!(2)$) circle (2pt);
\fill ($(1)!0.5!(3)$) circle (2pt);
\draw ($(1)!0.5!(3)$) -- ($(1)!0.5!(2)$);
\draw[interface,thick] ($(1)!0.5!(2)$)--(3);
\draw[interface,thick] (1)-- ($(1)!0.5!(2)$);
\node at (0,-0.5) {blue$_\ell$(b)};
\node at (0.4,0.35) {b3};
\node at (-0.6,0.35) {b2};
\node at (-0.3,0.75) {b1};
\end{tikzpicture}
\end{minipage}

\caption{\revision{From left to right: Initial triangle (none) and its possible }refinements red, green, blue$_r$, and blue$_\ell$. Reference edges are highlighted by hatched lines. \revision{The letters r,g,b denote the type of refinement and} the numbers indicate the storage sequence of the \revision{newly created} elements.}
\label{fig:TrefineRGBpattern}
\end{figure}

Reference edges are chosen such that during the refinement process all formed triangles starting from an initial triangle $T_0$ fall into at most four similarity classes \cite{aschmidt,sewell}. This ensures that degeneracies are avoided and Ciarlet's inscribed ball condition \cite{ciarlet} is satisified for a family of triangulations $\mathcal{T}_h$ formed by the refinements. This property is often referred to as \emph{shape regularity} of a triangulation. The assignment of reference edges is clearly prescribed through the refinement process. There still remains the question of how to select the reference edges in the initial triangulation $\mathcal{T}_0$. Obviously, the choice has some impact on the locality of the adaptive mesh. An intuitive choice is, e.g., the longest edge. Further possibilities are discussed in \cite{carstensen2004,stevensoncompletion}. In depictions we refrain from labeling the reference edges whenever it is irrelevant for the context.

\subsection{\revision{Requirements for the Data Structure of an RGB Refinement Implementation}}

\begin{figure}[!htb]
\begin{center}
\begin{minipage}[t]{0.23\linewidth}
\hspace*{6mm}\begin{tikzpicture}[ interface/.style={
        postaction={draw,decorate,decoration={border,angle=45,
                    amplitude=0.13cm,segment length=1mm}}},
    ]
\coordinate (1) at (-1.2,0);
\coordinate (2) at (1.2,0);
\coordinate (3) at (0,1.5);
\draw(1)--(2) -- (3)--(1);
\draw[interface,thick] (1) --(2);
\fill (1) circle (2pt); 
\fill (2) circle (2pt); 
\fill (3) circle (2pt); 
\draw ($(1)!0.5!(2)$) -- ($(2)!0.5!(3)$)--($(1)!0.5!(3)$)--($(1)!0.5!(2)$);
\draw[interface,thick] ($(2)!0.5!(3)$)--($(1)!0.5!(3)$);
\draw[interface,thick] ($(1)!0.5!(3)$)--($(2)!0.5!(3)$);
\node[rectangle,draw=black,fill=white,inner sep=0pt,minimum size=4pt] at ($(1)!0.5!(2)$) {};
\node[rectangle,draw=black,fill= white,inner sep=0pt,minimum size=4pt] at ($(3)!0.5!(2)$) {};
\node[rectangle,draw=black,fill=white,inner sep=0pt,minimum size=4pt] at ($(1)!0.5!(3)$) {};
\node at (0,-0.5) {red};
\end{tikzpicture}
\end{minipage}
\begin{minipage}[t]{0.23\linewidth}
\hspace*{6mm}\begin{tikzpicture}[ interface/.style={
        postaction={draw,decorate,decoration={border,angle=45,
                    amplitude=0.13cm,segment length=1mm}}},
    ]
\coordinate (1) at (-1.2,0);
\coordinate (2) at (1.2,0);
\coordinate (3) at (0,1.5);
\draw(1) --(2)--(3)--(1);
\draw[interface,thick] (2) -- (3);
\draw[interface,thick] (3) -- (1);
\fill (1) circle (2pt); 
\fill (2) circle (2pt); 
\fill (3) circle (2pt); 
\draw ($(1)!0.5!(2)$) --(3);
\node[rectangle,draw=black,fill=white,inner sep=0pt,minimum size=4pt] at ($(1)!0.5!(2)$) {};
\node at (0,-0.5) {green};
\end{tikzpicture}
\end{minipage}
\begin{minipage}[t]{0.23\linewidth}
\hspace*{6mm}\begin{tikzpicture}[ interface/.style={
        postaction={draw,decorate,decoration={border,angle=45,
                    amplitude=0.13cm,segment length=1mm}}},
    ]
\coordinate (1) at (-1.2,0);
\coordinate (2) at (1.2,0);
\coordinate (3) at (0,1.5);
\draw (1)--(2) -- (3)--(1);
\draw[interface,thick] (3) -- (1);
\fill (1) circle (2pt); 
\fill (2) circle (2pt); 
\fill (3) circle (2pt); 
\draw ($(2)!0.5!(3)$) -- ($(1)!0.5!(2)$);
\draw[interface,thick] (3)--($(1)!0.5!(2)$);
\draw[interface,thick] ($(1)!0.5!(2)$) --(2);
\node[rectangle,draw=black,fill=white,inner sep=0pt,minimum size=4pt] at ($(1)!0.5!(2)$) {};
\node[rectangle,draw=black,fill=white,inner sep=0pt,minimum size=4pt] at ($(3)!0.5!(2)$) {};
\node at (0,-0.5) {blue$_r$};
\end{tikzpicture}
\end{minipage}
\begin{minipage}[t]{0.23\linewidth}
\hspace*{6mm}\begin{tikzpicture}[ interface/.style={
        postaction={draw,decorate,decoration={border,angle=45,
                    amplitude=0.13cm,segment length=1mm}}},
    ]
\coordinate (1) at (-1.2,0);
\coordinate (2) at (1.2,0);
\coordinate (3) at (0,1.5);
\draw(1)--(2) -- (3)--(1);
\draw[interface,thick] (2) -- (3);
\fill (1) circle (2pt); 
\fill (2) circle (2pt); 
\fill (3) circle (2pt); 
\draw ($(1)!0.5!(3)$) -- ($(1)!0.5!(2)$);
\draw[interface,thick] ($(1)!0.5!(2)$)--(3);
\draw[interface,thick] (1)-- ($(1)!0.5!(2)$);
\node[rectangle,draw=black,fill=white,inner sep=0pt,minimum size=4pt] at ($(1)!0.5!(2)$) {};
\node[rectangle,draw=black,fill=white,inner sep=0pt,minimum size=4pt] at ($(1)!0.5!(3)$) {};
\node at (0,-0.5) {blue$_\ell$};
\end{tikzpicture}
\end{minipage}
\end{center}
\caption{Newest vertices per element (\revision{white squares}) displayed for each refinement pattern. \revision{The nodes that were last added in an element are the newest vertices. The set of newest nodes does not include any nodes from the initial triangulation.}}
\label{fig:newestvertex}
\end{figure}
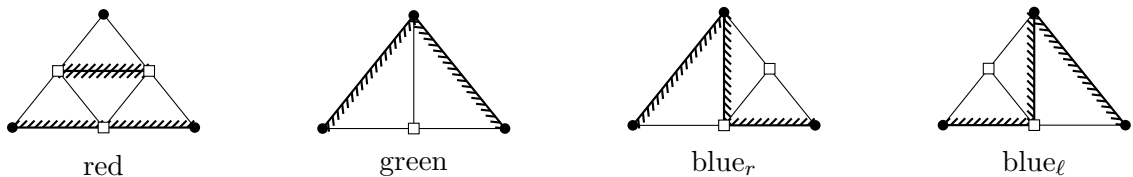

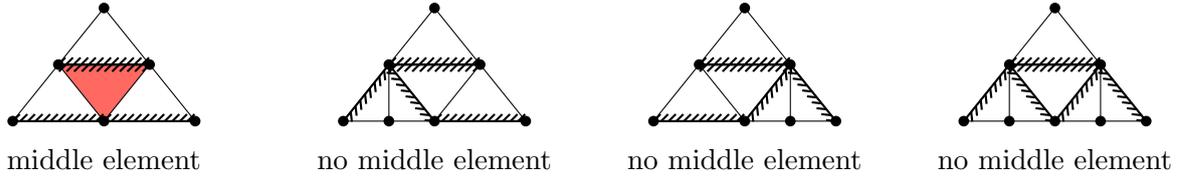
\begin{figure}[!htb]
\begin{center}
\begin{minipage}[t]{0.23\linewidth}
\hspace*{6mm}\begin{tikzpicture}[ interface/.style={
        postaction={draw,decorate,decoration={border,angle=45,
                    amplitude=0.13cm,segment length=1mm}}},
    ]
\coordinate (1) at (-1.2,0);
\coordinate (2) at (1.2,0);
\coordinate (3) at (0,1.5);
\draw (1)--(2) -- (3)--(1);
\draw [fill=pastelred]($(1)!0.5!(2)$) -- ($(2)!0.5!(3)$)--($(1)!0.5!(3)$)--($(1)!0.5!(2)$);
\draw[interface,thick] (1) --(2);
\fill (1) circle (2pt); 
\fill (2) circle (2pt); 
\fill (3) circle (2pt); 
\fill ($(1)!0.5!(2)$) circle (2pt);
\fill ($(2)!0.5!(3)$) circle (2pt);
\fill ($(1)!0.5!(3)$) circle (2pt);
\draw[interface,thick] ($(2)!0.5!(3)$)--($(1)!0.5!(3)$);
\draw[interface,thick] ($(1)!0.5!(3)$)--($(2)!0.5!(3)$);
\node at (0,-0.5) {middle element};
\end{tikzpicture}
\end{minipage}
\begin{minipage}[t]{0.23\linewidth}
\hspace*{6mm}\begin{tikzpicture}[ interface/.style={
        postaction={draw,decorate,decoration={border,angle=45,
                    amplitude=0.13cm,segment length=1mm}}},
    ]
\coordinate (1) at (-1.2,0);
\coordinate (2) at (1.2,0);
\coordinate (3) at (0,1.5);
\draw (1)--(2) -- (3)--(1);
\draw[interface,thick] ($(1)!0.5!(2)$) --(2);
\fill (1) circle (2pt); 
\fill (2) circle (2pt); 
\fill (3) circle (2pt); 
\fill ($(1)!0.5!(2)$) circle (2pt);
\fill ($(2)!0.5!(3)$) circle (2pt);
\fill ($(1)!0.5!(3)$) circle (2pt);
\fill ($(1)!0.25!(2)$) circle (2pt);
\draw($(1)!0.5!(2)$) -- ($(2)!0.5!(3)$)--($(1)!0.5!(3)$)--($(1)!0.5!(2)$);
\draw[interface,thick] ($(2)!0.5!(3)$)--($(1)!0.5!(3)$);
\draw[interface,thick] ($(1)!0.5!(3)$)--($(2)!0.5!(3)$);
\draw[interface,thick] ($(1)!0.5!(3)$)--(1);
\draw[interface,thick] ($(1)!0.5!(2)$)--($(1)!0.5!(3)$);
\draw ($(3)!0.5!(1)$) -- ($(1)!0.25!(2)$);
\node at (0,-0.5) {no middle element};
\end{tikzpicture}
\end{minipage}
\begin{minipage}[t]{0.23\linewidth}
\hspace*{6mm}\begin{tikzpicture}[ interface/.style={
        postaction={draw,decorate,decoration={border,angle=45,
                    amplitude=0.13cm,segment length=1mm}}},
    ]
\coordinate (1) at (-1.2,0);
\coordinate (2) at (1.2,0);
\coordinate (3) at (0,1.5);
\draw (1)--(2) -- (3)--(1);
\draw[interface,thick] (1) --($(1)!0.5!(2)$);
\fill (1) circle (2pt); 
\fill (2) circle (2pt); 
\fill (3) circle (2pt); 
\fill ($(1)!0.5!(2)$) circle (2pt);
\fill ($(2)!0.5!(3)$) circle (2pt);
\fill ($(1)!0.5!(3)$) circle (2pt);
\fill ($(1)!0.75!(2)$) circle (2pt);
\draw($(1)!0.5!(2)$) -- ($(2)!0.5!(3)$)--($(1)!0.5!(3)$)--($(1)!0.5!(2)$);
\draw[interface,thick] ($(2)!0.5!(3)$)--($(1)!0.5!(3)$);
\draw[interface,thick] ($(1)!0.5!(3)$)--($(2)!0.5!(3)$);
\draw[interface,thick] ($(2)!0.5!(3)$)--($(1)!0.5!(2)$);
\draw[interface,thick] (2)--($(2)!0.5!(3)$);
\draw ($(3)!0.5!(2)$) -- ($(2)!0.25!(1)$);
\node at (0,-0.5) {no middle element};
\end{tikzpicture}
\end{minipage}
\begin{minipage}[t]{0.23\linewidth}
\hspace*{6mm}\begin{tikzpicture}[ interface/.style={
        postaction={draw,decorate,decoration={border,angle=45,
                    amplitude=0.13cm,segment length=1mm}}},
    ]
\coordinate (1) at (-1.2,0);
\coordinate (2) at (1.2,0);
\coordinate (3) at (0,1.5);
\draw (1)--(2) -- (3)--(1);
\fill (1) circle (2pt); 
\fill (2) circle (2pt); 
\fill (3) circle (2pt); 
\fill ($(1)!0.5!(2)$) circle (2pt);
\fill ($(2)!0.5!(3)$) circle (2pt);
\fill ($(1)!0.5!(3)$) circle (2pt);
\fill ($(1)!0.75!(2)$) circle (2pt);
\fill ($(1)!0.25!(2)$) circle (2pt);
\draw ($(1)!0.5!(2)$) -- ($(2)!0.5!(3)$)--($(1)!0.5!(3)$)--($(1)!0.5!(2)$);
\draw[interface,thick] ($(2)!0.5!(3)$)--($(1)!0.5!(3)$);
\draw[interface,thick] ($(1)!0.5!(3)$)--($(2)!0.5!(3)$);
\draw ($(3)!0.5!(1)$) -- ($(1)!0.25!(2)$);
\draw ($(3)!0.5!(2)$) -- ($(2)!0.25!(1)$);
\draw[interface,thick] ($(1)!0.5!(3)$)--(1);
\draw[interface,thick] ($(1)!0.5!(2)$)--($(1)!0.5!(3)$);
\draw[interface,thick] ($(2)!0.5!(3)$)--($(1)!0.5!(2)$);
\draw[interface,thick] (2)--($(2)!0.5!(3)$);
\node at (0,-0.5) {no middle element};
\end{tikzpicture}
\end{minipage}
\caption{Red middle element is \revision{painted in color}. Middle elements which neighbours have a different refinement level are not considered in the set of red middle elements. \revision{To determine a red middle element, for each element all three neighbours that share an edge with this element are determined and the location of the reference edge is compared. If it matches the refinement pattern on the left, it is a red middle element. If the surrounding leads to other combinations as shown, e.\,g., in the other three patterns, it is not a red middle element.}}
\label{fig:redmiddleelement}\end{center}
\end{figure}

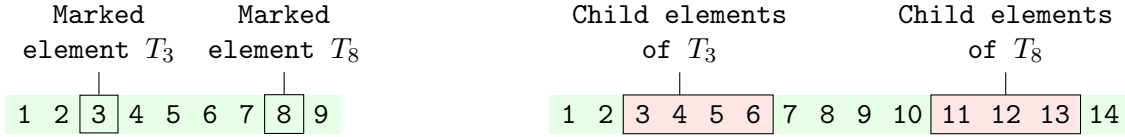
\begin{figure}[h!]
\begin{minipage}{0.38\textwidth}
\centering
\begin{tikzpicture}[font=\ttfamily,
array/.style={matrix of nodes,nodes={draw, minimum size=7mm, fill=green!30},column sep=-\pgflinewidth, row sep=0.5mm, nodes in empty cells,
row 1/.style={nodes={draw=none, fill=none, minimum size=5mm}},
row 1 column 3/.style={nodes={draw}}}]

\matrix[array] (array) {
1 & 2 & 3 & 4 & 5 & 6 & 7 & 8 & 9 \\};

\begin{scope}[on background layer]
\fill[green!10] (array-1-1.north west) rectangle (array-1-9.south east);
\end{scope}
\draw (array-1-3.north)--++(90:3mm) node [above,align=center, anchor=south] (first) {Marked \\ element $T_3$};
\draw (array-1-8.north)--++(90:3mm) node [above,align=center, anchor=south] (first) {Marked \\ element $T_8$};
\draw (array-1-8.north west) rectangle (array-1-8.south east);
\end{tikzpicture}
\end{minipage}\hfill
\begin{minipage}{0.6\textwidth}
\centering
\begin{tikzpicture}[font=\ttfamily,
array/.style={matrix of nodes,nodes={draw, minimum size=7mm, fill=green!30},column sep=-\pgflinewidth, row sep=0.5mm, nodes in empty cells,
row 1/.style={nodes={draw=none, fill=none, minimum size=5mm}},
row 0 column 0/.style={nodes={draw}}}]
\matrix[array] (array) {
1 & 2 & 3 & 4 & 5 & 6 & 7 & 8 & 9 & 10 & 11 & 12 & 13 & 14\\
};
\begin{scope}[on background layer]
\fill[green!10] (array-1-1.north west) rectangle (array-1-2.south east);
\end{scope}
\begin{scope}[on background layer]
\fill[red!10] (array-1-3.north west) rectangle (array-1-6.south east);
\end{scope}
\begin{scope}[on background layer]
\fill[green!10] (array-1-7.north west) rectangle (array-1-10.south east);
\end{scope}
\begin{scope}[on background layer]
\fill[red!10] (array-1-11.north west) rectangle (array-1-13.south east);
\end{scope}
\begin{scope}[on background layer]
\fill[green!10] (array-1-14.north west) rectangle (array-1-14.south east);
\end{scope}
\draw (array-1-4.north)--++(90:3mm) node [above,align=center, anchor=south] (first) {Child elements \\of $T_3$};
\draw (array-1-12.north)--++(90:3mm) node [above,align=center, anchor=south] (first) {Child elements \\of $T_8$};
\draw (array-1-3.north west) rectangle (array-1-6.south east);
\draw (array-1-11.north west) rectangle (array-1-13.south east);
\end{tikzpicture}
\end{minipage}
\caption{Numbering of elements before (left) and after refinement (right). \revision{Subtriangles of an element are stored at the previous position of this element and successive positions, rather than appending the new elements at the end of the array. The position of other elements is then shifted by the number of newly created elements.}}
\label{fig:enumeration}
\end{figure}

In this work, the implementation of the proposed coarsening algorithm is built on the RGB refinement implementation in the \texttt{ameshref}-package  \cite{web,funkenschmidt}. This is why we focus on this concrete data structure. However, the proposed coarsening algorithm can also be based on other RGB refinement implementations without explicit refinement history. For this to work, the following must be ensured: 
\revision{\begin{axioms}
\item \label{itm:1} Reference edges must be incorporated in the data structure. 
\item \label{itm:2} The data structure needs to be designed such that newest vertices, cf.~Figure~\ref{fig:newestvertex}, can easily \hspace*{1em} be determined. 
\item \label{itm:3} Elements have to be numbered in a way that a blue refinement leads to the same numbering \hspace*{1em} as an application of two green refinements. 
\item \label{itm:4} The middle element of red refinement patterns, cf.~Figure~\ref{fig:redmiddleelement}, must be identifiable.
\item \label{itm:5} Child elements are to be stored consecutively at the former position of the father element,  \hspace*{1em} cf.~Figure~\ref{fig:TrefineRGBpattern} and Figure~\ref{fig:enumeration}.
\end{axioms}}

\revision{Let us examine each of these listed points in more detail. \ref{itm:1} Reference edges play a crucial role in RGB refinement and are thus also important for coarsening. We know how the reference edges are chosen during the refinement process. Thus, this is one key information in joining elements back together and determining the reference edge of the father element, cf. the pattern \emph{none} in Figure~\ref{fig:TrefineRGBpattern}. However, this information on its own is not sufficient. \ref{itm:2} In each refinement step, new nodes are added. To this end, the newest nodes are the nodes that are first removed in a coarsening step. It is therefore important information as it provides the node candidates for removal, cf. Figure~\ref{fig:newestvertex}. These node candidates do not give any information about the patterns that lie around these nodes. However, they are important when it comes to joining elements back together. Our algorithm distinguishes between a red and a green pattern. \ref{itm:3} As a blue pattern is created by a green refinement of a green-refined element, we can deal with a blue pattern via a two-step removal of green patterns. For this to work, it is required that the numbering of a blue pattern leads to the same numbering as the application of two green refinements. Thus, it only remains to distinguish red and green patterns from each other. Where a green refinement is the result of a bisection of the element, a red refinement creates four subtriangles by joining its midpoints together. \ref{itm:4} This means that there is one triangle with new nodes only, that we call a red middle element, cf.~Figure~\ref{fig:redmiddleelement}. Such a middle element does not exist for green or blue patterns and therefore distinguishes red patterns uniquely from green or blue patterns. \ref{itm:5} To determine the color of the pattern, we make use of the property that subelements of one and the same element are stored consecutively, cf.~Figure~\ref{fig:TrefineRGBpattern}. This ensures that elements that need to be joined together when coarsening are implicitly given in the data structure. So far, we have come to a point where we determined node candidates for removal and know their surrounding patterns. With this information, elements can be coarsened once. To join elements together in a subsequent coarsening step, it is necessary that the local information of former joined elements can still be reconstructed. We can ensure this by storing the subtriangles of a refined element at the previous position of the father element and successive positions, rather than appending the new elements at the end of an array. The position of other triangles is then shifted by the number of newly created elements, cf. Figure~\ref{fig:enumeration}. If we reverse this operation when coarsening, we make sure that elements with the same father element are stored consecutively after a coarsening step. Therefore, our coarsening routine is able to coarsen back to the initial triangulation and no explicit history tree is needed to invert the refinement.\\
In the following subsection, we introduce our MATLAB implementation of the RGB refinement. The fact that the requirements \ref{itm:1}, \ref{itm:3} and \ref{itm:5} are met is made clear by the following explanation of the data structure used as well as the RGB refinement and the corresponding storage of new coordinates and elements. How to meet the requirements \ref{itm:2} and \ref{itm:4} is discussed in Section~\ref{sect:Matlab}.}

\subsection{\revision{Our} MATLAB Implementation of RGB Refinement}\label{sect:RGBrefine}

In this section, we give some insights into the implementation of the RGB refinement in the \texttt{ameshref}-package \cite{web,funkenschmidt,ameshref}. We focus on the parts that are essential for our coarsening routine. For a more thorough explanation, we refer to \cite{funkenschmidt,ameshref}.
We represent a triangulation $\mathcal{T}=\left\{T_1,\ldots,T_M\right\}$ with nodes $\mathcal{N}=\left\{v_1,\ldots,v_N\right\}$ as follows: The $x$- and $y$-coordinates of the nodes $\mathcal{N}$ are stored within an $N \times 2$ array \texttt{coordinates}. Furthermore, we represent the element-connectivity within an $M\times 3$ array \texttt{elements} where one row stores the indices of the element's three vertices $v_i,v_j,v_k \in \mathcal{N}$ with $i,j,k \in \left\{1,\ldots,N\right\}$. Optionally, boundary edges can be stored in an additional array with indices of the edge's two vertices. As depicted in Figure~\ref{fig:TrefineRGBpattern}, reference edges play a crucial role. Instead of storing this information in an additional data structure, we capture the reference edge implicitly as the edge between the first two vertices of an element indexed by the first two entries in the array \texttt{elements} \revision{(\ref{itm:1})}. Elements are numbered counterclockwise.

\begin{figure}[!htb]
\centering
\begin{minipage}{0.325\textwidth}
\begin{tikzpicture}[ interface/.style={
        postaction={draw,decorate,decoration={border,angle=45,
                    amplitude=0.1cm,segment length=1mm}}},
    ]
\coordinate (1) at (0,0);
\coordinate (2) at (2,0);
\coordinate (3) at (2,2);
\coordinate (4) at (0,2);
\fill[lightblue] (1)--(3)--(4);
\fill[lightgreen] (1)--(2)--(3);
\draw[interface,thick] (1)--(2)--(3)--(4); 
\draw (4)--($(4)!0.5!(1)$);   
\draw (1)--($(4)!0.5!(1)$)--($(3)!0.5!(1)$)--(2);
\draw ($(3)!0.5!(1)$)--(3);
\draw[interface,thick] ($(3)!0.5!(1)$)--(4);
\draw[interface,thick] (1)--($(3)!0.5!(1)$);
\node at (1,0.5) {g5};
\node at (1.5,1) {g4};
\node at (0.25,1.25) {b1};
\node at (0.25,0.75) {b2};
\node at (1,1.5) {b3};
\fill (1) circle (2pt);
\fill (2) circle (2pt);
\fill (3) circle (2pt);
\fill (4) circle (2pt);
\fill ($(3)!0.5!(1)$) circle (2pt);
\fill ($(4)!0.5!(1)$) circle (2pt);
\draw[<->,thick] (2.25,1)--(2.75,1);
\coordinate (1) at (3,0);
\coordinate (2) at (5,0);
\coordinate (3) at (5,2);
\coordinate (4) at (3,2);
\fill[lightgreen] (1)--($(3)!0.5!(1)$)--(4);
\fill[lightgreen] (4)--(2)--(3);
\draw[interface,thick] (1)--(2)--(3)--(4); 
\draw (4)--($(4)!0.5!(1)$);   
\draw (1)--($(4)!0.5!(1)$)--($(3)!0.5!(1)$)--(2);
\draw ($(3)!0.5!(1)$)--(3);
\draw[interface,thick] ($(3)!0.5!(1)$)--(4);
\draw[interface,thick] (1)--($(3)!0.5!(1)$);
\node at (4,0.5) {5};
\node at (4.5,1) {g4};
\node at (3.25,1.25) {g1};
\node at (3.25,0.75) {g2};
\node at (4,1.5) {g3};
\fill (1) circle (2pt);
\fill (2) circle (2pt);
\fill (3) circle (2pt);
\fill (4) circle (2pt);
\fill ($(3)!0.5!(1)$) circle (2pt);
\fill ($(4)!0.5!(1)$) circle (2pt);
\end{tikzpicture}
\end{minipage}
\begin{minipage}{0.325\textwidth}
\begin{tikzpicture}[ interface/.style={
        postaction={draw,decorate,decoration={border,angle=45,
                    amplitude=0.1cm,segment length=1mm}}},
    ]
\coordinate (1) at (0,0);
\coordinate (2) at (2,0);
\coordinate (3) at (2,2);
\coordinate (4) at (0,2);
\fill[lightblue] (1)--(2)--(3)--(4);
\draw[interface,thick] (2)--(3); 
\draw[interface,thick] (4)--(1); 
\draw (1)--(2);  
\draw (3)--(4);   
\draw ($(3)!0.5!(4)$)--($(2)!0.5!(1)$);    
\draw[interface,thick] (2)--($(3)!0.5!(1)$)--(1);
\draw[interface,thick] (4)--($(3)!0.5!(1)$)--(3);
\node at (0.5,1) {b1};
\node at (0.75,1.75) {b3};
\node at (1.25,1.75) {b2};
\node at (1.5,1) {b4};
\node at (0.75,0.25) {b5};
\node at (1.25,0.25) {b6};
\fill (1) circle (2pt);
\fill (2) circle (2pt);
\fill (3) circle (2pt);
\fill (4) circle (2pt);
\fill ($(2)!0.5!(1)$) circle (2pt);
\fill ($(4)!0.5!(3)$) circle (2pt);
\fill ($(1)!0.5!(3)$) circle (2pt);
\draw[<->,thick] (2.25,1)--(2.75,1);
\coordinate (1) at (3,0);
\coordinate (2) at (5,0);
\coordinate (3) at (5,2);
\coordinate (4) at (3,2);
\fill[lightblue] (2)--(3)--(4);
\fill[lightgreen] (1)--(2)--($(1)!0.5!(3)$);
\draw[interface,thick] (2)--(3); 
\draw[interface,thick] (4)--(1); 
\draw (1)--(2);  
\draw (3)--(4);   
\draw ($(3)!0.5!(4)$)--($(2)!0.5!(1)$);    
\draw[interface,thick] (2)--($(3)!0.5!(1)$)--(1);
\draw[interface,thick] (4)--($(3)!0.5!(1)$)--(3);
\node at (3.5,1) {1};
\node at (3.75,1.75) {b3};
\node at (4.25,1.75) {b2};
\node at (4.5,1) {b4};
\node at (3.75,0.25) {g5};
\node at (4.25,0.25) {g6};
\fill (1) circle (2pt);
\fill (2) circle (2pt);
\fill (3) circle (2pt);
\fill (4) circle (2pt);
\fill ($(2)!0.5!(1)$) circle (2pt);
\fill ($(4)!0.5!(3)$) circle (2pt);
\fill ($(1)!0.5!(3)$) circle (2pt);
\end{tikzpicture}
\end{minipage}
\begin{minipage}{0.325\textwidth}
\begin{tikzpicture}[ interface/.style={
        postaction={draw,decorate,decoration={border,angle=45,
                    amplitude=0.1cm,segment length=1mm}}},
    ]
\coordinate (1) at (0,0);
\coordinate (2) at (2,0);
\coordinate (3) at (2,2);
\coordinate (4) at (0,2);
\fill[lightgreen] (1)--(3)--(4);
\fill[lightgreen] (1)--(2)--(3);
\draw[interface,thick] (1)--(2)--(3)--(4)--(1); 
\draw ($(3)!0.5!(1)$)--(2);
\draw ($(3)!0.5!(1)$)--(3);
\draw ($(3)!0.5!(1)$)--(4);
\draw (1)--($(3)!0.5!(1)$);
\node at (1,0.5) {g4};
\node at (1.5,1) {g3};
\node at (0.5,1) {g1};
\node at (1,1.5) {g2};
\fill (1) circle (2pt);
\fill (2) circle (2pt);
\fill (3) circle (2pt);
\fill (4) circle (2pt);
\fill ($(3)!0.5!(1)$) circle (2pt);
\draw[<->,thick] (2.25,1)--(2.75,1);
\coordinate (1) at (3,0);
\coordinate (2) at (5,0);
\coordinate (3) at (5,2);
\coordinate (4) at (3,2);
\fill[lightgreen] (2)--(3)--(4);
\draw[interface,thick] (1)--(2)--(3)--(4)--(1); 
\draw ($(3)!0.5!(1)$)--(2);
\draw ($(3)!0.5!(1)$)--(3);
\draw ($(3)!0.5!(1)$)--(4);
\draw (1)--($(3)!0.5!(1)$);
\node at (4,0.5) {4};
\node at (4.5,1) {g3};
\node at (3.5,1) {1};
\node at (4,1.5) {g2};
\fill (1) circle (2pt);
\fill (2) circle (2pt);
\fill (3) circle (2pt);
\fill (4) circle (2pt);
\fill ($(3)!0.5!(1)$) circle (2pt);
\end{tikzpicture}
\end{minipage}

\caption{Three examples to show that adjacent blue\revision{(b)} and green\revision{(g)} patterns cannot be distinguished on the basis of numbering: [Actual connection of elements] $\leftrightarrow$ [Another conceivable connection based on this numbering].}
\label{fig:nodisting}
\end{figure}
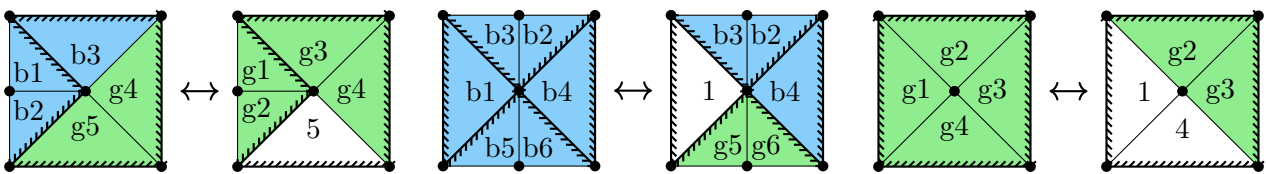

In adaptive procedures, a set of marked elements $\mathcal{M}$ is given. We flag elements $T \in \mathcal{M}$ by marking each edge of the element for bisection. Obviously, neighbouring elements $T \not\in \mathcal{M}$ are affected indirectly by this marking. A CLOSURE step is performed to avoid creating hanging nodes. As mentioned, the assignment of reference edges ensures the shape regularity of the triangulation. To this end, the reference edge needs to be bisected before any other edge of this element is refined. For this reason, we mark edges according to the hash map shown in Table~\ref{tab:mappingRGB} and loop through this CLOSURE step until no further markers are added. Then, we refine the elements according to Figure~\ref{fig:TrefineRGBpattern} \revision{(\ref{itm:3})} and save the new elements at the corresponding position in the array as depicted exemplarily in Figure~\ref{fig:enumeration} \revision{(\ref{itm:5})}. This is essential for coarsening of more than just one layer as the recursive information is implicitly given in the array \texttt{elements}. Without this, any hierarchical information is lost and thus, coarsening of more than one layer is impossible. As direct consequence of this way of storing the refined elements, adjacent blue and green patterns can no longer be distinguished from each other, cf. Figure~\ref{fig:nodisting}. \revision{It thus makes sense to consider inverting green and red refinement only. }
{
\centering
\begin{table}[!htb]
     \caption[Mapping of eight possible markings to the five patterns allowed in RGB refinement.]{Mapping of eight possible markings to the five patterns allowed in RGB refinement. For each hash a binary number is given. }
       \label{tab:mappingRGB}
\begin{tabular}{ccccccccc}\hline
&&&&&&&\\
mark& \begin{minipage}{0.085\textwidth}
\begin{tikzpicture}
\vspace*{2mm}
\coordinate (1) at (0,0);
\coordinate (2) at (1.2,0);
\coordinate (3) at (0.6,0.8);
\draw(1) --(2)--(3)--(1);
\fill (1) circle (2pt); 
\fill (2) circle (2pt); 
\fill (3) circle (2pt); 
\end{tikzpicture}
\end{minipage} &   \begin{minipage}{0.085\textwidth}
\begin{tikzpicture}
\vspace*{2mm}
\coordinate (1) at (0,0);
\coordinate (2) at (1.2,0);
\coordinate (3) at (0.6,0.8);
\draw(1) --(2)--(3)--(1);
\draw ($(1)!0.5!(2)$) node[cross,very thick] {};
\fill (1) circle (2pt); 
\fill (2) circle (2pt); 
\fill (3) circle (2pt); 
\end{tikzpicture}
\end{minipage} &   \begin{minipage}{0.085\textwidth}
\begin{tikzpicture}
\vspace*{2mm}
\coordinate (1) at (0,0);
\coordinate (2) at (1.2,0);
\coordinate (3) at (0.6,0.8);
\draw(1) --(2)--(3)--(1);
\draw ($(2)!0.5!(3)$) node[cross,very thick, rotate = 135] {};
\fill (1) circle (2pt); 
\fill (2) circle (2pt); 
\fill (3) circle (2pt); 
\end{tikzpicture}
\end{minipage} &   \begin{minipage}{0.085\textwidth}
\begin{tikzpicture}
\vspace*{2mm}
\coordinate (1) at (0,0);
\coordinate (2) at (1.2,0);
\coordinate (3) at (0.6,0.8);
\draw(1) --(2)--(3)--(1);
\draw ($(1)!0.5!(2)$) node[cross,very thick] {};
\draw ($(2)!0.5!(3)$) node[cross,very thick, rotate = 135] {};
\fill (1) circle (2pt); 
\fill (2) circle (2pt); 
\fill (3) circle (2pt); 
\end{tikzpicture}
\end{minipage} &   \begin{minipage}{0.085\textwidth}
\begin{tikzpicture}
\vspace*{2mm}
\coordinate (1) at (0,0);
\coordinate (2) at (1.2,0);
\coordinate (3) at (0.6,0.8);
\draw(1) --(2)--(3)--(1);
\draw ($(3)!0.5!(1)$) node[cross,very thick, rotate = 45] {};
\fill (1) circle (2pt); 
\fill (2) circle (2pt); 
\fill (3) circle (2pt); 
\end{tikzpicture}
\end{minipage} &  \begin{minipage}{0.085\textwidth}
\begin{tikzpicture}
\vspace*{2mm}
\coordinate (1) at (0,0);
\coordinate (2) at (1.2,0);
\coordinate (3) at (0.6,0.8);
\draw(1) --(2)--(3)--(1);
\draw ($(1)!0.5!(2)$) node[cross,very thick] {};
\draw ($(3)!0.5!(1)$) node[cross,very thick, rotate = 45] {};
\fill (1) circle (2pt); 
\fill (2) circle (2pt); 
\fill (3) circle (2pt); 
\end{tikzpicture}
\end{minipage}
& 
  \begin{minipage}{0.085\textwidth}
\begin{tikzpicture}
\vspace*{2mm}
\coordinate (1) at (0,0);
\coordinate (2) at (1.2,0);
\coordinate (3) at (0.6,0.8);
\draw(1) --(2)--(3)--(1);
\draw ($(2)!0.5!(3)$) node[cross,very thick, rotate = 135] {};
\draw ($(3)!0.5!(1)$) node[cross,very thick, rotate = 45] {};
\fill (1) circle (2pt); 
\fill (2) circle (2pt); 
\fill (3) circle (2pt); 
\end{tikzpicture}
\end{minipage} &  
  \begin{minipage}{0.085\textwidth}
\begin{tikzpicture}
\vspace*{2mm}
\coordinate (1) at (0,0);
\coordinate (2) at (1.2,0);
\coordinate (3) at (0.6,0.8);
\draw(1) --(2)--(3)--(1);
\draw ($(1)!0.5!(2)$) node[cross,very thick] {};
\draw ($(2)!0.5!(3)$) node[cross,very thick, rotate = 135] {};
\draw ($(3)!0.5!(1)$) node[cross,very thick, rotate = 45] {};
\fill (1) circle (2pt); 
\fill (2) circle (2pt); 
\fill (3) circle (2pt); 
\end{tikzpicture}
\end{minipage} \\
&&&&&&&&\\ 
bin &000 & 100 & 010 & 110 & 001 & 101 & 011 & 111\\ 
&\begin{tikzpicture}
\draw[->] (0,0)--(0,-0.5);
\end{tikzpicture}&\begin{tikzpicture}
\draw[->] (0,0)--(0,-0.5);
\end{tikzpicture}&\begin{tikzpicture}
\draw[->] (0,0)--(0,-0.5);
\end{tikzpicture}&\begin{tikzpicture}
\draw[->] (0,0)--(0,-0.5);
\end{tikzpicture}&\begin{tikzpicture}
\draw[->] (0,0)--(0,-0.5);
\end{tikzpicture}&\begin{tikzpicture}
\draw[->] (0,0)--(0,-0.5);
\end{tikzpicture}&\begin{tikzpicture}
\draw[->] (0,0)--(0,-0.5);
\end{tikzpicture}&\begin{tikzpicture}
\draw[->] (0,0)--(0,-0.5);
\end{tikzpicture}\\
&&&&&&&&\\
hash& \begin{minipage}{0.085\textwidth}
\begin{tikzpicture}
\vspace*{2mm}
\coordinate (1) at (0,0);
\coordinate (2) at (1.2,0);
\coordinate (3) at (0.6,0.8);
\draw(1) --(2)--(3)--(1);
\fill (1) circle (2pt); 
\fill (2) circle (2pt); 
\fill (3) circle (2pt); 
\end{tikzpicture}
\end{minipage} &   \begin{minipage}{0.085\textwidth}
\begin{tikzpicture}
\vspace*{2mm}
\coordinate (1) at (0,0);
\coordinate (2) at (1.2,0);
\coordinate (3) at (0.6,0.8);
\coordinate (4) at ($(1)!0.5!(2)$);
\draw(1) --(2)--(3)--(1);
\draw (4)--(3);
\fill (1) circle (2pt); 
\fill (2) circle (2pt); 
\fill (3) circle (2pt); 
\fill (4) circle (2pt); 
\end{tikzpicture}
\end{minipage} &   \begin{minipage}{0.085\textwidth}
\begin{tikzpicture}
\vspace*{2mm}
\coordinate (1) at (0,0);
\coordinate (2) at (1.2,0);
\coordinate (3) at (0.6,0.8);
\coordinate (4) at ($(1)!0.5!(2)$);
\coordinate (5) at ($(2)!0.5!(3)$);
\draw(1) --(2)--(3)--(1);
\draw (3)--(4)--(5);
\fill (1) circle (2pt); 
\fill (2) circle (2pt); 
\fill (3) circle (2pt); 
\fill (4) circle (2pt); 
\fill (5) circle (2pt); 
\end{tikzpicture}
\end{minipage} &   \begin{minipage}{0.085\textwidth}
\begin{tikzpicture}
\vspace*{2mm}
\coordinate (1) at (0,0);
\coordinate (2) at (1.2,0);
\coordinate (3) at (0.6,0.8);
\coordinate (4) at ($(1)!0.5!(2)$);
\coordinate (5) at ($(2)!0.5!(3)$);
\draw(1) --(2)--(3)--(1);
\draw (3)--(4)--(5);
\fill (1) circle (2pt); 
\fill (2) circle (2pt); 
\fill (3) circle (2pt); 
\fill (4) circle (2pt); 
\fill (5) circle (2pt); 
\end{tikzpicture}
\end{minipage} &   \begin{minipage}{0.085\textwidth}
\begin{tikzpicture}
\vspace*{2mm}
\coordinate (1) at (0,0);
\coordinate (2) at (1.2,0);
\coordinate (3) at (0.6,0.8);
\coordinate (4) at ($(1)!0.5!(2)$);
\coordinate (5) at ($(1)!0.5!(3)$);
\draw(1) --(2)--(3)--(1);
\draw (3)--(4)--(5);
\fill (1) circle (2pt); 
\fill (2) circle (2pt); 
\fill (3) circle (2pt); 
\fill (4) circle (2pt); 
\fill (5) circle (2pt); 
\end{tikzpicture}
\end{minipage} &  \begin{minipage}{0.085\textwidth}
\begin{tikzpicture}
\vspace*{2mm}
\coordinate (1) at (0,0);
\coordinate (2) at (1.2,0);
\coordinate (3) at (0.6,0.8);
\coordinate (4) at ($(1)!0.5!(2)$);
\coordinate (5) at ($(1)!0.5!(3)$);
\draw(1) --(2)--(3)--(1);
\draw (3)--(4)--(5);
\fill (1) circle (2pt); 
\fill (2) circle (2pt); 
\fill (3) circle (2pt); 
\fill (4) circle (2pt); 
\fill (5) circle (2pt); 
\end{tikzpicture}
\end{minipage}
& 
  \begin{minipage}{0.085\textwidth}
\begin{tikzpicture}
\vspace*{2mm}
\coordinate (1) at (0,0);
\coordinate (2) at (1.2,0);
\coordinate (3) at (0.6,0.8);
\coordinate (4) at ($(1)!0.5!(2)$);
\coordinate (5) at ($(2)!0.5!(3)$);
\coordinate (6) at ($(1)!0.5!(3)$);
\draw(1) --(2)--(3)--(1);
\draw (6)--(4)--(5)--(6);
\fill (1) circle (2pt); 
\fill (2) circle (2pt); 
\fill (3) circle (2pt); 
\fill (4) circle (2pt); 
\fill (5) circle (2pt); 
\fill (6) circle (2pt);
\end{tikzpicture}
\end{minipage} &  
  \begin{minipage}{0.085\textwidth}
\begin{tikzpicture}
\vspace*{2mm}
\coordinate (1) at (0,0);
\coordinate (2) at (1.2,0);
\coordinate (3) at (0.6,0.8);
\coordinate (4) at ($(1)!0.5!(2)$);
\coordinate (5) at ($(2)!0.5!(3)$);
\coordinate (6) at ($(1)!0.5!(3)$);
\draw(1) --(2)--(3)--(1);
\draw (6)--(4)--(5)--(6);
\fill (1) circle (2pt); 
\fill (2) circle (2pt); 
\fill (3) circle (2pt); 
\fill (4) circle (2pt); 
\fill (5) circle (2pt); 
\fill (6) circle (2pt);
\end{tikzpicture}
\end{minipage} \\&&&&&&&&\\
bin &000 & 100 & 110 & 110 & 101 & 101 & 111 & 111\\
type & none & green & blue$_r$ & blue$_r$ & blue$_\ell$ & blue$_\ell$ &red & red\\
&&&&&&&&\\ \hline
    \end{tabular}

  \end{table}}
  
\revision{Figure~\ref{fig:example_mesh} serves to illustrate the data structure used in the \texttt{ameshref}-package as well as the RGB refinement and the corresponding storage of the new coordinates and elements.}

\revision{
\begin{figure}[!htb]
\centering
\scalebox{0.7}{
\begin{minipage}{\textwidth}
 \begin{tikzpicture}
 \draw[thin] (0,0)--(18.5,0);
 \node[right] at (0,0.25) {{\bfseries(A)} \small Initial triangulation with reference edges displayed as hatched lines.};
  \draw[thin] (0,0.5)--(18.5,0.5);
 \end{tikzpicture}
 \begin{minipage}{\textwidth}
   \vspace*{0.25cm}
 \begin{minipage}{.50\textwidth}
 \scalebox{1.7}{
\begin{tikzpicture}[ interface/.style={
        postaction={draw,decorate,decoration={border,angle=45,
                    amplitude=0.08cm,segment length=1mm}}},
    ]
\draw[style=help lines,step=0.5cm,gray, very thin] (0,0) grid (2,2);
\draw[->] (0,0) -- (2.4,0); \draw[->] (0,0) -- (0,2.4);
\foreach \x in {-1,0,1}
\node[below] (1) at (0,0) {{\tiny 1}};
\node[below] (2) at (2,0) {{\tiny 2}};
\node[right] (3) at (2,2) {{\tiny 3}};
\node[left] (4) at (0,2) {{\tiny 4}};
\coordinate (1) at (0,0);
\coordinate (2) at (2,0);
\coordinate (3) at (2,2);
\coordinate (4) at (0,2);
\draw[thick] (1)--(2)--(3)--(4)--(1);
\draw[thick,interface] (1)--(3);
\draw[thick,interface] (3)--(1);
\fill (1) circle (1.5pt);
\fill (2) circle (1.5pt);
\fill (3) circle (1.5pt);
\fill (4) circle (1.5pt);
\end{tikzpicture}}
 \end{minipage}
  \begin{minipage}{.15\textwidth}
  \lstinline[basicstyle=\small\lstbasicfont]|coordinates|
\begin{lstlisting}[basicstyle=\small\lstbasicfont,aboveskip=2mm,belowskip=7.35mm,frame=lines,nolol=true]
   0   0
   2   0
   2   2
   0   2
\end{lstlisting}
  \end{minipage}\hspace*{2cm}
  \begin{minipage}{0.15\textwidth}
  \vspace*{-1.2cm}
    \lstinline[basicstyle=\small\lstbasicfont]|elements|
\begin{lstlisting}[basicstyle=\small\lstbasicfont,aboveskip=2mm,frame=lines,nolol=true]
   1   3   4
   3   1   2
\end{lstlisting}
 \end{minipage}
   \vspace*{0.25cm}
 \end{minipage}\\
 \begin{tikzpicture}
 \draw[thin] (0,0)--(18.5,0);
 \node[right] at (0,0.25) {{\bfseries(B)} \small Refined mesh obtained by marking both elements in the initial triangulation and performing an RGB refinement.};
  \draw[thin] (0,0.5)--(18.5,0.5);
 \end{tikzpicture}
  \begin{minipage}{\textwidth}
    \vspace*{0.25cm}
  \begin{minipage}{.50\textwidth}
 \scalebox{1.7}{
\begin{tikzpicture}[ interface/.style={
        postaction={draw,decorate,decoration={border,angle=45,
                    amplitude=0.08cm,segment length=1mm}}},
    ]
\draw[style=help lines,step=0.5cm,gray, very thin] (0,0) grid (2,2);
\draw[->] (0,0) -- (2.4,0); \draw[->] (0,0) -- (0,2.4);
\foreach \x in {-1,0,1}
\node[below] (1) at (0,0) {{\tiny 1}};
\node[below] (2) at (2,0) {{\tiny 2}};
\node[right] (3) at (2,2) {{\tiny 3}};
\node[left] (4) at (0,2) {{\tiny 4}};
\node[below] (5) at (1,0) {{\tiny 5}};
\node[above left] (6) at (1,1) {{\tiny 6}};
\node[left] (7) at (0,1) {{\tiny 7}};
\node[right] (8) at (2,1) {{\tiny 8}};
\node[above] (9) at (1,2) {{\tiny 9}};
\coordinate (1) at (0,0);
\coordinate (2) at (2,0);
\coordinate (3) at (2,2);
\coordinate (4) at (0,2);
\coordinate (5) at (1,0);
\coordinate (6) at (1,1);
\coordinate (7) at (0,1);
\coordinate (8) at (2,1);
\coordinate (9) at (1,2);
\draw[thick] (1)--(2)--(3)--(4)--(1);
\draw[thick] (6)--(9)--(7)--(6)--(5)--(8)--(6);
\draw[thick,interface] (1)--(3);
\draw[thick,interface] (3)--(1);
\draw[thick,interface] (9)--(7);
\draw[thick,interface] (7)--(9);
\draw[thick,interface] (8)--(5);
\draw[thick,interface] (5)--(8);
\fill (1) circle (1.5pt);
\fill (2) circle (1.5pt);
\fill (3) circle (1.5pt);
\fill (4) circle (1.5pt);
\fill (5) circle (1.5pt);
\fill (6) circle (1.5pt);
\fill (7) circle (1.5pt);
\fill (8) circle (1.5pt);
\fill (9) circle (1.5pt);
\end{tikzpicture}}
 \end{minipage}
  \begin{minipage}{.15\textwidth}
  \lstinline[basicstyle=\small\lstbasicfont]|coordinates|
\begin{lstlisting}[basicstyle=\small\lstbasicfont,aboveskip=2mm,belowskip=7.35mm,frame=lines,nolol=true]
   0   0
   2   0
   2   2
   0   2
   1   0
   1   1
   0   1
   2   1
   1   2
\end{lstlisting}
  \end{minipage}\hspace*{2cm}
  \begin{minipage}{0.15\textwidth}
  \vspace*{-0.85cm}
  \lstinline[basicstyle=\small\lstbasicfont]|elements|
\begin{lstlisting}[basicstyle=\small\lstbasicfont,aboveskip=2mm,frame=lines,nolol=true]
   1   6   7
   6   3   9
   7   9   4
   9   7   6
   3   6   8
   6   1   5
   8   5   2
   5   8   6
\end{lstlisting}
 \end{minipage}
  \vspace*{0.25cm}
 \end{minipage}\\
 \begin{tikzpicture}
 \draw[thin] (0,0)--(18.5,0);
 \node[right] at (0,0.25) {{\bfseries(C)} \small Adaptive triangulation obtained by RGB refining the mesh in (B) for the marked element number 8.};
  \draw[thin] (0,0.5)--(18.5,0.5);
 \end{tikzpicture}
  \begin{minipage}{\textwidth}
    \vspace*{0.25cm}
  \begin{minipage}{.50\textwidth}
    \vspace*{-2.3cm}
 \scalebox{1.7}{
\begin{tikzpicture}[ interface/.style={
        postaction={draw,decorate,decoration={border,angle=45,
                    amplitude=0.08cm,segment length=1mm}}},
    ]
\draw[style=help lines,step=0.5cm,gray, very thin] (0,0) grid (2,2);
\draw[->] (0,0) -- (2.4,0); \draw[->] (0,0) -- (0,2.4);
\foreach \x in {-1,0,1}
\node[below] (1) at (0,0) {{\tiny 1}};
\node[below] (2) at (2,0) {{\tiny 2}};
\node[right] (3) at (2,2) {{\tiny 3}};
\node[left] (4) at (0,2) {{\tiny 4}};
\node[left] (10) at (0.5, 0.5) {{\tiny 10}};
\node[above] (11) at (1.5,1.5) {{\tiny 11}};
\node (12) at (0.85,0.65) {{\tiny 12}};
\node[below] (13) at (1.5,0.5) {{\tiny 13}};
\node (14) at (1.65,1.15) {{\tiny 14}};
\coordinate (1) at (0,0);
\coordinate (2) at (2,0);
\coordinate (3) at (2,2);
\coordinate (4) at (0,2);
\coordinate (5) at (1,0);
\coordinate (6) at (1,1);
\coordinate (7) at (0,1);
\coordinate (8) at (2,1);
\coordinate (9) at (1,2);
\coordinate (10) at (0.5, 0.5);
\coordinate (11) at (1.5,1.5);
\coordinate (12) at (1,0.5);
\coordinate (13) at (1.5,0.5);
\coordinate (14) at (1.5,1);
\draw[thick] (1)--(2)--(3)--(4)--(1);
\draw[thick] (6)--(9)--(7)--(6)--(5)--(8)--(6);
\draw[thick] (5)--(7);
\draw[thick] (10)--(12);
\draw[thick] (13)--(14)--(12)--(13)--(2);
\draw[thick](14)--(11);
\draw[thick] (9)--(8);
\draw[thick] (1)--(3);
\draw[thick,interface] (1)--(5);
\draw[thick,interface] (8)--(3);
\draw[thick,interface] (9)--(7);
\draw[thick,interface] (7)--(9);
\draw[thick,interface] (5)--(8);
\draw[thick,interface] (11)--(10);
\draw[thick,interface] (8)--(11);
\draw[thick,interface] (10)--(5);
\draw[thick,interface] (12)--(14)--(12);
\draw[thick,interface] (6)--(7)--(1);
\draw[thick,interface] (3)--(9)--(6);
\draw[thick,interface] (5)--(2)--(8);
\fill (1) circle (1.5pt);
\fill (2) circle (1.5pt);
\fill (3) circle (1.5pt);
\fill (4) circle (1.5pt);
\fill (5) circle (1.5pt);
\fill (6) circle (1.5pt);
\fill (7) circle (1.5pt);
\fill (8) circle (1.5pt);
\fill (9) circle (1.5pt);
\fill (10) circle (1.5pt);
\fill (11) circle (1.5pt);
\fill (12) circle (1.5pt);
\fill (13) circle (1.5pt);
\fill (14) circle (1.5pt);
\end{tikzpicture}}
 \end{minipage}
   \begin{minipage}{.15\textwidth}
   \vspace*{-0.65cm}
  \lstinline[basicstyle=\small\lstbasicfont]|coordinates|
\begin{lstlisting}[basicstyle=\small\lstbasicfont,aboveskip=2mm,belowskip=7.35mm,frame=lines,nolol=true]
   0    0
   2    0
   2    2
   0    2
   1    0
   1    1
   0    1
   2    1
   1    2
 0.5  0.5
 1.5  1.5
   1  0.5
 1.5  0.5
 1.5    1
\end{lstlisting}
  \end{minipage}\hspace*{2cm}
  \begin{minipage}{0.2\textwidth}
  \vspace*{0.15cm}
    \lstinline[basicstyle=\small\lstbasicfont]|elements|
\begin{lstlisting}[basicstyle=\small\lstbasicfont,aboveskip=2mm,frame=lines,nolol=true]
   7    1   10
   6    7   10
   9    6   11
   3    9   11
   7    9    4
   9    7    6
   8    3   11
  11    6   14
   8   11   14
  10    5   12
   6   10   12
   1    5   10
   2    8   13
   5    2   13
   5   13   12
  13    8   14
  12   14    6
  14   12   13
\end{lstlisting}
 \end{minipage}
 \end{minipage}
\end{minipage}}
 \caption{\revision{{\bfseries(A)} Initial triangulation with reference edges displayed as hatched lines. The array \matlab{coordinates} lists the x- and y-coordinates of the nodes, the corresponding indices are labeled in the meshes. The array \matlab{elements} specifies the element-connectivities by indexing the corresponding coordinates. The edge between the first two nodes in an element corresponds to the reference edge. {\bfseries(B)} Refined mesh obtained by marking both elements in the initial triangulation and performing an RGB refinement. New coordinates are appended to \matlab{coordinates} whereas new elements are stored in \matlab{elements} at the previous position of the unrefined element and successive positions. The rest is shifted by the amount of new included elements. Here, element numbers 1 to 4 are the red refinement of element 1 in (A), and element number 5 to 8 correspond to a red refinement of element number 2 in (A). New reference edges are highlighted and stored analogously. {\bfseries(C)} Adaptive triangulation obtained by RGB refining the mesh in (B) for the marked element number 8. This causes a CLOSURE step to eliminate arising hanging nodes as shown in Table~\ref{tab:mappingRGB}. New coordinates are appended to \matlab{coordinates}, reference edges are highlighted and stored as the edge between the first two nodes of an element, newly generated elements are stored at the previous position of the father element and the rest is shifted, e.g., the green refinement of element 1 in (B) corresponds to element numbers 1 and 2, the green refinement of element 2 in (B) corresponds to element number 3 and 4, etc.. }
 }
 \label{fig:example_mesh}
\end{figure}
}

\section{Coarsening Requirements: Red-Green-Blue Refinement vs. Newest Vertex Bisection}\label{sect:requirements}

The goal of geometric refinement is to reduce the element size by adding further nodes to a given triangulation. In other words, one wants to increase the number of degrees of freedom. Coarsening, conversely, decreases the number of degrees of freedom in a triangulation, i.\,e.\,, eliminates nodes of a triangulation. However, there are still some questions remaining. Let $\tilde{\mathcal{T}}$ be a refinement of a triangulation $\mathcal{T}$ satisfying shape regularity. How to eliminate nodes 
\begin{itemize}
\item to receive a triangular mesh (i.e., quadrilateral elements are not part of the triangulation)?
\item to receive a shape regular mesh (i.e., the inscribed ball condition is satisfied)?
\item to receive a \revision{conforming} triangulation (i.e., a triangulation without hanging nodes)?
\item to undo/invert a refinement without knowing the refinement history explicitly?
\end{itemize}
In literature, there are different approaches on coarsening - dealing with these details in different manners. The most common approach is to use edge collapsing known from Delaunay algorithms. This does not require to know the refinement history at all and the mesh quality is assured in the process of edge collapsing. We refer to \cite{bankxu,gooch}. Coarsening can also be done by clustering into regions via the centroidal Voronoi tesselation, cf. \cite{Shu2009}. The new mesh is then constructed via its dual -- a Delaunay triangulation. As a further coarsening algorithm, we would like to mention the concept of using the refinement history. More precisely, the history is used to invert the refinement procedure. Most works based on this approach use a hierarchical structure, i.e., store the refinement history explicitly, see e.\,g.\, \cite{kossa,alberta}. To the best of our knowledge, in 2D, there is only one work on non-hierarchical coarsening for the refinement procedure \emph{newest vertex bisection} by Chen and Zhang \cite{chenzhang}. \revision{The newest vertex bisection (NVB) differs from RGB refinement in one pattern. Instead of a red refinement, a bisec(3)-operation is used. I.\,e.\,instead of joining midpoints of the element's edges, the element is divided into four subtriangles by joining the midpoint of the reference edge to the vertex opposite to this edge and the midpoints of the remaining edges, cf. Figure~\ref{fig:NVBpattern}. Chen and Zhang} found an easy-to-verify criterion to determine whether nodes are allowed to be eliminated or not. This works well because NVB is implemented by a sequence of bisections and those can easily be undone. In other words, NVB consists of successive green refinements. The same is not true for the above mentioned red-green-blue refinement. We still see that the green pattern emerges from a bisection and the blue pattern arises from two subsequent bisections, cf. Figure~\ref{fig:TrefineRGBpattern}. However, the red pattern does not originate from a bisection of elements and thus the criterion proposed by Chen and Zhang fails to work for the RGB refinement. In this work, we discuss this issue and propose an easy-to-verify criterion to determine nodes for elimination in an RGB refined triangulation.
\revision{
\begin{figure}[!htb]
\centering
\begin{minipage}[t]{0.19\linewidth}
\hspace*{6mm}\begin{tikzpicture}[ interface/.style={
        postaction={draw,decorate,decoration={border,angle=45,
                    amplitude=0.13cm,segment length=1mm}}},
    ]
\coordinate (1) at (-1.2,0);
\coordinate (2) at (1.2,0);
\coordinate (3) at (0,1.5);
\draw (1)--(2) -- (3)--(1);
\fill (1) circle (2pt); 
\fill (2) circle (2pt); 
\fill (3) circle (2pt); 
\fill ($(1)!0.5!(2)$) circle (2pt);
\fill ($(2)!0.5!(3)$) circle (2pt);
\fill ($(1)!0.5!(3)$) circle (2pt);
\draw ($(2)!0.5!(3)$) -- ($(1)!0.5!(2)$);
\draw ($(1)!0.5!(3)$) -- ($(1)!0.5!(2)$);
\draw[interface,thick] (3)--($(1)!0.5!(2)$);
\draw[interface,thick] ($(1)!0.5!(2)$) --(2);
\draw[interface,thick] ($(1)!0.5!(2)$)--(3);
\draw[interface,thick] (1)--($(1)!0.5!(2)$);
\node at (0,-0.5) {bisec(3)};
\node at (-0.6,0.35) {2};
\node at (-0.3,0.75) {1};
\node at (0.6,0.35) {3};
\node at (0.3,0.75) {4};
\end{tikzpicture}
\end{minipage}
\caption{\revision{Newest vertex bisection differs from RGB refinement through the use of a bisec(3)-operation instead of the red pattern.  A bisec(3)-operation is essentially a green refinement of an element and each of the child elements is again green-refined.}}
\label{fig:NVBpattern}
\end{figure}
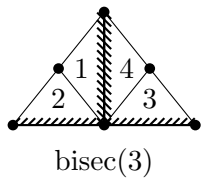}

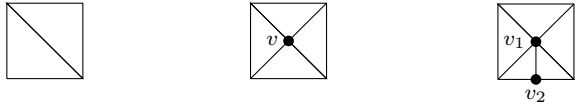
\begin{figure}[h]
\centering
\hfill
\begin{minipage}{0.15\textwidth}
\begin{tikzpicture}
\coordinate (1) at (0,0);
\coordinate (2) at (1,0);
\coordinate (3) at (1,1);
\coordinate (4) at (0,1);
\draw (1)--(2)--(4)--(1);
\draw (2)--(3)--(4)--(2);
\node[white,below] at (0.5,0) {\scriptsize 2};
\end{tikzpicture}
\end{minipage}
\begin{minipage}{0.15\textwidth}
\begin{tikzpicture}
\hspace*{5mm}
\coordinate (1) at (0,0);
\coordinate (2) at (1,0);
\coordinate (3) at (1,1);
\coordinate (4) at (0,1);
\draw (1)--(2)--(4)--(1);
\draw (2)--(3)--(4)--(2);
\draw (1)--(3);
\node[left] at (0.5,0.5) {\scriptsize $v$};
\fill (0.5,0.5) circle (2pt);
\node[white,below] at (0.5,0) {\scriptsize 2};
\end{tikzpicture}
\end{minipage}
\begin{minipage}{0.15\textwidth}
\begin{tikzpicture}
\hspace*{10mm}
\coordinate (1) at (0,0);
\coordinate (2) at (1,0);
\coordinate (3) at (1,1);
\coordinate (4) at (0,1);
\draw (1)--(2)--(4)--(1);
\draw (2)--(3)--(4)--(2);
\draw (1)--(3);
\draw (0.5,0.5)--(0.5,0);
\node[left] at (0.5,0.5) {\scriptsize $v_1$};
\fill (0.5,0.5) circle (2pt);
\fill (0.5,0) circle (2pt);
\node[below] at (0.5,0) {\scriptsize $v_2$};
\end{tikzpicture}
\end{minipage}\hfill
\caption{Left: Initial mesh. Middle: Newest vertex $v$ with $\#\mathcal{R}_v=4$. This vertex can be removed. Right: Newest vertices $v_1,v_2$ with $\#\mathcal{R}_{v_1}=5$ and $\#\mathcal{R}_{v_2}=2$. Only  $v_2$ can be removed.}
\label{fig:valence_chen}
\end{figure}

To this end, let us first investigate Chen and Zhang's approach to determine admissible nodes for the newest vertex bisection (in their paper called ``good-for-coarsening node''), cf.~\cite{chenzhang}. Defining the patch \[\mathcal{R}_v\coloneqq \left\{T \in \mathcal{T}~|~v \in T\right\},\] the \emph{valence} $\#\mathcal{R}_v$ counts the elements that are contiguous to a node $v \in \mathcal{N}$. Let \[\mathcal{N}_{\mathrm{new}}\coloneqq \left\{v \in \mathcal{N}\revision{\setminus\mathcal{N}_0}~|~v \text{ is newest vertex of some }T \in \mathcal{T} \right\}\] be the set of newest nodes in a triangulation $\mathcal{T}$. Chen and Zhang claim that the set of admissible nodes is characterized by the set \[\mathcal{N}_{\mathrm{adm}} \coloneqq \left\{v \in \mathcal{N}_{\mathrm{new}}~:~\#\mathcal{R}_v=4 \text{ or } \#\mathcal{R}_v=2\right\}.\] In Figure~\ref{fig:valence_chen}, this idea is illustrated. The set $\mathcal{N}_{\mathrm{adm}}$ is shown to be non-empty. This is an important requirement if one wants to assure that this criterion is useful in practical implementations. In short, a set of admissible nodes $\mathcal{N}_{\mathrm{adm}}$ is determined with this criterion. Adaptive coarsening can then be pursued by elimination of the set of nodes $\mathcal{N}_{\mathrm{adm}}\cap\mathcal{N}_{\mathrm{mark}}$ where $\mathcal{N}_{\mathrm{mark}}$ is the set of nodes that are marked through a given marking strategy. \revision{This method is powerful as it determines a set of nodes $\mathcal{N}_{\mathrm{adm}}\cap\mathcal{N}_{\mathrm{mark}}$ for which it is ensured that eliminating these nodes by joining elements together to its father element does not introduce any hanging nodes.} 

\begin{figure}[h]
\hfill
\centering
\begin{minipage}{0.15\textwidth}
\begin{tikzpicture}
\coordinate (1) at (0,0);
\coordinate (2) at (1,0);
\coordinate (3) at (1,1);
\coordinate (4) at (0,1);
\draw (1)--(2)--(4)--(1);
\draw (2)--(3)--(4)--(2);
\node[white,below] at (0.5,0) {\scriptsize 2};
\node[white,above] at (0.5,1) {\scriptsize 2};
\end{tikzpicture}
\end{minipage}
\begin{minipage}{0.15\textwidth}
\begin{tikzpicture}
\coordinate (1) at (0,0);
\coordinate (2) at (1,0);
\coordinate (3) at (1,1);
\coordinate (4) at (0,1);
\draw (1)--(2)--(4)--(1);
\draw (2)--(3)--(4)--(2);
\draw ($(2)!0.5!(4)$)--(3);
\draw ($(2)!0.5!(4)$)--($(1)!0.5!(4)$)--($(2)!0.5!(1)$)--($(2)!0.5!(4)$);
\node[above] at (0.5,0.5) {\scriptsize $v_1$};
\node[below] at (0.5,0) {\scriptsize $v_3$};
\node[left] at (0,0.5) {\scriptsize $v_2$};
\node[white,above] at (0.5,1) {\scriptsize 2};
\fill (0.5,0.5) circle (2pt);
\fill (0,0.5) circle (2pt);
\fill (0.5,0) circle (2pt);
\end{tikzpicture}
\end{minipage}
\begin{minipage}{0.15\textwidth}
\begin{tikzpicture}
\hspace*{5mm}
\coordinate (1) at (0,0);
\coordinate (2) at (1,0);
\coordinate (3) at (1,1);
\coordinate (4) at (0,1);
\draw (1)--(2)--(4)--(1);
\draw (2)--(3)--(4)--(2);
\draw ($(2)!0.5!(4)$)--($(3)!0.5!(4)$)--($(2)!0.5!(3)$)--($(2)!0.5!(4)$);
\draw ($(2)!0.5!(4)$)--($(1)!0.5!(4)$)--($(2)!0.5!(1)$)--($(2)!0.5!(4)$);
\node[above right] at (0.43,0.43) {\scriptsize $v_1$};
\node[below] at (0.5,0) {\scriptsize $v_3$};
\node[left] at (0,0.5) {\scriptsize $v_2$};
\node[above] at (0.5,1) {\scriptsize $v_5$};
\node[right] at (1,0.5) {\scriptsize $v_4$};
\fill (0.5,0.5) circle (2pt);
\fill (0,0.5) circle (2pt);
\fill (0.5,0) circle (2pt);
\fill (1,0.5) circle (2pt);
\fill (0.5,1) circle (2pt);
\end{tikzpicture}
\end{minipage}\hfill
\caption{Left: Initial mesh. Middle: Newest vertices $v_1,v_2,v_3$ with $\#\mathcal{R}_{v_1}=5$ and $\#\mathcal{R}_{v_2}=\#\mathcal{R}_{v_3}=3$. Right: Newest vertices $v_i, i=1,\ldots,5$ with $\#\mathcal{R}_{v_1}=6$ and $\#\mathcal{R}_{v_i}=3$ for $i=2,\ldots,5$. Although all nodes in both pictures could be removed the NVB-criteria cannot cover these cases, because this criterion demands $\#\mathcal{R}_v$ to be $2$ or $4$.}
\label{fig:valence_rgb}
\end{figure}
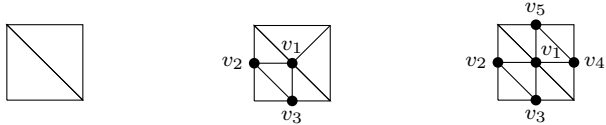

Let us now apply this easy-to-verify criterion for the red-green-blue refinement. One can easily see that green and blue refinements carry over, i.e., green refinements are removed for a valence of two or four and blue refinements are removed in a two-step-process - deleting one green refinement and then the subsequent one. This is favorable because, as already mentioned, adjacent green and blue patterns cannot be distinguished in our data structure. However, as blue patterns are not considered separately, but only as a sequence of green patterns, this does not pose any implementation problems. Applying this criterion to red patterns, it fails to detect the nodes that can be deleted, cf. Figure~\ref{fig:valence_rgb}. 

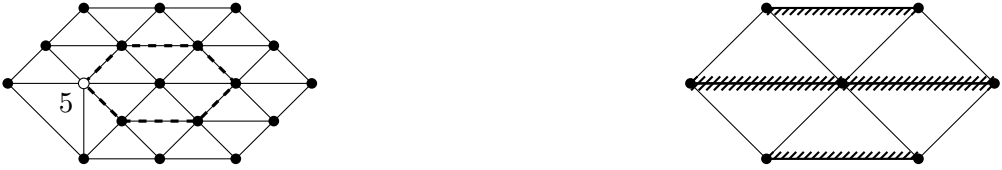
\begin{figure}
\begin{minipage}{0.47\textwidth}
\centering
\begin{tikzpicture}[ interface/.style={
        postaction={draw,decorate,decoration={border,angle=45,
                    amplitude=0.13cm,segment length=1mm}}},
    ]
\coordinate (1) at (0,0);
\coordinate (2) at (1,-1);
\coordinate (3) at (2,0);
\coordinate (4) at (1,1);
\coordinate (5) at (-1,1);
\coordinate (6) at (-2,0);
\coordinate (7) at (-1,-1);

\draw(1)--(2) -- (3)--(1)--(4)--(3);
\draw(4)--(5) -- (6)--(1)--(5);
\draw(2)--(7)--(6);
\draw (7)--(1);
\fill (1) circle (2pt); 
\fill (2) circle (2pt); 
\fill (3) circle (2pt); 
\fill (4) circle (2pt); 
\fill (5) circle (2pt); 
\fill (6) circle (2pt); 
\fill (7) circle (2pt); 
\draw ($(1)!0.5!(2)$) -- ($(2)!0.5!(3)$)--($(1)!0.5!(3)$)--($(1)!0.5!(2)$);
\fill ($(2)!0.5!(3)$) circle (2pt);
\draw ($(1)!0.5!(3)$) -- ($(4)!0.5!(3)$)--($(1)!0.5!(4)$)--($(1)!0.5!(3)$);
\fill ($(4)!0.5!(3)$) circle (2pt);
\draw ($(1)!0.5!(4)$) -- ($(4)!0.5!(5)$)--($(1)!0.5!(5)$)--($(1)!0.5!(4)$);
\fill ($(4)!0.5!(5)$) circle (2pt);
\draw ($(1)!0.5!(5)$) -- ($(5)!0.5!(6)$)--($(1)!0.5!(6)$)--($(1)!0.5!(5)$);
\fill ($(5)!0.5!(6)$) circle (2pt);
\draw ($(1)!0.5!(7)$) -- ($(1)!0.5!(6)$)--(7);
\draw ($(1)!0.5!(7)$) -- ($(2)!0.5!(7)$)--($(1)!0.5!(2)$)--($(1)!0.5!(7)$);
\fill ($(7)!0.5!(2)$) circle (2pt);
\draw[dashed,very thick] ($(1)!0.5!(6)$)--($(1)!0.5!(5)$)--($(1)!0.5!(4)$)--($(1)!0.5!(3)$)--($(1)!0.5!(2)$)--($(1)!0.5!(7)$)--($(1)!0.5!(6)$);
\fill ($(1)!0.5!(2)$) circle (2pt);
\fill ($(1)!0.5!(3)$) circle (2pt);
\fill ($(1)!0.5!(4)$) circle (2pt);
\fill ($(1)!0.5!(5)$) circle (2pt);
\fill ($(1)!0.5!(7)$) circle (2pt);
\node[below left] at ($(1)!0.5!(6)$) {5};
\node[circle,draw=black, fill=white,inner sep=0pt,minimum size=4pt] at ($(1)!0.5!(6)$) {};
\end{tikzpicture}
\end{minipage}\hfill
\begin{minipage}{0.47\textwidth}
\centering
\begin{tikzpicture}[ interface/.style={
        postaction={draw,decorate,decoration={border,angle=45,
                    amplitude=0.13cm,segment length=1mm}}},
    ]
\coordinate (1) at (0,0);
\coordinate (2) at (1,-1);
\coordinate (3) at (2,0);
\coordinate (4) at (1,1);
\coordinate (5) at (-1,1);
\coordinate (6) at (-2,0);
\coordinate (7) at (-1,-1);

\draw(1)--(2) -- (3)--(1)--(4)--(3);
\draw(4)--(5) -- (6)--(1)--(5);
\draw(2)--(7)--(6);
\draw (7)--(1);
\draw[interface,thick] (7) --(2);
\draw[interface,thick] (1) --(3);
\draw[interface,thick] (6) --(1);
\draw[interface,thick] (4) --(5);
\draw[interface,thick] (1) --(6);
\draw[interface,thick] (3) --(1);
\fill (1) circle (2pt); 
\fill (2) circle (2pt); 
\fill (3) circle (2pt); 
\fill (4) circle (2pt); 
\fill (5) circle (2pt); 
\fill (6) circle (2pt); 
\fill (7) circle (2pt); 
\end{tikzpicture}
\end{minipage}
\caption{\revision{Left: }Exemplary triangulation $(\mathcal{T},\refe_{\mathcal{T}})$. \revision{The vertex $v$ with $\# \tilde{\mathcal{R}}_v=5$ (in white) blocks all vertices along the loop (dotted) from deletion. These vertices are connected through red middle elements. Right: }Initial triangulation $(\mathcal{T}_0,\refe_{\mathcal{T}_0})$ of the left mesh.}
\label{fig:bsp_init}
\end{figure}

For this reason, a new criterion needs to be developed to cover red and green patterns at once. \revision{We closely follow the ideas from Chen and Zhang for NVB but incorporate the red middle element in our computations. We again consider the set of newest vertices \[ \mathcal{N}_{\mathrm{new}} \coloneqq \left\{v \in \mathcal{N}\setminus\mathcal{N}_0~|~v\text{ is \emph{newest vertex} of some } T \in \mathcal{T}\right\},\] cf.~Figure~\ref{fig:newestvertex}. Let further
\[\mathcal{M} \coloneqq \left\{T \in \mathcal{T}~|~T\text{ is a \emph{red middle element}} \right\} \] be the set of red middle elements in our triangulation $\mathcal{T}$. Red middle elements are determined by comparing the position of reference edges of neighbouring elements sharing an edge with this middle element, cf. Figure~\ref{fig:redmiddleelement}. The complement \[M^C \coloneqq \mathcal{T}\setminus \mathcal{M}\] is then used to define an adapted patch \[\tilde{\mathcal{R}}_v \coloneqq \left\{T \in \mathcal{T}~|~ v \in \mathcal{N}(T), T \in \mathcal{M}^C\right\}\] and with 
 \[\mathcal{N}_{\mathrm{candidates}} \coloneqq \left\{v \in \mathcal{N}_{\mathrm{new}}~:~\#\tilde{\mathcal{R}}_v=4 \text{ or } \#\tilde{\mathcal{R}}_v=2\right\}\] a set of candidates for elimination is found. In both cases of Figure~\ref{fig:valence_rgb} the adapted valence $\#\tilde{\mathcal{R}}_{v_1}=4$ and thus $v_1$ is considered a candidate for elimination. The valences for other vertices stay the same with this adapted definition of a patch.} Let us remark that checking only one node for a red pattern is inadequate. During a red refinement, three new nodes are created at once. To this end, it does not suffice to look at only one node, \revision{in contrast to NVB}. Moreover, checking all three nodes of a red pattern is not enough either, because the neighbouring pattern may be a red pattern and thus additional two nodes need to be taken into account. In the process of eliminating nodes, a whole chain of red patterns has to be followed to determine which nodes can actually be eliminated without creating a hanging node. \revision{To determine the set of admissible nodes $\mathcal{N}_{\mathrm{adm}}$ (here we consider $\mathcal{N}_\mathrm{adm}$ to be the set of nodes that does not create hanging nodes when eliminated) we have to follow this chain of red patterns and if for all nodes $v$ laying along this chain holds $v \in \mathcal{N}_{\mathrm{candidates}}$, it follows that $v \in \mathcal{N_\mathrm{adm}}$. An easy example shows that the so determined set $\mathcal{N}_\mathrm{adm}$ }may possibly be empty. Consider the triangulation $(\mathcal{T},\refe_{\mathcal{T}})$ shown in Figure~\ref{fig:bsp_init}. Here, $\mathcal{N}_{\mathrm{adm}}=\emptyset$ since the vertex $v$ with weight \revision{$\# \tilde{\mathcal{R}}_v=5$} blocks all vertices \revision{along the chain} from deletion. These vertices are connected through red middle elements along the loop.

Due to \revision{this property}, this method is not suitable for practical purposes as we may end up in a case where the mesh is not coarsened at all. In addition, we have not even considered adaptivity here. In contrast to Chen and Zhang's method we cannot use the set of nodes $\mathcal{N}_{\mathrm{adm}}\cap\mathcal{N}_{\mathrm{mark}}$ for adaptive deletion. In our case, we need to include $\mathcal{N}_{\mathrm{mark}}$ within the determination of admissible nodes $\mathcal{N}_{\mathrm{adm}}$ since we considered a whole chain of red patterns to avoid the creation of hanging nodes. If a node in this chain is not marked for deletion, it causes the same blockage as a node with valence unequal to two or four. 
\revision{In order to design a practically useful algorithm, we drop the requirement to avoid hanging nodes and rather work with a CLOSURE step. 
\begin{remark}
We see that even though NVB and RGB refinement only differ by one pattern and both refinement methods are easily implemented, finding a coarsening strategy for RGB is more involved without explicit knowledge of the refinement history. This is due to the red pattern and the resulting loss of a binary structure of the refinement history. In numerical experiments for refinement, no differences were found that would place one method above the other. However, the loss of a binary structure has more consequences, for example in the analysis of adaptive finite element methods. The analysis of convergence rates rely on a mesh overlay property, see \cite{bdd,cascon,stevensoncompletion} for the first contributions and \cite{axiomsofadaptivity} for an axiomatic contribution with a historical overview. This mesh overlay property is automatically fulfilled for binary tree refinement structures but does not hold for the RGB refinement as shown in \cite{pavlicek}.
\end{remark}
}

\section{The RGB Coarsening Algorithm}\label{sect:algorithm}

In this section, we present our RGB coarsening algorithm. We use the ideas presented in Section~\ref{sect:requirements} \revision{but loosen the conditions to the set $\mathcal{N}_\mathrm{adm}$. In the previous section, we considered $\mathcal{N}_\mathrm{adm}$ to be the set of nodes that does not create hanging nodes when eliminated. To this end, it was necessary to look at the chain of red patterns. Now, $\mathcal{N}_\mathrm{adm}$ declares the set of candidates for removal that are marked, cf. $\mathcal{N}_\mathrm{candidates}$ in Section~\ref{sect:requirements} with additional marking parameter, i.\,e.\,, the approach is local.} The main goal is to \revision{determine this set} efficiently with a non-hierarchical data structure. As soon as the pattern is determined, deleting the pattern is an easy task. \revision{This may introduce some hanging nodes. A subsequent CLOSURE step eliminates hanging nodes to obtain a conforming triangulation. This is a practical approach and guarantees} that coarsening is done locally. Algorithm~\ref{alg:coarsen} describes our coarsening algorithm with an additional CLOSURE step. Figure~\ref{fig:algo_illustr2} illustrates this procedure.

\begin{algorithm}[h!]
\caption{Coarsen a conforming RGB \revision{refined} triangulation $(\mathcal{T},\refe_{\mathcal{T}})$ \revision{of $(\mathcal{T}_0,\refe_{\mathcal{T}_0})$} locally at the nodes in $\mathcal{N}_{\mathrm{mark}}$.}\label{alg:coarsen}
\begin{algorithmic}[1]
\Procedure{Coarsen}{$\mathcal{T},\refe_{\mathcal{T}},\mathcal{N}_{\mathrm{mark}}$}
     \State $\mathcal{N}_{\mathrm{new}} \gets \left\{v \in \mathcal{N}\revision{\setminus\mathcal{N}_0}~|~v\text{ is \emph{newest vertex} of some }T \in \mathcal{T}\right\}$\Comment{see Figure~\ref{fig:newestvertex}}

   \State  $ \mathcal{M} \gets \left\{T \in \mathcal{T}~|~T\text{ is a \emph{red middle element}} \right\}$ \Comment{see Figure~\ref{fig:redmiddleelement}}
 \State $\mathcal{M}^C \gets  \mathcal{T}\setminus \mathcal{M}$
 \revision{\State $\tilde{\mathcal{R}}_v \gets \left\{T \in \mathcal{T}~|~ v \in \mathcal{N}(T), T \in \mathcal{M}^C\right\}$}
   \State $\mathcal{N}_{\mathrm{adm}} \gets\left\{v\in \mathcal{N}~\big|~\revision{\#\tilde{\mathcal{R}}}_v \in \left\{2,4\right\} \text{ and } v \in \mathcal{N}_{\mathrm{new}} \text{ and } v \in \mathcal{N}_{\mathrm{mark}}\right\}$
    \State $\mathcal{N}_{\mathrm{block}} \gets \mathcal{N}\setminus\mathcal{N}_{\mathrm{adm}}$
    \While{$\mathcal{N}_{\mathrm{block}}$ changes}\Comment{CLOSURE}
    \State $\mathcal{M}_{\mathrm{block}} \gets \left\{T \in \mathcal{M}~|~v \in \mathcal{N}_{\mathrm{block}}\cap\mathcal{N}(T)\right\}$
    \State $\mathcal{N}_{\mathrm{block}} \gets \mathcal{N}_{\mathrm{block}}\cup \left\{v \in \mathcal{N}(\mathcal{M}_{\mathrm{block}})~|~v \text{ is opposite to }\refe_{\mathcal{T}}(M), M \in \mathcal{M}_{\mathrm{block}}\right\}$ 
    \EndWhile
  \State $\mathcal{N}_{\mathrm{hang}} \gets \mathcal{N}_{\mathrm{block}} \cap \mathcal{N}_{\mathrm{new}}$
   \State $\mathcal{N}_{\mathrm{adm}} \gets \mathcal{N}_{\mathrm{adm}}\setminus\mathcal{N}_{\mathrm{hang}}$
   \State Create new elements and reference edges $(\hat{\mathcal{T}},\refe_{\hat{\mathcal{T}}})$ according to Figure~\ref{fig:creation}.
         \State \textbf{return} $(\hat{\mathcal{T}},\refe_{\hat{\mathcal{T}}})$
\EndProcedure
\end{algorithmic}
\end{algorithm}

\begin{figure}[h]
\centering
\begin{minipage}[t]{0.23\linewidth}
\hspace*{6mm}\begin{tikzpicture}[ interface/.style={
        postaction={draw,decorate,decoration={border,angle=45,
                    amplitude=0.13cm,segment length=1mm}}},
    ]
\coordinate (1) at (-1.2,0);
\coordinate (2) at (1.2,0);
\coordinate (3) at (0,1.5);
\draw(1)--(2) -- (3)--(1);
\draw[interface,thick] (1) --(2);
\fill (1) circle (2pt); 
\fill (2) circle (2pt); 
\fill (3) circle (2pt); 
\draw ($(1)!0.5!(2)$) -- ($(2)!0.5!(3)$)--($(1)!0.5!(3)$)--($(1)!0.5!(2)$);
\draw[interface,thick] ($(2)!0.5!(3)$)--($(1)!0.5!(3)$);
\draw[interface,thick] ($(1)!0.5!(3)$)--($(2)!0.5!(3)$);
\node[rectangle,draw=black, fill=white,inner sep=0pt,minimum size=4pt] at ($(1)!0.5!(2)$) {};
\node[rectangle,draw=black, fill=white,inner sep=0pt,minimum size=4pt] at ($(2)!0.5!(3)$) {};
\node[rectangle,draw=black, fill=white,inner sep=0pt,minimum size=4pt] at ($(1)!0.5!(3)$) {};
\draw[->,thick] (1.5,0.75)--(2.75,0.75);
\end{tikzpicture}
\end{minipage}
\begin{minipage}[t]{0.23\linewidth}
\hspace*{6mm}\begin{tikzpicture}[ interface/.style={
        postaction={draw,decorate,decoration={border,angle=45,
                    amplitude=0.13cm,segment length=1mm}}},
    ]
\coordinate (1) at (-1.2,0);
\coordinate (2) at (1.2,0);
\coordinate (3) at (0,1.5);
\draw(1) --(2)--(3)--(1);
\draw[interface,thick] (1) -- (2);
\fill (1) circle (2pt); 
\fill (2) circle (2pt); 
\fill (3) circle (2pt); 
\end{tikzpicture}
\end{minipage}\hspace{1cm}
\begin{minipage}[t]{0.23\linewidth}
\hspace*{6mm}\begin{tikzpicture}[ interface/.style={
        postaction={draw,decorate,decoration={border,angle=45,
                    amplitude=0.13cm,segment length=1mm}}},
    ]
\coordinate (1) at (-1.2,0);
\coordinate (2) at (1.2,0);
\coordinate (3) at (0,1.5);
\draw(1)--(2) -- (3)--(1);
\draw[interface,thick] (1) --(2);
\fill (1) circle (2pt); 
\fill (2) circle (2pt); 
\fill (3) circle (2pt); 
\draw ($(1)!0.5!(2)$) -- ($(2)!0.5!(3)$)--($(1)!0.5!(3)$)--($(1)!0.5!(2)$);
\draw[interface,thick] ($(2)!0.5!(3)$)--($(1)!0.5!(3)$);
\draw[interface,thick] ($(1)!0.5!(3)$)--($(2)!0.5!(3)$);
\node[circle,draw=black, fill=white,inner sep=0pt,minimum size=4pt] at ($(1)!0.5!(2)$) {};
\node[circle,draw=black, fill=white,inner sep=0pt,minimum size=4pt] at ($(2)!0.5!(3)$) {};
\node[circle,draw=black, fill=white,inner sep=0pt,minimum size=4pt] at ($(1)!0.5!(3)$) {};
\draw[->,thick] (1.5,0.75)--(2.75,0.75);
\end{tikzpicture}
\end{minipage}
\begin{minipage}[t]{0.23\linewidth}
\hspace*{6mm}\begin{tikzpicture}[ interface/.style={
        postaction={draw,decorate,decoration={border,angle=45,
                    amplitude=0.13cm,segment length=1mm}}},
    ]
\coordinate (1) at (-1.2,0);
\coordinate (2) at (1.2,0);
\coordinate (3) at (0,1.5);
\draw(1)--(2) -- (3)--(1);
\draw[interface,thick] (1) --(2);
\fill (1) circle (2pt); 
\fill (2) circle (2pt); 
\fill (3) circle (2pt); 
\fill ($(1)!0.5!(2)$) circle (2pt); 
\fill ($(3)!0.5!(2)$) circle (2pt); 
\fill ($(1)!0.5!(3)$) circle (2pt);
\draw ($(1)!0.5!(2)$) -- ($(2)!0.5!(3)$)--($(1)!0.5!(3)$)--($(1)!0.5!(2)$);
\draw[interface,thick] ($(2)!0.5!(3)$)--($(1)!0.5!(3)$);
\draw[interface,thick] ($(1)!0.5!(3)$)--($(2)!0.5!(3)$);
\end{tikzpicture}
\end{minipage}\\
\vspace*{5mm}
\begin{minipage}[t]{0.23\linewidth}
\hspace*{6mm}\begin{tikzpicture}[ interface/.style={
        postaction={draw,decorate,decoration={border,angle=45,
                    amplitude=0.13cm,segment length=1mm}}},
    ]
\coordinate (1) at (-1.2,0);
\coordinate (2) at (1.2,0);
\coordinate (3) at (0,1.5);
\draw(1)--(2) -- (3)--(1);
\draw[interface,thick] (1) --(2);
\fill (1) circle (2pt); 
\fill (2) circle (2pt); 
\fill (3) circle (2pt); 
\draw ($(1)!0.5!(2)$) -- ($(2)!0.5!(3)$)--($(1)!0.5!(3)$)--($(1)!0.5!(2)$);
\draw[interface,thick] ($(2)!0.5!(3)$)--($(1)!0.5!(3)$);
\draw[interface,thick] ($(1)!0.5!(3)$)--($(2)!0.5!(3)$);
\node[circle,draw=black, fill=white,inner sep=0pt,minimum size=4pt] at ($(1)!0.5!(2)$) {};
\node[rectangle,draw=black, fill=white,inner sep=0pt,minimum size=4pt] at ($(1)!0.5!(3)$) {};
\node[rectangle,draw=black, fill=white,inner sep=0pt,minimum size=4pt] at ($(2)!0.5!(3)$) {};
\draw[->,thick] (1.5,0.75)--(2.75,0.75);
\end{tikzpicture}
\end{minipage}
\begin{minipage}[t]{0.23\linewidth}
\hspace*{6mm}\begin{tikzpicture}[ interface/.style={
        postaction={draw,decorate,decoration={border,angle=45,
                    amplitude=0.13cm,segment length=1mm}}},
    ]
\coordinate (1) at (-1.2,0);
\coordinate (2) at (1.2,0);
\coordinate (3) at (0,1.5);
\draw(1) --(2)--(3)--(1);
\draw[interface,thick] (2) -- (3);
\draw[interface,thick] (3) -- (1);
\fill (1) circle (2pt); 
\fill (2) circle (2pt); 
\fill (3) circle (2pt); 
\draw ($(1)!0.5!(2)$) --(3);
\fill ($(1)!0.5!(2)$) circle (2pt);
\end{tikzpicture}
\end{minipage}\hspace{1cm}
\begin{minipage}[t]{0.23\linewidth}
\hspace*{6mm}\begin{tikzpicture}[ interface/.style={
        postaction={draw,decorate,decoration={border,angle=45,
                    amplitude=0.13cm,segment length=1mm}}},
    ]
\coordinate (1) at (-1.2,0);
\coordinate (2) at (1.2,0);
\coordinate (3) at (0,1.5);
\draw(1) --(2)--(3)--(1);
\draw[interface,thick] (2) -- (3);
\draw[interface,thick] (3) -- (1);
\fill (1) circle (2pt); 
\fill (2) circle (2pt); 
\fill (3) circle (2pt); 
\draw ($(1)!0.5!(2)$) --(3);
\node[rectangle,draw=black, fill=white,inner sep=0pt,minimum size=4pt] at ($(1)!0.5!(2)$) {};
\draw[->,thick] (1.5,0.75)--(2.75,0.75);
\end{tikzpicture}
\end{minipage}
\begin{minipage}[t]{0.23\linewidth}
\hspace*{6mm}\begin{tikzpicture}[ interface/.style={
        postaction={draw,decorate,decoration={border,angle=45,
                    amplitude=0.13cm,segment length=1mm}}},
    ]
\coordinate (1) at (-1.2,0);
\coordinate (2) at (1.2,0);
\coordinate (3) at (0,1.5);
\draw(1) --(2)--(3)--(1);
\draw[interface,thick] (1) -- (2);
\fill (1) circle (2pt); 
\fill (2) circle (2pt); 
\fill (3) circle (2pt); 
\end{tikzpicture}
\end{minipage}\\
\vspace*{5mm}
\begin{minipage}[t]{0.23\linewidth}
\hspace*{6mm}\begin{tikzpicture}[ interface/.style={
        postaction={draw,decorate,decoration={border,angle=45,
                    amplitude=0.13cm,segment length=1mm}}},
    ]
\coordinate (1) at (-1.2,0);
\coordinate (2) at (1.2,0);
\coordinate (3) at (0,1.5);
\draw(1)--(2) -- (3)--(1);
\draw[interface,thick] (1) --(2);
\fill (1) circle (2pt); 
\fill (2) circle (2pt); 
\fill (3) circle (2pt); 
\draw ($(1)!0.5!(2)$) -- ($(2)!0.5!(3)$)--($(1)!0.5!(3)$)--($(1)!0.5!(2)$);
\draw[interface,thick] ($(2)!0.5!(3)$)--($(1)!0.5!(3)$);
\draw[interface,thick] ($(1)!0.5!(3)$)--($(2)!0.5!(3)$);
\node[circle,draw=black, fill=white,inner sep=0pt,minimum size=4pt] at ($(1)!0.5!(2)$) {};
\node[rectangle,draw=black, fill=white,inner sep=0pt,minimum size=4pt] at ($(2)!0.5!(3)$) {};
\node[circle,draw=black, fill=white,inner sep=0pt,minimum size=4pt] at ($(1)!0.5!(3)$) {};
\draw[->,thick] (1.5,0.75)--(2.75,0.75);
\end{tikzpicture}
\end{minipage}
\begin{minipage}[t]{0.23\linewidth}
\hspace*{6mm}\begin{tikzpicture}[ interface/.style={
        postaction={draw,decorate,decoration={border,angle=45,
                    amplitude=0.13cm,segment length=1mm}}},
    ]
\coordinate (1) at (-1.2,0);
\coordinate (2) at (1.2,0);
\coordinate (3) at (0,1.5);
\draw(1) --(2)--(3)--(1);
\draw[interface,thick] (2) -- (3);
\draw[interface,thick] (1) -- ($(1)!0.5!(2)$);
\draw[interface,thick] ($(1)!0.5!(2)$)--(3);
\fill (1) circle (2pt); 
\fill (2) circle (2pt); 
\fill (3) circle (2pt); 
\draw ($(1)!0.5!(2)$) --(3);
\draw ($(1)!0.5!(2)$) --($(1)!0.5!(3)$);
\fill ($(1)!0.5!(2)$) circle (2pt);
\fill ($(3)!0.5!(1)$) circle (2pt);
\end{tikzpicture}
\end{minipage}\hspace{1cm}
\begin{minipage}[t]{0.23\linewidth}
\hspace*{6mm}\begin{tikzpicture}[ interface/.style={
        postaction={draw,decorate,decoration={border,angle=45,
                    amplitude=0.13cm,segment length=1mm}}},
    ]
\coordinate (1) at (-1.2,0);
\coordinate (2) at (1.2,0);
\coordinate (3) at (0,1.5);
\draw(1) --(2)--(3)--(1);
\draw[interface,thick] (2) -- (3);
\draw[interface,thick] (3) -- (1);
\fill (1) circle (2pt); 
\fill (2) circle (2pt); 
\fill (3) circle (2pt); 
\draw ($(1)!0.5!(2)$) --(3);
\node[circle,draw=black, fill=white,inner sep=0pt,minimum size=4pt] at ($(1)!0.5!(2)$) {};
\draw[->,thick] (1.5,0.75)--(2.75,0.75);
\end{tikzpicture}
\end{minipage}
\begin{minipage}[t]{0.23\linewidth}
\hspace*{6mm}\begin{tikzpicture}[ interface/.style={
        postaction={draw,decorate,decoration={border,angle=45,
                    amplitude=0.13cm,segment length=1mm}}},
    ]
\coordinate (1) at (-1.2,0);
\coordinate (2) at (1.2,0);
\coordinate (3) at (0,1.5);
\draw(1) --(2)--(3)--(1);
\draw[interface,thick] (2) -- (3);
\draw[interface,thick] (3) -- (1);
\fill (1) circle (2pt); 
\fill (2) circle (2pt); 
\fill (3) circle (2pt); 
\draw ($(1)!0.5!(2)$) --(3);
\fill ($(1)!0.5!(2)$) circle (2pt);
\end{tikzpicture}
\end{minipage}\\
\vspace*{5mm}
\begin{minipage}[t]{0.23\linewidth}
\hspace*{6mm}\begin{tikzpicture}[ interface/.style={
        postaction={draw,decorate,decoration={border,angle=45,
                    amplitude=0.13cm,segment length=1mm}}},
    ]
\coordinate (1) at (-1.2,0);
\coordinate (2) at (1.2,0);
\coordinate (3) at (0,1.5);
\draw(1)--(2) -- (3)--(1);
\draw[interface,thick] (1) --(2);
\fill (1) circle (2pt); 
\fill (2) circle (2pt); 
\fill (3) circle (2pt); 
\draw ($(1)!0.5!(2)$) -- ($(2)!0.5!(3)$)--($(1)!0.5!(3)$)--($(1)!0.5!(2)$);
\draw[interface,thick] ($(2)!0.5!(3)$)--($(1)!0.5!(3)$);
\draw[interface,thick] ($(1)!0.5!(3)$)--($(2)!0.5!(3)$);
\node[circle,draw=black, fill=white,inner sep=0pt,minimum size=4pt] at ($(1)!0.5!(2)$) {};
\node[circle,draw=black, fill=white,inner sep=0pt,minimum size=4pt] at ($(3)!0.5!(2)$) {};
\node[rectangle,draw=black, fill=white,inner sep=0pt,minimum size=4pt] at ($(1)!0.5!(3)$) {};
\draw[->,thick] (1.5,0.75)--(2.75,0.75);
\end{tikzpicture}
\end{minipage}
\begin{minipage}[t]{0.23\linewidth}
\hspace*{6mm}\begin{tikzpicture}[ interface/.style={
        postaction={draw,decorate,decoration={border,angle=45,
                    amplitude=0.13cm,segment length=1mm}}},
    ]
\coordinate (1) at (-1.2,0);
\coordinate (2) at (1.2,0);
\coordinate (3) at (0,1.5);
\draw(1) --(2)--(3)--(1);
\draw[interface,thick] ($(1)!0.5!(2)$) -- (2);
\draw[interface,thick] ($(1)!0.5!(2)$) -- (3);
\draw[interface,thick] (3) -- (1);
\fill (1) circle (2pt); 
\fill (2) circle (2pt); 
\fill (3) circle (2pt); 
\draw ($(1)!0.5!(2)$) --(3);
\draw ($(1)!0.5!(2)$) --($(3)!0.5!(2)$);
\fill ($(1)!0.5!(2)$) circle (2pt);
\fill ($(3)!0.5!(2)$) circle (2pt);
\end{tikzpicture}
\end{minipage}\hfill
\begin{minipage}[t]{0.23\linewidth}
\hspace*{6mm}\begin{tikzpicture}[ interface/.style={
        postaction={draw,decorate,decoration={border,angle=45,
                    amplitude=0.13cm,segment length=1mm}}},
    ]
\coordinate (1) at (-1.2,0);
\coordinate (2) at (1.2,0);
\coordinate (3) at (0,1.5);
\draw(1) --(2)--(3)--(1);
\draw[interface,thick] (1) -- (2);
\fill (1) circle (2pt); 
\fill (2) circle (2pt); 
\node[circle,draw=black, fill=white,inner sep=0pt,minimum size=4pt] at (3) {};
\draw[->,thick] (1.5,0.75)--(2.75,0.75);
\end{tikzpicture}
\end{minipage}
\begin{minipage}[t]{0.23\linewidth}
\hspace*{6mm}\begin{tikzpicture}[ interface/.style={
        postaction={draw,decorate,decoration={border,angle=45,
                    amplitude=0.13cm,segment length=1mm}}},
    ]
\coordinate (1) at (-1.2,0);
\coordinate (2) at (1.2,0);
\coordinate (3) at (0,1.5);
\draw(1) --(2)--(3)--(1);
\draw(1) --(2)--(3)--(1);
\draw[interface,thick] (1) -- (2);
\fill (1) circle (2pt); 
\fill (2) circle (2pt); 
\fill (3) circle (2pt); 
\end{tikzpicture}
\end{minipage}
\caption{Create new elements and reference edges $(\hat{\mathcal{T}},\refe_{\hat{\mathcal{T}}})$ according to the depiction for nodes in $\mathcal{N}_{\mathrm{adm}}$ \revision{(white squares)} and $\mathcal{N}_{\mathrm{hang}}$ \revision{(white dots). The patterns used are the same as allowed in the RGB refinement process. Properties such as the shape regularity are thus preserved.}}
\label{fig:creation}
\end{figure}
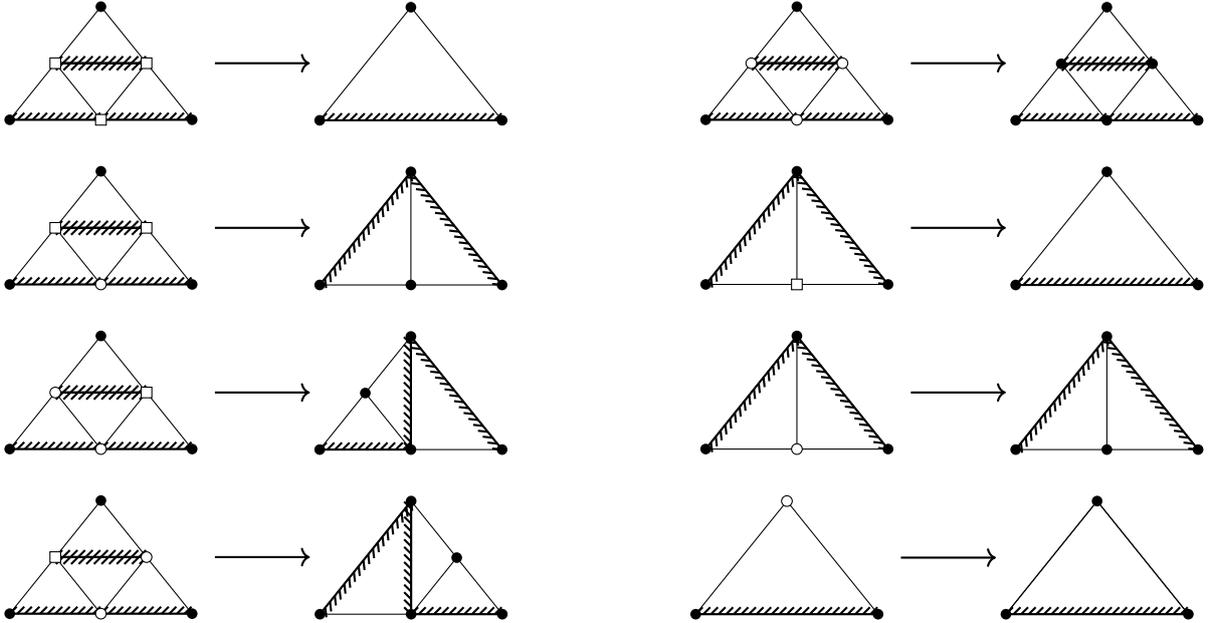

\input{algo_illus2}

\revision{
Let us elaborate whether a hanging node can be created through Algorithm~\ref{alg:coarsen}. A hanging node can be created by coarsening a pattern present at that node. Green patterns are fully removed. Red patterns may be coarsened into a green or blue pattern or they are fully removed. If the node is eliminated from the boundary, no hanging node is introduced. If this node is in the interior of the domain, it is shared by two neighbouring patterns. A hanging node is introduced if the connection to this node is eliminated on one side but not on the other. This might be the case if one pattern falls into the presented cases shown in Figure~\ref{fig:creation} and the other does not. We further argue that this cannot happen. First, all newest nodes $\mathcal{N}_{\mathrm{new}} \coloneqq \left\{v \in \mathcal{N}\revision{\setminus\mathcal{N}_0}~|~v\text{ is \emph{newest vertex} of some }T \in \mathcal{T}\right\}$ are collected. It might happen that a node $v$ is the newest node of one pattern but not of the other. The definition still includes this node, as it is the newest node of \emph{some} $T \in \mathcal{T}$. Figure~\ref{fig:possible} shows possible cases where this occurs. We see that for those cases holds that $\#\tilde{\mathcal{R}}_v \not\in \left\{2,4\right\}$ and $v$ is thus not considered for elimination. As a consequence neither the one nor the other pattern is coarsened at that node and thus no hanging node is introduced. Further, this can happen for a case shown in Figure~\ref{fig:possiblehanging}. But those cases cannot arise as only conforming triangulations are considered as input. There still remains the case, where this node is the newest node for both patterns. In the refinement process, we have red, green and blue patterns. In the coarsening step, we only coarsen red and green patterns. To this end, we examine what happens if blue patterns are present at a new node, see Figure~\ref{fig:possibleblue}. In this case, $\#\mathcal{R}_v \not\in \left\{2,4\right\}$ thus the connection to this node is not eliminated. The set $\mathcal{N}_\mathrm{adm}$ calculated first will only get smaller during the algorithm and therefore no new cases causing hanging nodes will occur. If a node is blocked, it is not eliminated by either pattern. If it is admissible, it is eliminated by both patterns. Overall, the algorithm generates a conforming triangulation. The shape regularity is preserved, as Lines~8-11 ensure that the reference edge is always marked and thus an element $T$ is coarsened into triangles that are in four similarity classes only. The inscribed ball condition is thus satisfied. This shows
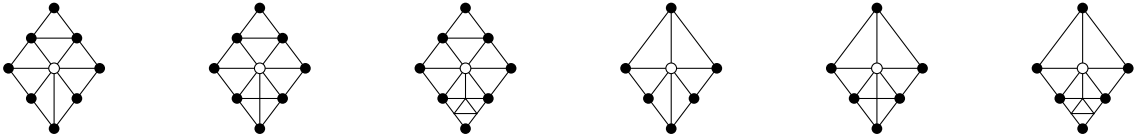
\begin{figure}
\centering
\hspace*{8mm}
\begin{minipage}{0.15\textwidth}
\begin{tikzpicture}
\coordinate (1) at (0,0);
\coordinate (2) at (1.2,0);
\coordinate (3) at (0.6,0.8);
\coordinate (7) at (0.6,-0.8);
\coordinate (4) at ($(1)!0.5!(2)$);
\coordinate (5) at ($(2)!0.5!(3)$);
\coordinate (6) at ($(1)!0.5!(3)$);
\coordinate (8) at ($(2)!0.5!(7)$);
\coordinate (9) at ($(1)!0.5!(7)$);
\draw(1) --(2)--(3)--(1);
\draw (6)--(4)--(5)--(6);
\draw (4)--(7)--(1);
\draw (7)--(2);
\draw (9)--(4)--(8);
\fill (1) circle (2pt); 
\fill (2) circle (2pt); 
\fill (3) circle (2pt); 
\node[circle,draw=black, fill=white,inner sep=0pt,minimum size=4pt] at (4) {};
\fill (5) circle (2pt); 
\fill (6) circle (2pt);
\fill (7) circle (2pt);
\fill (9) circle (2pt);
\fill (8) circle (2pt);
\end{tikzpicture}
\end{minipage}
\begin{minipage}{0.15\textwidth}
\begin{tikzpicture}
\coordinate (1) at (0,0);
\coordinate (2) at (1.2,0);
\coordinate (3) at (0.6,0.8);
\coordinate (7) at (0.6,-0.8);
\coordinate (4) at ($(1)!0.5!(2)$);
\coordinate (5) at ($(2)!0.5!(3)$);
\coordinate (6) at ($(1)!0.5!(3)$);
\coordinate (8) at ($(1)!0.5!(7)$);
\coordinate (9) at ($(2)!0.5!(7)$);
\draw(1) --(2)--(3)--(1);
\draw (6)--(4)--(5)--(6);
\draw (4)--(8)--(9)--(4);
\draw (1)--(7)--(2);
\draw (4)--(7);
\fill (1) circle (2pt); 
\fill (2) circle (2pt); 
\fill (3) circle (2pt); 
\node[circle,draw=black, fill=white,inner sep=0pt,minimum size=4pt] at (4) {};
\fill (5) circle (2pt); 
\fill (6) circle (2pt);
\fill (7) circle (2pt);
\fill (8) circle (2pt);
\fill (9) circle (2pt);
\end{tikzpicture}
\end{minipage}
\begin{minipage}{0.15\textwidth}
\begin{tikzpicture}
\coordinate (1) at (0,0);
\coordinate (2) at (1.2,0);
\coordinate (3) at (0.6,0.8);
\coordinate (7) at (0.6,-0.8);
\coordinate (4) at ($(1)!0.5!(2)$);
\coordinate (5) at ($(2)!0.5!(3)$);
\coordinate (6) at ($(1)!0.5!(3)$);
\coordinate (8) at ($(1)!0.5!(7)$);
\coordinate (9) at ($(2)!0.5!(7)$);
\coordinate (10) at ($(8)!0.5!(7)$);
\coordinate (11) at ($(9)!0.5!(7)$);
\coordinate (12) at ($(8)!0.5!(9)$);
\draw(1) --(2)--(3)--(1);
\draw (6)--(4)--(5)--(6);
\draw (4)--(8)--(9)--(4);
\draw (1)--(7)--(2);
\draw (4)--(12);
\draw (12)--(11)--(10)--(12);
\fill (1) circle (2pt); 
\fill (2) circle (2pt); 
\fill (3) circle (2pt); 
\node[circle,draw=black, fill=white,inner sep=0pt,minimum size=4pt] at (4) {};
\fill (5) circle (2pt); 
\fill (6) circle (2pt);
\fill (7) circle (2pt);
\fill (8) circle (2pt);
\fill (9) circle (2pt);
\end{tikzpicture}
\end{minipage}
\begin{minipage}{0.15\textwidth}
\begin{tikzpicture}
\coordinate (1) at (0,0);
\coordinate (2) at (1.2,0);
\coordinate (3) at (0.6,0.8);
\coordinate (7) at (0.6,-0.8);
\coordinate (4) at ($(1)!0.5!(2)$);
\coordinate (5) at ($(2)!0.5!(3)$);
\coordinate (6) at ($(1)!0.5!(3)$);
\coordinate (8) at ($(2)!0.5!(7)$);
\coordinate (9) at ($(1)!0.5!(7)$);
\draw(1) --(2)--(3)--(1);
\draw (4)--(3);
\draw (4)--(7)--(1);
\draw (7)--(2);
\draw (9)--(4)--(8);
\fill (1) circle (2pt); 
\fill (2) circle (2pt); 
\fill (3) circle (2pt); 
\node[circle,draw=black, fill=white,inner sep=0pt,minimum size=4pt] at (4) {};
\fill (7) circle (2pt);
\fill (9) circle (2pt);
\fill (8) circle (2pt);
\end{tikzpicture}
\end{minipage}
\begin{minipage}{0.15\textwidth}
\begin{tikzpicture}
\coordinate (1) at (0,0);
\coordinate (2) at (1.2,0);
\coordinate (3) at (0.6,0.8);
\coordinate (7) at (0.6,-0.8);
\coordinate (4) at ($(1)!0.5!(2)$);
\coordinate (5) at ($(2)!0.5!(3)$);
\coordinate (6) at ($(1)!0.5!(3)$);
\coordinate (8) at ($(1)!0.5!(7)$);
\coordinate (9) at ($(2)!0.5!(7)$);
\draw(1) --(2)--(3)--(1);
\draw (4)--(3);
\draw (4)--(8)--(9)--(4);
\draw (1)--(7)--(2);
\draw (4)--(7);
\fill (1) circle (2pt); 
\fill (2) circle (2pt); 
\fill (3) circle (2pt); 
\node[circle,draw=black, fill=white,inner sep=0pt,minimum size=4pt] at (4) {};
\fill (7) circle (2pt);
\fill (8) circle (2pt);
\fill (9) circle (2pt);
\end{tikzpicture}
\end{minipage}
\begin{minipage}{0.15\textwidth}
\begin{tikzpicture}
\coordinate (1) at (0,0);
\coordinate (2) at (1.2,0);
\coordinate (3) at (0.6,0.8);
\coordinate (7) at (0.6,-0.8);
\coordinate (4) at ($(1)!0.5!(2)$);
\coordinate (5) at ($(2)!0.5!(3)$);
\coordinate (6) at ($(1)!0.5!(3)$);
\coordinate (8) at ($(1)!0.5!(7)$);
\coordinate (9) at ($(2)!0.5!(7)$);
\coordinate (10) at ($(8)!0.5!(7)$);
\coordinate (11) at ($(9)!0.5!(7)$);
\coordinate (12) at ($(8)!0.5!(9)$);
\draw(1) --(2)--(3)--(1);
\draw (4)--(3);
\draw (4)--(8)--(9)--(4);
\draw (1)--(7)--(2);
\draw (4)--(12);
\draw (12)--(11)--(10)--(12);
\fill (1) circle (2pt); 
\fill (2) circle (2pt); 
\fill (3) circle (2pt); 
\node[circle,draw=black, fill=white,inner sep=0pt,minimum size=4pt] at (4) {}; 
\fill (7) circle (2pt);
\fill (8) circle (2pt);
\fill (9) circle (2pt);
\end{tikzpicture}
\end{minipage}
\caption{\revision{Possible cases for which a node $v$ (white dot) is the newest node of one pattern but not of the other. Here, the node $v$ is the newest node for the upper pattern but not for the lower one. We see that for all cases $\#\mathcal{R}_v \not\in \left\{2,4\right\}$ and thus those nodes are not considered for elimination.}}
\label{fig:possible}
\end{figure}

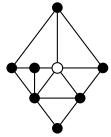
\begin{figure}
\centering
\hspace*{8mm}
\begin{minipage}{0.15\textwidth}
\begin{tikzpicture}
\coordinate (1) at (0,0);
\coordinate (2) at (1.2,0);
\coordinate (3) at (0.6,0.8);
\coordinate (7) at (0.6,-0.8);
\coordinate (4) at ($(1)!0.5!(2)$);
\coordinate (5) at ($(2)!0.5!(3)$);
\coordinate (6) at ($(1)!0.5!(3)$);
\coordinate (8) at ($(1)!0.5!(7)$);
\coordinate (9) at ($(2)!0.5!(7)$);
\draw(1) --(2)--(3)--(1);
\draw (4)--(3);
\draw (4)--(8)--(9)--(4);
\draw (1)--(7)--(2);
\draw (8)--($(1)!0.5!(4)$);
\fill (1) circle (2pt); 
\fill (2) circle (2pt); 
\fill (3) circle (2pt); 
\node[circle,draw=black, fill=white,inner sep=0pt,minimum size=4pt] at (4) {}; 
\fill (7) circle (2pt);
\fill (8) circle (2pt);
\fill (9) circle (2pt);
\fill ($(1)!0.5!(4)$) circle (2pt);
\end{tikzpicture}
\end{minipage}
\caption{\revision{Impossible situation for which the node (white dot) is the newest node for the upper pattern but not for the lower one. This situation cannot occur because only conforming triangulations are allowed as input parameters, i.\,e.\,, hanging nodes cannot exist.}}
\label{fig:possiblehanging}
\end{figure}

\begin{figure}
\centering
\hspace*{8mm}
\begin{minipage}{0.15\textwidth}
\begin{tikzpicture}
\coordinate (1) at (0,0);
\coordinate (2) at (1.2,0);
\coordinate (3) at (0.6,0.8);
\coordinate (7) at (0.6,-0.8);
\coordinate (4) at ($(1)!0.5!(2)$);
\coordinate (5) at ($(2)!0.5!(3)$);
\coordinate (6) at ($(1)!0.5!(3)$);
\coordinate (8) at ($(2)!0.5!(7)$);
\draw(1) --(2)--(3)--(1);
\draw (6)--(4)--(5)--(6);
\draw (4)--(7)--(1);
\draw (7)--(2);
\draw (8)--(4);
\fill (1) circle (2pt); 
\fill (2) circle (2pt); 
\fill (3) circle (2pt); 
\node[circle,draw=black, fill=white,inner sep=0pt,minimum size=4pt] at (4) {}; 
\fill (5) circle (2pt); 
\fill (6) circle (2pt);
\fill (7) circle (2pt);
\fill (8) circle (2pt);
\end{tikzpicture}
\end{minipage}
\begin{minipage}{0.15\textwidth}
\begin{tikzpicture}
\coordinate (1) at (0,0);
\coordinate (2) at (1.2,0);
\coordinate (3) at (0.6,0.8);
\coordinate (7) at (0.6,-0.8);
\coordinate (4) at ($(1)!0.5!(2)$);
\coordinate (5) at ($(2)!0.5!(3)$);
\coordinate (6) at ($(1)!0.5!(3)$);
\coordinate (8) at ($(2)!0.5!(7)$);
\draw(1) --(2)--(3)--(1);
\draw (4)--(3);
\draw (4)--(7)--(1);
\draw (7)--(2);
\draw (8)--(4);
\fill (1) circle (2pt); 
\fill (2) circle (2pt); 
\fill (3) circle (2pt); 
\node[circle,draw=black, fill=white,inner sep=0pt,minimum size=4pt] at (4) {};
\fill (7) circle (2pt);
\fill (8) circle (2pt);
\end{tikzpicture}
\end{minipage}
\begin{minipage}{0.15\textwidth}
\begin{tikzpicture}
\coordinate (1) at (0,0);
\coordinate (2) at (1.2,0);
\coordinate (3) at (0.6,0.8);
\coordinate (7) at (0.6,-0.8);
\coordinate (4) at ($(1)!0.5!(2)$);
\coordinate (5) at ($(2)!0.5!(3)$);
\coordinate (6) at ($(1)!0.5!(3)$);
\coordinate (8) at ($(2)!0.5!(7)$);
\draw(1) --(2)--(3)--(1);
\draw (4)--(3);
\draw (4)--(7)--(1);
\draw (7)--(2);
\draw (8)--(4);
\draw (4) -- ($(2)!0.5!(3)$);
\fill (1) circle (2pt); 
\fill (2) circle (2pt); 
\fill (3) circle (2pt); 
\fill ($(2)!0.5!(3)$) circle (2pt); 
\node[circle,draw=black, fill=white,inner sep=0pt,minimum size=4pt] at (4) {};
\fill (7) circle (2pt);
\fill (8) circle (2pt);
\end{tikzpicture}
\end{minipage}
\caption{\revision{Blue patterns present at a node $v$ (white dot). In these cases holds that $\#\mathcal{R}_v \not\in \left\{2,4\right\}$, i.\,e.\,, the node $v$ is not considered for elimination.}}
\label{fig:possibleblue}
\end{figure}
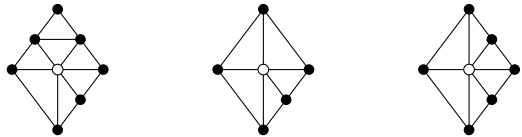
\begin{theorem}[Output COARSEN]
Let $(\mathcal{T},\refe_{\mathcal{T}})$ be a conforming triangulation obtained by RGB refinement of an initial conforming triangulation $(\mathcal{T}_0,\refe_{\mathcal{T}_0})$ and $\mathcal{N}_\mathrm{mark} \subset \mathcal{N}$. Then $\COARSEN(\mathcal{T}, \refe_{\mathcal{T}},\mathcal{N}_\mathrm{mark})$ from Algorithm~\ref{alg:coarsen} generates a conforming and shape regular triangulation. 
\end{theorem}
\begin{remark}
The CLOSURE step in Algorithm~\ref{alg:coarsen} might introduce new connectivities in the coarsened mesh. The shape regularity of the mesh is still preserved because the new connectivities are the same as the ones allowed in the refinement process. However, in adaptive methods, the required interpolation process is more involved by creating new connectivities, especially when evaluating non-nodal values. \cite{praetorius_wissgott} shows an interpolation approach for the red-green refinement. A similar ansatz can be used in our case.
\end{remark}
Our coarsening operation is not completely inverse to RGB refinement. A blue refinement of an element results in three child elements. One application of the coarsening algorithm does not coarsen the blue refinement to its father element but to two child elements of the father element. Additional patterns created during the CLOSURE step also do not follow an inverse operation of the RGB refinement. However, we can relate these patterns to a corresponding mesh obtained by NVB refinement. More specifically, as soon as one or two nodes of a red pattern are eliminated, our CLOSURE step does not go back to the father element but to blue or green children of this father element. We thus handle the mesh as if the mesh was NVB refined with this set of newly added nodes. The corresponding mesh can be defined as a bijective function and is discussed in detail in \cite{karkulik_extended}.\\
Even though the coarsening operation is not inverse the following result applies: Algorithm~\ref{alg:coarsen} can fully recover the initial triangulation provided that the initial triangulation $(\mathcal{T}_0,\refe_{\mathcal{T}_0})$ is of weak BDD-type.
\begin{definition}[weak BDD-property, cf.\cite{carstensen2004,karkulik}]\label{def:weak}
An element $T \in \mathcal{T}$ is called \emph{isolated} if the reference edge $\refe_\mathcal{T}(T)$ is shared with another element $\tilde{T} \in \mathcal{T}$ and $\refe_\mathcal{T}(\tilde{T}) \neq  \refe_\mathcal{T}(T)$. A mesh $(\mathcal{T},\refe_{\mathcal{T}})$ has the \emph{weak BDD-property} if two distinct isolated elements $T_1,T_2 \in \mathcal{T}$ do not share an edge.
\end{definition}
\begin{remark}
In adaptive meshing, it is common to impose conditions on the distribution of reference edges. In fact, the NVB coarsening algorithm of Chen and Zhang relies on a stricter condition on the initial triangulation, namely that there are no isolated elements at all, cf.~\cite{chenzhang}. This is the same condition as Binev, Dahmen and DeVore (BDD) imposed on the initial triangulation to prove optimal convergence rates for adaptive finite element methods with NVB \cite{bdd}. Carstensen weakened this condition to the above Definition~\ref{def:weak} to prove the $H^1$-stability of the $L^2$-projection for RGB refined meshes \cite{carstensen2004}. Later, Karkulik, Pavlicek and Praetorius improved these results in the sense that the condition of assignment of reference edges in the initial triangulation was removed \cite{karkulik}. For coarsening, we still need the weakened condition to guarantee the existence of nodes for elimination. Without any assumptions on the initial mesh, a loop of isolated elements may be formed that cannot be eliminated with our coarsening criteria, cf.~\cite{chenzhang}. This is not restrictive. In fact, Carstensen provides an algorithm that generates an extended conforming triangulation of weak BDD-type for an arbitrary conforming triangulation in linear complexity \cite{carstensen2004}. The results of Chen and Zhang for NVB remain also true under the weak BDD-assumption.
\end{remark}
With the weak BDD-property we can show
\begin{theorem}[Coarsening]
Let $(\mathcal{T},\refe_{\mathcal{T}})$ be an arbitrary RGB refinement of an initial conforming triangulation $(\mathcal{T}_0,\refe_{\mathcal{T}_0})$ where $(\mathcal{T}_0,\refe_{\mathcal{T}_0})$ is of weak BDD-type. Let $(\mathcal{T}^{(i)})_{i=0,1,\ldots}$ be a sequence of triangulations generated by Algorithm~\ref{alg:coarsen}, i.\,e.\,, 
\begin{equation*}\mathcal{T}^{(0)}\coloneqq \mathcal{T} \text{ and }\Big(\mathcal{T}^{(i+1)},\refe_{\mathcal{T}^{(i+1)}}\Big)\coloneqq \COARSEN\Big(\mathcal{T}^{(i)}, \refe_{\mathcal{T}^{(i)}},\mathcal{N}\big(\mathcal{T}^{(i)}\big)\Big).\end{equation*} Then, after a finite number of steps $M\in \mathbb{N}_0$, we obtain \begin{equation*} \mathcal{T}^{(M)} = \mathcal{T}_0.\end{equation*}
\end{theorem}
\begin{remark}
In practice, we would like to ensure that $M\leq cN$ for a small $c\geq 1$, where $N$ is the number of adaptive RGB refinement steps performed to obtain $(\mathcal{T},\refe_{\mathcal{T}})$. In Section~\ref{sect:MN}, some numerical experiments to estimate the ratio $\frac{M}{N}$ are performed.
\end{remark}
\begin{proof}
For the proof we consider a mesh refined by NVB as the coarsening patterns are chosen as if we would trace back the refinement history of a NVB refined mesh with the same nodes. RGB can then be related to NVB via a corresponding mesh function, cf.~\cite{karkulik}, i.\,e\,, the results are also valid for RGB. For a NVB refined mesh, there holds that $\tilde{\mathcal{R}}_v = \mathcal{R}_v$ for each node $v$ and thus relates to the work of Chen and Zhang \cite{chenzhang}. As shown in their work, only compatible patches are eliminated, see Figure~\ref{fig:compatible}. \\
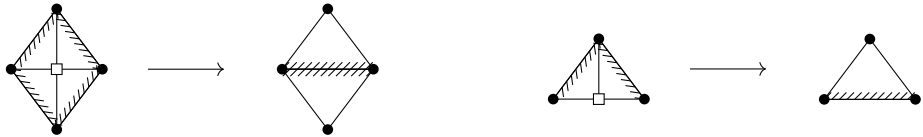
\begin{figure}[!htb]
\centering
\hspace*{15mm}
\begin{minipage}{0.2\textwidth}
\begin{tikzpicture}[ interface/.style={
       postaction={draw,decorate,decoration={border,angle=45,
                    amplitude=0.13cm,segment length=1mm}}},
    ]
\coordinate (1) at (0,0);
\coordinate (2) at (1.2,0);
\coordinate (3) at (0.6,0.8);
\coordinate (7) at (0.6,-0.8);
\coordinate (4) at ($(1)!0.5!(2)$);
\coordinate (5) at ($(2)!0.5!(3)$);
\coordinate (6) at ($(1)!0.5!(3)$);
\coordinate (8) at ($(2)!0.5!(7)$);
\draw[->] (1.8,0)--(2.8,0);
\draw(1) --(2)--(3)--(1);
\draw (4)--(3);
\draw (4)--(7)--(1);
\draw (7)--(2);
\draw[interface] (1)--(7)--(2)--(3)--(1);
\fill (1) circle (2pt); 
\fill (2) circle (2pt); 
\fill (3) circle (2pt); 
\node[rectangle,draw=black, fill=white,inner sep=0pt,minimum size=4pt] at (4) {};
\fill (7) circle (2pt);
\end{tikzpicture}
\end{minipage}
\begin{minipage}{0.2\textwidth}
\begin{tikzpicture}[ interface/.style={
       postaction={draw,decorate,decoration={border,angle=45,
                    amplitude=0.13cm,segment length=1mm}}},
    ]
\coordinate (1) at (0,0);
\coordinate (2) at (1.2,0);
\coordinate (3) at (0.6,0.8);
\coordinate (7) at (0.6,-0.8);
\coordinate (4) at ($(1)!0.5!(2)$);
\coordinate (5) at ($(2)!0.5!(3)$);
\coordinate (6) at ($(1)!0.5!(3)$);
\coordinate (8) at ($(2)!0.5!(7)$);
\draw(1) --(2)--(3)--(1);
\draw (1)--(7)--(2);
\draw[interface] (1)--(2);
\draw[interface] (2)--(1);
\fill (1) circle (2pt); 
\fill (2) circle (2pt); 
\fill (3) circle (2pt); 
\fill (7) circle (2pt);
\end{tikzpicture}
\end{minipage}
\begin{minipage}{0.2\textwidth}
\begin{tikzpicture}[ interface/.style={
       postaction={draw,decorate,decoration={border,angle=45,
                    amplitude=0.13cm,segment length=1mm}}},
    ]
\coordinate (1) at (0,0);
\coordinate (2) at (1.2,0);
\coordinate (3) at (0.6,0.8);
\coordinate (7) at (0.6,-0.8);
\coordinate (4) at ($(1)!0.5!(2)$);
\coordinate (5) at ($(2)!0.5!(3)$);
\coordinate (6) at ($(1)!0.5!(3)$);
\coordinate (8) at ($(2)!0.5!(7)$);
\draw(1) --(2)--(3)--(1);
\draw (4)--(3);
\draw[interface] (2)--(3)--(1);
\draw[->] (1.8,0.4)--(2.8,0.4);
\fill (1) circle (2pt); 
\fill (2) circle (2pt); 
\fill (3) circle (2pt); 
\node[rectangle,draw=black, fill=white,inner sep=0pt,minimum size=4pt] at (4) {};
\end{tikzpicture}
\end{minipage}
\begin{minipage}{0.2\textwidth}
\begin{tikzpicture}[ interface/.style={
       postaction={draw,decorate,decoration={border,angle=45,
                    amplitude=0.13cm,segment length=1mm}}},
    ]
\coordinate (1) at (0,0);
\coordinate (2) at (1.2,0);
\coordinate (3) at (0.6,0.8);
\coordinate (7) at (0.6,-0.8);
\coordinate (4) at ($(1)!0.5!(2)$);
\coordinate (5) at ($(2)!0.5!(3)$);
\coordinate (6) at ($(1)!0.5!(3)$);
\coordinate (8) at ($(2)!0.5!(7)$);
\draw(1) --(2)--(3)--(1);
\draw[interface] (1)--(2);
\fill (1) circle (2pt); 
\fill (2) circle (2pt); 
\fill (3) circle (2pt); 
\end{tikzpicture}
\end{minipage}
\caption{\revision{Compatible patches as shown in \cite{chenzhang}. The nodes $v$ (white squares) have the property $v \in \mathcal{N}_\mathrm{new}$ and $\#\mathcal{R}_v \in \left\{2,4\right\}$ and can thus be eliminated in a coarsening step. Reference edges are shown as hatched lines.}}
\label{fig:compatible}
\end{figure}
Let $(\mathcal{T}_0,\refe_{\mathcal{T}_0})$ be of weak BDD-type. We first show that any uniform bisec(3)-refinement of a weak BDD triangulation results into a weak BDD triangulation. For this purpose, we first recognize that elements with a reference edge as inner edge are not isolated due to the allowed refinement patterns. Thus, we only have to consider elements that share their reference edge with a triangle of another parent element. As illustrated in Figure~\ref{fig:weakBDD} the weak BDD property is inherited at these edges.\\
\begin{figure}[!htb]
\centering
\hspace*{15mm}
\begin{minipage}{0.43\textwidth}
\scalebox{0.8}{
\begin{tikzpicture}[ interface/.style={
       postaction={draw,decorate,decoration={border,angle=45,
                    amplitude=0.13cm,segment length=1mm}}},
    ]
\coordinate (1) at (0,0);
\coordinate (2) at (1.5,-3);
\coordinate (3) at (3,0);
\coordinate (4) at (1.5,3);
\coordinate (5) at (-1.5,3);
\coordinate (6) at (-3,0);
\coordinate (7) at (-1.5,-3);
\draw[fill=yellow] (1)--(5)--(6);
\draw[fill=yellow] (1)--(2)--(3);
\draw (2)--(3)--(4)--(5)--(6)--(7)--(2);
\draw (1)--(3);
\draw (1)--(6);
\draw[interface] (1)--(2);
\draw[interface] (4)--(7);
\draw[interface] (7)--(4);
\draw[interface] (1)--(5);
\fill (1) circle (2pt);
\fill (2) circle (2pt);
\fill (3) circle (2pt);
\fill (4) circle (2pt);
\fill (5) circle (2pt);
\fill (6) circle (2pt);
\fill (7) circle (2pt);
\node at ($($(1)!0.5!(6)$)!0.3!(5)$) {isolated};
\node at ($($(1)!0.5!(3)$)!0.3!(2)$) {isolated};
\end{tikzpicture}}
\end{minipage}
\begin{minipage}{0.43\textwidth}
\scalebox{0.8}{
\begin{tikzpicture}[ interface/.style={
       postaction={draw,decorate,decoration={border,angle=45,
                    amplitude=0.13cm,segment length=1mm}}},
    ]
\coordinate (1) at (0,0);
\coordinate (2) at (1.5,-3);
\coordinate (3) at (3,0);
\coordinate (4) at (1.5,3);
\coordinate (5) at (-1.5,3);
\coordinate (6) at (-3,0);
\coordinate (7) at (-1.5,-3);
\draw[fill=yellow] ($(1)!0.5!(2)$)--($(1)!0.5!(3)$)--(1);
\draw[fill=yellow] ($(3)!0.5!(2)$)--($(1)!0.5!(2)$)--(2);
\draw[fill=yellow] ($(5)!0.5!(6)$)--($(1)!0.5!(5)$)--(5);
\draw[fill=yellow] ($(1)!0.5!(5)$)--($(1)!0.5!(6)$)--(1);
\draw (2)--(3)--(4)--(5)--(6)--(7)--(2);
\draw (1)--(3);
\draw (1)--(6);
\draw[interface] (6)--($(1)!0.5!(5)$);
\draw[interface] ($(1)!0.5!(5)$)--(6);
\draw[interface] (3)--($(1)!0.5!(2)$);
\draw[interface] ($(1)!0.5!(2)$)--(3);
\draw[interface] (5)--(3);
\draw[interface] (3)--(5);
\draw[interface] (6)--(2);
\draw[interface] (2)--(6);
\draw[interface] (1)--(2);
\draw[interface] (4)--(7);
\draw[interface] (7)--(4);
\draw[interface] (1)--(5);
\draw ($(6)!0.5!(5)$)--($(1)!0.5!(5)$)--($(1)!0.5!(6)$);
\draw ($(4)!0.5!(5)$)--($(1)!0.5!(4)$)--($(1)!0.5!(5)$);
\draw ($(4)!0.5!(3)$)--($(1)!0.5!(4)$)--($(1)!0.5!(3)$);
\draw ($(3)!0.5!(2)$)--($(1)!0.5!(2)$)--($(1)!0.5!(3)$);
\draw ($(2)!0.5!(7)$)--($(1)!0.5!(7)$)--($(1)!0.5!(2)$);
\draw ($(6)!0.5!(7)$)--($(1)!0.5!(7)$)--($(1)!0.5!(6)$);
\fill (1) circle (2pt);
\fill (2) circle (2pt);
\fill (3) circle (2pt);
\fill (4) circle (2pt);
\fill (5) circle (2pt);
\fill (6) circle (2pt);
\fill (7) circle (2pt);
\fill ($(1)!0.5!(2)$) circle (2pt);
\fill ($(1)!0.5!(3)$) circle (2pt);
\fill ($(1)!0.5!(4)$) circle (2pt);
\fill ($(1)!0.5!(5)$) circle (2pt);
\fill ($(1)!0.5!(6)$) circle (2pt);
\fill ($(1)!0.5!(7)$) circle (2pt);
\fill ($(3)!0.5!(2)$) circle (2pt);
\fill ($(4)!0.5!(3)$) circle (2pt);
\fill ($(5)!0.5!(4)$) circle (2pt);
\fill ($(6)!0.5!(5)$) circle (2pt);
\fill ($(7)!0.5!(6)$) circle (2pt);
\fill ($(2)!0.5!(7)$) circle (2pt);
\node at ($($(1)!0.5!(6)$)!0.6!(5)$) {\tiny isolated};
\node at ($($(1)!0.5!(3)$)!0.6!(2)$) {\tiny isolated};
\node at ($($(1)!0.5!(2)$)!0.8!($(1)!0.25!(3)$)$) {\tiny isolated};
\node at ($($(1)!0.5!(5)$)!0.8!($(1)!0.25!(6)$)$) {\tiny isolated};
\end{tikzpicture}}
\end{minipage}
\caption{\revision{The weak BDD-property is inherited if the mesh is uniformly bisec(3)- or red-refined. Isolated elements are marked, reference edges are shown as hatched lines. Left: Initial mesh of weak BDD-type. Right: Uniform bisec(3)-refinement of the left mesh. The resulting mesh is still of weak BDD-type. Isolated elements exist only at the reference edges previously shared with another parent element with a different reference edge. A RGB refined mesh has the same property.}}
\label{fig:weakBDD}
\end{figure}
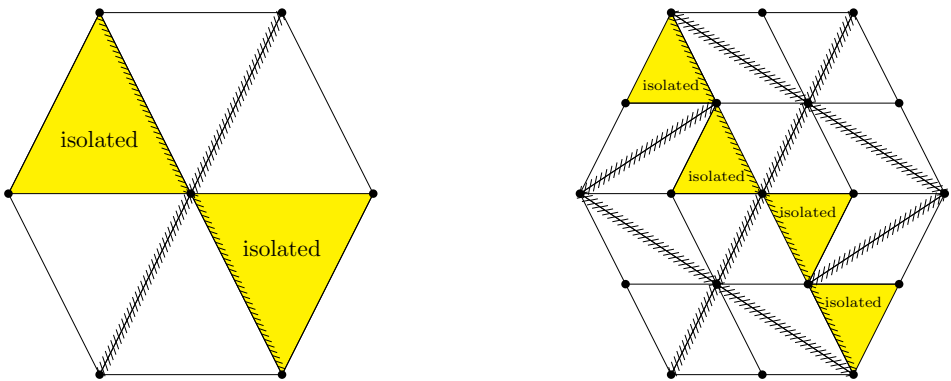
Therefore, we can restrict ourselves to a piece of the whole refinement and see if we can eliminate nodes to achieve the uniform refinement of a lower level. Figure~\ref{fig:eli} shows what happens for a piece of the whole refinement. For the highlighted nodes $\#\mathcal{R}_v \in \left\{2,4\right\}$ applies, otherwise two isolated elements would have shared an edge in the initial triangulation, i.\,e.\,, the mesh would not be of weak BDD type. Therefore, these nodes can be eliminated. In a further coarsening step, we again have a compatible patch that can be coarsened. In a last step, with the same arguments as above, $\#\mathcal{R}_v \in \left\{2,4\right\}$ and can thus be eliminated. As long as $\mathcal{T}^{(i)} \neq \mathcal{T}_0$, this process can be repeated until there are no more newest nodes. In this way, the initial mesh is recovered after a finite number of steps. \qed
\begin{figure}[!htb]
\centering
\hspace*{8mm}
\begin{minipage}{0.4\textwidth}
\scalebox{0.8}{
\begin{tikzpicture}[ interface/.style={
       postaction={draw,decorate,decoration={border,angle=45,
                    amplitude=0.13cm,segment length=1mm}}},
    ]
\coordinate (1) at (0,0);
\coordinate (2) at (1.5,-3);
\coordinate (3) at (3,0);
\coordinate (4) at (1.5,3);
\coordinate (5) at (-1.5,3);
\draw (3)--(4)--(5)--(6)--(1)--(3);
\draw ($(6)!0.5!(5)$)--($(1)!0.5!(5)$)--($(1)!0.5!(6)$);
\draw ($(1)!0.5!(5)$)--($(5)!0.5!(4)$)--($(1)!0.5!(4)$);
\draw ($(3)!0.5!(4)$)--($(1)!0.5!(4)$)--($(1)!0.5!(3)$);
\draw[interface] (1)--(5);
\draw[interface] (4)--(1);
\draw[interface] (4)--(5);
\draw[interface] ($(1)!0.5!(5)$)--(6);
\draw[interface] ($(5)!0.5!(4)$)--(1);
\draw[interface] ($(1)!0.5!(4)$)--(3);
\draw[interface] (6)--($(1)!0.5!(5)$);
\draw[interface] (1)--($(5)!0.5!(4)$);
\draw[interface] (3)--($(1)!0.5!(4)$);
\fill (1) circle (2pt);
\fill (3) circle (2pt);
\fill (4) circle (2pt);
\fill (5) circle (2pt);
\fill (6) circle (2pt);
\fill ($(1)!0.5!(5)$) circle (2pt);
\fill ($(1)!0.5!(4)$) circle (2pt);
\fill ($(4)!0.5!(5)$) circle (2pt);
\node[rectangle,draw=black, fill=white,inner sep=0pt,minimum size=4pt] at ($(1)!0.5!(6)$) {};
\node[rectangle,draw=black, fill=white,inner sep=0pt,minimum size=4pt] at ($(5)!0.5!(6)$) {};
\node[rectangle,draw=black, fill=white,inner sep=0pt,minimum size=4pt] at ($(3)!0.5!(4)$) {};
\node[rectangle,draw=black, fill=white,inner sep=0pt,minimum size=4pt] at ($(3)!0.5!(1)$) {};
\node at ($($(1)!0.5!(6)$)!0.6!(5)$) {\tiny isolated};
\node at ($($(1)!0.5!(3)$)!0.6!(4)$) {\tiny isolated};
\node at ($($(1)!0.5!(4)$)!0.8!($(1)!0.25!(3)$)$) {\tiny isolated};
\node at ($($(1)!0.5!(5)$)!0.8!($(1)!0.25!(6)$)$) {\tiny isolated};
\draw[->] (3.5,1.5)--(5,1.5);
\end{tikzpicture}}
\end{minipage}
\begin{minipage}{0.4\textwidth}
\scalebox{0.8}{
\begin{tikzpicture}[ interface/.style={
       postaction={draw,decorate,decoration={border,angle=45,
                    amplitude=0.13cm,segment length=1mm}}},
    ]
\coordinate (1) at (0,0);
\coordinate (2) at (1.5,-3);
\coordinate (3) at (3,0);
\coordinate (4) at (1.5,3);
\coordinate (5) at (-1.5,3);
\draw (3)--(4)--(5)--(6)--(1)--(3);
\draw[interface] (5)--(6);
\draw[interface] (6)--(1);
\draw[interface] (1)--(3);
\draw[interface] (3)--(4);
\draw[interface] (4)--(5);
\draw (5)--(1)--(4);
\draw (6)--($(5)!0.5!(1)$)--($(5)!0.5!(4)$)--($(1)!0.5!(4)$)--(3);
\fill (1) circle (2pt);
\fill (3) circle (2pt);
\fill (4) circle (2pt);
\fill (5) circle (2pt);
\fill (6) circle (2pt);
\fill ($(5)!0.5!(4)$) circle (2pt);
\draw[interface] ($(5)!0.5!(4)$)--(1);
\draw[interface] (1)--($(5)!0.5!(4)$);
\node[rectangle,draw=black, fill=white,inner sep=0pt,minimum size=4pt] at ($(1)!0.5!(4)$) {};
\node[rectangle,draw=black, fill=white,inner sep=0pt,minimum size=4pt] at ($(5)!0.5!(1)$) {};
\draw[->] (3.5,1.5)--(5,1.5);
\end{tikzpicture}}
\end{minipage}\\
\vspace*{5mm}
\hspace*{8mm}
\begin{minipage}{0.4\textwidth}
\scalebox{0.8}{
\begin{tikzpicture}[ interface/.style={
       postaction={draw,decorate,decoration={border,angle=45,
                    amplitude=0.13cm,segment length=1mm}}},
    ]
\coordinate (1) at (0,0);
\coordinate (2) at (1.5,-3);
\coordinate (3) at (3,0);
\coordinate (4) at (1.5,3);
\coordinate (5) at (-1.5,3);
\draw (3)--(4)--(5)--(6)--(1)--(3);
\draw[interface] (1)--(5);
\draw[interface] (4)--(1);
\draw[interface] (5)--(1);
\draw[interface] (1)--(4);
\fill (1) circle (2pt);
\fill (3) circle (2pt);
\fill (4) circle (2pt);
\fill (5) circle (2pt);
\fill (6) circle (2pt);
\draw ($(5)!0.5!(4)$)--(1);
\node[rectangle,draw=black, fill=white,inner sep=0pt,minimum size=4pt] at ($(5)!0.5!(4)$) {};
\draw[->] (3.5,1.5)--(5,1.5);
\end{tikzpicture}}
\end{minipage}
\begin{minipage}{0.4\textwidth}
\scalebox{0.8}{
\begin{tikzpicture}[ interface/.style={
       postaction={draw,decorate,decoration={border,angle=45,
                    amplitude=0.13cm,segment length=1mm}}},
    ]
\coordinate (1) at (0,0);
\coordinate (2) at (1.5,-3);
\coordinate (3) at (3,0);
\coordinate (4) at (1.5,3);
\coordinate (5) at (-1.5,3);
\draw (3)--(4)--(5)--(6)--(1)--(3);
\draw[interface] (1)--(5);
\draw[interface] (4)--(1);
\draw[interface] (4)--(5);
\fill (1) circle (2pt);
\fill (3) circle (2pt);
\fill (4) circle (2pt);
\fill (5) circle (2pt);
\fill (6) circle (2pt);
\node at ($($(1)!0.5!(6)$)!0.3!(5)$) {isolated};
\node at ($($(1)!0.5!(3)$)!0.3!(4)$) {isolated};
\end{tikzpicture}}
\end{minipage}
\caption{\revision{Piece of a uniformly bisec(3)-refined triangulation of weak BDD-type. Here, two isolated elements were connected by one element. From top left to bottom right: Uniformly refined mesh with nodes (white squares) that can be coarsened. Subsequently coarsened meshes where the nodes marked in the previous mesh have been eliminated and further nodes (white squares) are determined for the next elimination. After three coarsening steps, a uniform refinement of a lower level is reached. }}
\label{fig:eli}
\end{figure}
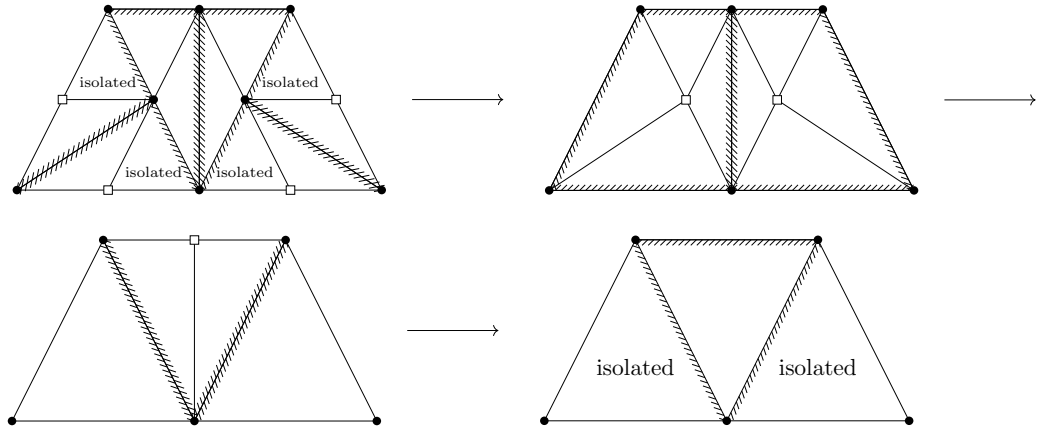
\end{proof}
\begin{remark}
At most $M= \#\mathcal{N}(\mathcal{T})-\#\mathcal{N}_0$ are needed to recover the initial triangulation. Numerical experiments show that this bound is very pessimistic. 
\end{remark}
}

\section{Matlab Implementation of the Coarsening Algorithm}\label{sect:Matlab}

\begin{figure}[h]\begin{minipage}{\linewidth}
\centering
\begin{tikzpicture}
\coordinate (1) at (-1.5,0);
\coordinate (2) at (1.5,0);
\coordinate (3) at (0,1.5);
\coordinate (4) at (0,0);
\draw(1)--(2) -- (3)--(1);
\draw (4)--(3);
\draw[fill=pastelred] ($(1)!0.5!(4)$)--($(4)!0.5!(3)$)--($(1)!0.5!(3)$)--($(1)!0.5!(4)$);
\draw[fill=pastelred] ($(2)!0.5!(4)$)--($(4)!0.5!(3)$)--($(2)!0.5!(3)$)--($(2)!0.5!(4)$);
\fill (4) circle (2pt);
\node at (-1,0.25) {2};
\node at (-0.25,0.25) {3};
\node at (-0.25,1) {1};
\node at (-0.5,0.5) {4};
\node at (1,0.25) {5};
\node at (0.25,0.25) {7};
\node at (0.25,1) {6};
\node at (0.5,0.5) {8};
\end{tikzpicture}
\end{minipage}
\caption{Due to the systematic error, situations arise where the valence $\revision{\#\tilde{\mathcal{R}}}_v=2$ and the elements containing $v$ are not numbered consecutively (3 and 7). Analogously, the same situation can arise for $\revision{\#\tilde{\mathcal{R}}}_v=4$. }
\label{fig:bsp_systematic}
\end{figure}
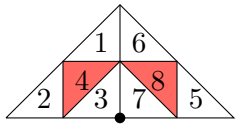

In the previous section, we have presented our RGB coarsening algorithm. In this section, we focus on the concrete implementation in MATLAB based on the refinement routine \texttt{TrefineRGB.m}, see \cite{web}. We have already discussed the data structure used in the refinement procedure in Section~\ref{sect:RGBrefine}. This will also play a major role in the implementation of the coarsening routine. Let us recap quickly the main structures: Elements $T \in \mathcal{T}$ are defined by their vertices $v_i, \,i = 1,2,3$ and numbered counterclockwise. The edge in between the first two vertices $v_1,v_2$ is the reference edge $\refe_\mathcal{T}(T)=v_1v_2$ of $T$. Refined elements are stored at the former position of the father element and subsequent positions, cf.~Figure~\ref{fig:enumeration}. Elements that belong together are therefore listed one after the other. RGB coarsening can then be implemented as follows (see Listing~\ref{lst:coarsen}):
\begin{itemize}
\item Lines 1--4: The function \texttt{TcoarsenRGB} expects the number of coordinates $N_0$ of the initial triangulation $\mathcal{T}_0$, mesh information such as \texttt{coordinates} and \texttt{elements} and optionally boundary data in \texttt{varargin} as input. The last argument in \texttt{varargin} is reserved for marked elements (Line 4). \revision{It is sufficient to pass the number of coordinates $N_0$ of the initial triangulation, since the coordinates added during the refinement process are appended to the end of the array \texttt{coordinates}. This means that the first $N_0$ entries in \texttt{coordinates} are the nodes $\mathcal{N}_0$.}
\item Lines 5--7: A triangulation is represented by the array \texttt{elements} and \texttt{coordinates}. The auxiliary functions \texttt{provideGeometricData} and \texttt{createEdge2Elements\_adv} provide more geometric information on the mesh. The array \texttt{element2edges} specifies the edges of each element, cf.~\cite{funkenschmidt}, and the array \texttt{edge2elements} specifies the elements containing this edge and the position of this edge within an element for all edges. E.\,g.\,, for \[\texttt{element2edges} = \begin{bmatrix} 3&1&2\\3&4&5\end{bmatrix}\text{ holds } \texttt{edge2elements} = \begin{bmatrix} 1&2&0&0\\1&3&0&0\\1&1&2&1\\2&2&0&0\\2&3&0&0\end{bmatrix}.\] \texttt{edge2elements($\ell$,1:2)} specifies the position (column, row) of the $\ell$-th edge in \texttt{element2edges}. If the $\ell-$th edge is a boundary edge, it holds \texttt{edge2elements($\ell$,3:4)}$=[0,0]$. Otherwise it specifies the position (column, row) of the $\ell$-th edge's second entry in \texttt{element2edges}.
\item Lines 9--10: \revision{(\ref{itm:4})} We determine all red middle elements as depicted in Figure~\ref{fig:redmiddleelement} by comparing the edges of four consecutive entries in \texttt{element2edges}. If the second edge of element one is equal to the second edge of element four, the third edge of element two is equal to the third edge of element four and the first edge of element three is equal to the first edge of element four, then \texttt{sum(abs($\ldots$))=0} and therefore the index of the red middle element is given by \texttt{find(sum(abs($\ldots$))=0)+3}, cf.~Figure~\ref{fig:rednumbering}. Note, that this characterizes a red pattern uniquely and thus this criterion to find red middle elements is appropriate.
\begin{figure}[h]
\centering
\begin{tikzpicture}[ interface/.style={
        postaction={draw,decorate,decoration={border,angle=45,
                    amplitude=0.13cm,segment length=1mm}}},
    ]
\coordinate (1) at (-3,0);
\coordinate (2) at (3,0);
\coordinate (3) at (0,3);
\draw(1)--(2) -- (3)--(1);
\draw[fill=pastelred] ($(1)!0.5!(2)$) -- ($(2)!0.5!(3)$)--($(1)!0.5!(3)$)--($(1)!0.5!(2)$);
\draw[interface,thick] ($(2)!0.5!(3)$)--($(1)!0.5!(3)$);
\draw[interface,thick] ($(1)!0.5!(3)$)--($(2)!0.5!(3)$);
\draw[interface,thick] (1) --(2);
\fill (1) circle (2pt); 
\fill (2) circle (2pt); 
\fill (3) circle (2pt); 
\fill ($(1)!0.5!(2)$) circle (2pt);
\fill ($(2)!0.5!(3)$) circle (2pt);
\fill ($(1)!0.5!(3)$) circle (2pt);
\node at (-1.5,0.75) {\bfseries 1};
\node at (1.5,0.75) {\bfseries 2};
\node at (0,2.25) {\bfseries 3};
\node at (0,0.75) {\bfseries 4}; 
\node[right] at ($(1)!0.5!($(1)!0.5!(3)$)$) {\scriptsize 3};
\node[above] at ($(1)!0.5!($(1)!0.5!(2)$)$) {\scriptsize 1};
\node[left] at ($($(1)!0.5!(2)$)!0.5!($(1)!0.5!(3)$)$) {\scriptsize 2};
\node[above] at ($(2)!0.5!($(1)!0.5!(2)$)$) {\scriptsize 1};
\node[left] at ($(2)!0.5!($(2)!0.5!(3)$)$) {\scriptsize 2};
\node[right] at ($($(2)!0.5!(3)$)!0.5!($(1)!0.5!(2)$)$) {\scriptsize 3};
\node[above] at ($($(1)!0.5!(3)$)!0.5!($(2)!0.5!(3)$)$) {\scriptsize 1};
\node[left] at ($(3)!0.5!($(2)!0.5!(3)$)$) {\scriptsize 2};
\node[right] at ($(3)!0.5!($(1)!0.5!(3)$)$) {\scriptsize 3};
\node[below] at ($($(1)!0.5!(3)$)!0.5!($(2)!0.5!(3)$)$) {\scriptsize 1};
\node[left] at ($($(2)!0.5!(3)$)!0.5!($(1)!0.5!(2)$)$) {\scriptsize 3};
\node[right] at ($($(2)!0.5!(1)$)!0.5!($(1)!0.5!(3)$)$) {\scriptsize 2};
\end{tikzpicture}
\caption{Numbering of red pattern: Element numbers are in bold, edge numbers per element are shown in smaller font size.}
\label{fig:rednumbering}
\end{figure}
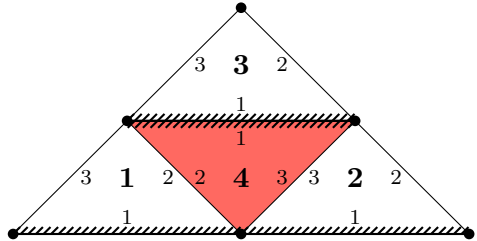
\item Lines 11--13: \revision{(\ref{itm:3})} In this part, the newest node of each element is detected. The newest node of an element is stored on position three in \texttt{elements}. We only consider nodes for coarsening that are not part of the initial triangulation $\mathcal{T}_0$ (Line 12). Note that we make a systematic error for red patterns, as an additional node is detected as newest node even though it is not a new one, cf.~Figure~\ref{fig:systematic}. We have to consider this systematic error in the course of our implementation. If an element is marked for coarsening, we mark all nodes of this element for coarsening (Line 13).
\begin{figure}[h]
\centering
\hspace*{20mm}
\begin{tikzpicture}[ interface/.style={
        postaction={draw,decorate,decoration={border,angle=45,
                    amplitude=0.13cm,segment length=1mm}}},
    ]
\coordinate (1) at (-1.2,0);
\coordinate (2) at (1.2,0);
\coordinate (3) at (0,1.5);
\draw(1)--(2) -- (3)--(1);
\draw[interface,thick] (1) --(2);
\fill (1) circle (2pt); 
\fill (2) circle (2pt); 
\fill (3) circle (2pt); 
\draw ($(1)!0.5!(2)$) -- ($(2)!0.5!(3)$)--($(1)!0.5!(3)$)--($(1)!0.5!(2)$);
\draw[interface,thick] ($(2)!0.5!(3)$)--($(1)!0.5!(3)$);
\draw[interface,thick] ($(1)!0.5!(3)$)--($(2)!0.5!(3)$);
\node[circle,draw=black, fill=red,inner sep=0pt,minimum size=4pt] at ($(1)!0.5!(2)$) {};
\node[circle,draw=black, fill=red,inner sep=0pt,minimum size=4pt] at ($(2)!0.5!(3)$) {};
\node[circle,draw=black, fill=red,inner sep=0pt,minimum size=4pt] at ($(1)!0.5!(3)$) {};
\node[circle,draw=black, fill=red,inner sep=0pt,minimum size=4pt] at (3) {};
\draw[->,thick] (0.5,1.5)--(0.2,1.5);
\node[right] at (0.5,1.5) {\small systematic error};
\end{tikzpicture}
\caption{Systematic error made by taking \texttt{elements(:,3)} as newest nodes.}
\label{fig:systematic}
\end{figure}
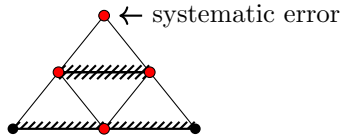 
\item Lines 15--21: Let \texttt{red\_node} be the set of newest nodes of a red pattern. The term \texttt{valence} computed in Line 17 is the \revision{number of elements $\#\tilde{\mathcal{R}}_v$ of the patch $\tilde{\mathcal{R}}_v$} defined in Algorithm~\ref{alg:coarsen}. Nodes of green patterns are then all new nodes that are not red. Note that the array \texttt{green\_node} includes the systematic errors.
\item Lines 22--29: The admissible set $\mathcal{N}_{\mathrm{adm}} =\left\{v\in \mathcal{N}~|~\revision{\#\tilde{\mathcal{R}}}_v \in \left\{2,4\right\} \text{ and } v \in \mathcal{N}_{\mathrm{new}} \text{ and } v \in \mathcal{N}_{\mathrm{mark}}\right\}$ and $\mathcal{N}_{\mathrm{block}}=\mathcal{N}\setminus\mathcal{N}_{\mathrm{adm}}$ are computed. As the reference edge plays a crucial role, we need to do a CLOSURE step to ensure that the shape regularity still holds when coarsening. We again consider the patterns shown in Figure~\ref{fig:TrefineRGBpattern}. We make sure that at least the reference edge of the father element is marked. Differently said, at least the third vertex of the red middle element needs to be blocked if any other vertex in this element is blocked (Lines 27--28). We loop through this process until no further changes are made.
\item Lines 30--37: Hanging nodes are then given by $\mathcal{N}_{\mathrm{hang}} = \mathcal{N}_{\mathrm{block}}\cap \mathcal{N}_{\mathrm{new}}$ (Line 30). With this, we determine the coarsening pattern regarding to these hanging nodes. For green patterns, this is either coarsen or not coarsen. For red patterns, we have more cases to consider, see Figure~\ref{fig:creation}. To this end, we form the weighted sum of hanging nodes. A red pattern can then be coarsened to the patterns: none $(000)$, green $(001)$, blue$_r$ $(101)$, blue$_\ell$ $(011)$ and red $(111)$. In the weighted sum computed in Line 33 these patterns correspond to the values $0,4,5,6$ and $7$. The value $7$ is not considered separately as in this case the elements are kept as they are and are not coarsened. 
\end{itemize}
We omit the presentation of the rest of the code (further 70 Lines), as it is a straightforward implementation of element updates. We distinguish between former red patterns, former green patterns neighbouring a red pattern and former green patterns not neighbouring red patterns. Note, that for the latter the corresponding array includes the systematic error made earlier. To this end, we only consider subsequent elements for coarsening, the case shown in Figure~\ref{fig:bsp_systematic} is out of question. The valence is $\revision{\#\tilde{\mathcal{R}}}_v=2$ but the elements containing $v$ are not numbered consecutively and thus are not considered. In a next step, the old triangulation is deleted. If provided, boundary data is updated. Again, nodes are eliminated only if they are not blocked and they do not stem from the systematic error shown in Figure~\ref{fig:systematic}. Lastly, surplus nodes are deleted and the new coordinates and elements are updated. The interested reader may download the full code from \cite{ameshcoars}.

\hspace*{3cm}
\begin{lstlisting}[xleftmargin=0.8cm,language=matlab,frame=tb,label=lst:coarsen,firstline = 1, lastline = 35,caption= {Lines 1--35 of \texttt{TcoarsenRGB.m}}]
function [coordinates,elements,varargout] = ...
                    TcoarsenRGB(N0,coordinates,elements,varargin)
nC = size(coordinates,1); nE = size(elements,1);
marked = varargin{end};
%*** obtain geometric information
[~,element2edges] = provideGeometricData(elements,zeros(0,4));
edge2elements = createEdge2Elements_adv(element2edges);
%*** determine admissible nodes for coarsening
midelement = find(sum(abs(element2edges(4:end,:)-[element2edges(3:end-1,1),...
   element2edges(1:end-3,2),element2edges(2:end-2,3)]),2)==0)+3;
newest_nodes = zeros(nC,1); marked_nodes = zeros(nC,1);
newest_nodes(elements(:,3)) = 1; newest_nodes(1:N0)=0;
marked_nodes(elements(marked,:)) = 1;
%*** define color of nodes
mid_nodes = elements(midelement,:);
red_node = accumarray(reshape(mid_nodes,[],1),1);
valence = accumarray(elements(:),ones(3*nE,1),[nC,1])...
            -[red_node;zeros(nC-size(red_node,1),1)];
newest_marked = elements(marked,3);
green_node = setdiff(newest_marked(newest_marked>N0),find(red_node));
green_val = valence(green_node);
%*** determine arising hanging nodes
blocked_nodes = ~(((valence == 2) | (valence == 4)) & newest_nodes & marked_nodes);
old = -1;
while ~isempty(find(old ~=blocked_nodes))
    old = blocked_nodes;
    [row,~] = find(reshape(blocked_nodes(mid_nodes),[],3));
    blocked_nodes(elements(midelement(row),3)) = 1; % mark reference edge
end
hanging_nodes = blocked_nodes & newest_nodes;
%*** determine coarsening pattern regarding hanging nodes
red_elements = reshape(hanging_nodes(mid_nodes),[],3);
val = sum(red_elements .* [1,2,4],2);
none = find(val==0);
gdx = find(val == 4);
\end{lstlisting}

\subsection{A Minimal Example}\label{sect:minimalexample}

Listing~\ref{lst:minimalexample} shows an exemplary code of how to embed the coarsening routine into a framework. We start with defining an initial mesh $\mathcal{T}_0$ (Lines 1--4). A refined mesh $\tilde{\mathcal{T}}$ is created via \texttt{TrefineRGB} (Lines 7--15). For a given triangulation $\mathcal{T}$ and a given discrete point set $\mathcal{P}$, the function \texttt{point2element} determines the elements of $\mathcal{T}$ that include $p$ for some $p \in \mathcal{P}$. Thus, for the defined discrete point set in Lines 16--22, elements in $\tilde{\mathcal{T}}$ are marked according to \texttt{point2element}. We coarsen the mesh via the function call \texttt{TcoarsenRGB} (Line 26--27) until no further change is made (Line 28). Lines 33--35 plot the locally coarsened mesh.

\begin{lstlisting}[xleftmargin=0.8cm,language=matlab,frame=tb,label=lst:minimalexample,caption= {A minimal example}]
%% Define initial mesh
coordinates = [0,0;1,0;1,1;0,1;2,0;2,1];
elements = [3,1,2;1,3,4;2,6,3;6,2,5];
boundary = [1,2;2,5;5,6;6,3;3,4;4,1];
N0 = size(coordinates,1);
c_old = 0;
%% Refine uniformly
while 1
    marked = 1:size(elements,1);
    [coordinates,elements,boundary] ...
        = TrefineRGB(coordinates,elements,boundary,marked);
    if isempty(marked)|| (size(coordinates,1)>1e3)
        break
    end
end
% define discrete points (here a disc)
phi = -pi:pi/50:pi;
r = 0.2:1/50:0.4;
[r,phi] = meshgrid(r,phi);
s = r.*cos(phi);
t = r.*sin(phi);
points = [s(:)+1,t(:)+0.5];
%% Coarsen at discrete points
while 1
    mark3 = point2element(coordinates,elements,points);
    [coordinates,elements,boundary] = ...
                TcoarsenRGB(N0,coordinates,elements,boundary,mark3);
    if size(coordinates,1) == c_old
        break
    end
    c_old = size(coordinates,1);
end
% plot results
clf, patch('Faces',elements,'Vertices',coordinates,'Facecolor','none')
axis equal, axis off

\end{lstlisting}

\subsection{Examples and Demo Files}
The coarsening routine for RGB meshes is part of the toolbox \texttt{ameshcoars} - Efficient Implementation of Adaptive Mesh Coarsening in 2D \cite{ameshcoars}. Numerical examples and demo files based on the interplay of refinement and coarsening are provided in subdirectories of the \texttt{ameshcoars}--toolbox:
\begin{itemize}
\item \texttt{example1/}: refinement along a moving circle,
\item \texttt{example2/}: adaptive finite element implementation following \cite{p1afem} for a quasi-stationary partial differential equation,
\item \texttt{example3/}: triangulation of a GIF, 
\item \texttt{example4/}: local coarsening of a uniformly refined triangulation.
\end{itemize}

\section{Numerical Experiments}\label{sect:experiments}

In this section, we test our coarsening routine with MATLAB 2018a. We present some results based on \texttt{example1/} and \texttt{example4/} of the \texttt{ameshcoars}-toolbox \cite{ameshcoars}. In particular, we look at the interplay of refinement and coarsening and how well moving singularities can be captured by this procedure. Further, we take a look at the efficiency of our coarsening algorithm. We know that our coarsening implementation is not inverse to the refinement but that it can fully recover the initial triangulation. To this end, we want to examine what element- and coordinate-ratios we get between each refinement/coarsening step to get a feeling of how efficient our coarsening routine is. \revision{Further, we examine our implementation for scalability and give a remark on how the implementation depends on the local refinement.} Lastly, we show that coarsening can be done locally.

\subsection{Interplay of Refinement and Coarsening}
\begin{figure}[!htb]
\centering
\begin{tikzpicture}
    \draw [<->,thick] (0,3) node (yaxis) [above] {\scriptsize \#$\mathcal{N}$}
        |- (3,0) node (xaxis) [right] {$t$};
\node[circle,draw=black, fill=white,inner sep=0pt,minimum size=4pt] at (1,0.3) {};  
\node[circle,draw=black, fill=white,inner sep=0pt,minimum size=4pt] at (1,0.6) {};  
\node[circle,draw=black, fill=white,inner sep=0pt,minimum size=4pt] at (1,0.9) {};  
\node[circle,draw=black, fill=white,inner sep=0pt,minimum size=4pt] at (1,1.4) {};  
\node[circle,draw=black, fill=white,inner sep=0pt,minimum size=4pt] at (1,2.0) {};  
\node[circle,draw=black, fill=white,inner sep=0pt,minimum size=4pt] at (1,2.7) {};  
\node[circle,draw=black, fill=white,inner sep=0pt,minimum size=4pt] at (2,2.1) {};
\node[circle,draw=black, fill=white,inner sep=0pt,minimum size=4pt] at (2,2.5) {};  
\node[circle,draw=black, fill=white,inner sep=0pt,minimum size=4pt] at (2,2.9) {}; 
\draw[->,thick] (1,2.7)--(2,2.1);
\node[circle,draw=black, fill=white,inner sep=0pt,minimum size=4pt] at (3.1,2.7) {};
\draw[->,thick] (3,2.5)--(3.3,2.3);
\node[right] at (3.3,2.7) {\scriptsize refinement steps};  
\node[right] at (3.3,2.4) {\scriptsize coarsening steps};  
\draw (2.8,2.2) rectangle (5.8,2.9);
\draw (1,-0.1)--(1,0.1);
\node[below] at (1,0) {$t_0$};
\draw (2,-0.1)--(2,0.1);
\draw[dashed] (0,2.2)--(2.5,2.2);
\draw[dashed] (0,2.65)--(2.5,2.65);
\node[left] at (0,2.2) {\scriptsize $N_{\mathrm{min}}$};
\node[left] at (0,2.65) {\scriptsize $N_{\mathrm{max}}$};
\node[below] at (2,0) {$t_1$};
\node at (1.5,-0.8) {refinement and coarsening};
    \draw [<->,thick] (7,3) node (yaxis) [above] {\scriptsize \#$\mathcal{N}$}
        |- (10,0) node (xaxis) [right] {$t$};
\node[circle,draw=black, fill=white,inner sep=0pt,minimum size=4pt] at (8,0.3) {};  
\node[circle,draw=black, fill=white,inner sep=0pt,minimum size=4pt] at (8,0.6) {};  
\node[circle,draw=black, fill=white,inner sep=0pt,minimum size=4pt] at (8,0.9) {};  
\node[circle,draw=black, fill=white,inner sep=0pt,minimum size=4pt] at (8,1.4) {};  
\node[circle,draw=black, fill=white,inner sep=0pt,minimum size=4pt] at (8,2.0) {};  
\node[circle,draw=black, fill=white,inner sep=0pt,minimum size=4pt] at (8,2.7) {};  
\node[circle,draw=black, fill=white,inner sep=0pt,minimum size=4pt] at (9,0.3) {};  
\node[circle,draw=black, fill=white,inner sep=0pt,minimum size=4pt] at (9,0.7) {};  
\node[circle,draw=black, fill=white,inner sep=0pt,minimum size=4pt] at (9,0.95) {};  
\node[circle,draw=black, fill=white,inner sep=0pt,minimum size=4pt] at (9,1.3) {};  
\node[circle,draw=black, fill=white,inner sep=0pt,minimum size=4pt] at (9,1.9) {};  
\node[circle,draw=black, fill=white,inner sep=0pt,minimum size=4pt] at (9,2.3) {}; 
\node[circle,draw=black, fill=white,inner sep=0pt,minimum size=4pt] at (9,3.0) {};
\node[right] at (10.3,2.7) {\scriptsize refinement steps};  
\node[circle,draw=black, fill=white,inner sep=0pt,minimum size=4pt] at (10.1,2.7) {};
\draw (9.8,2.5) rectangle (12.8,2.9);
\draw (8,-0.1)--(8,0.1);
\draw[dashed] (7,2.65)--(9.5,2.65);
\node[left] at (7,2.65) {\scriptsize $N_{\mathrm{max}}$};
\node[below] at (8,0) {$t_0$};
\draw (9,-0.1)--(9,0.1);
\node[below] at (9,0) {$t_1$};
\node at (8.5,-0.8) {refinement only};
\end{tikzpicture}
\caption{Left: Interplay between refinement and coarsening to capture moving circle in time step $t_1>t_0$. Right: Refinement only to capture moving circle in time step $t_1$.}
\label{fig:scheme}
\end{figure}
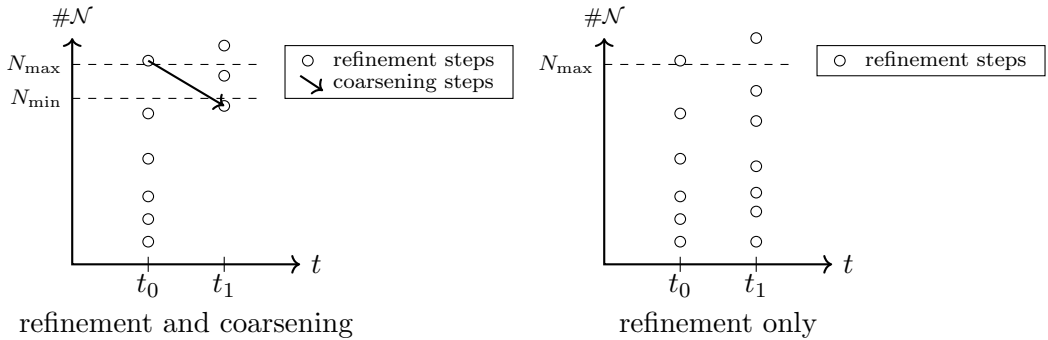

Let us start with a basic example. Adaptive coarsening is widely used to release degrees of freedom that are not needed anymore as, e\,.g.\,, a singularity advances. We imitate the behavior by a moving circle, which is supposed to represent the singularity. Starting off with an initial triangulation $(\mathcal{T}_0,\refe_{\mathcal{T}_0})$, we refine along the circle at time $t_0$. To capture the singularity in $t_1>t_0$, we could start off again from $(\mathcal{T}_0,\refe_{\mathcal{T}_0})$ and only use the refinement procedure. However, as the circle progresses steadily, only a few nodes need to be released and only a few nodes need to be added. To this end, we use our coarsening routine to set the corresponding coordinates free and the refinement routine to add further coordinates to capture the shape of the circle. The comparison of number of degrees of freedom between refinement only and refinement combined with coarsening is illustrated in Figure~\ref{fig:scheme}.

\begin{figure}[!htb]
\begin{minipage}{0.45\textwidth}
\includegraphics[width=\textwidth]{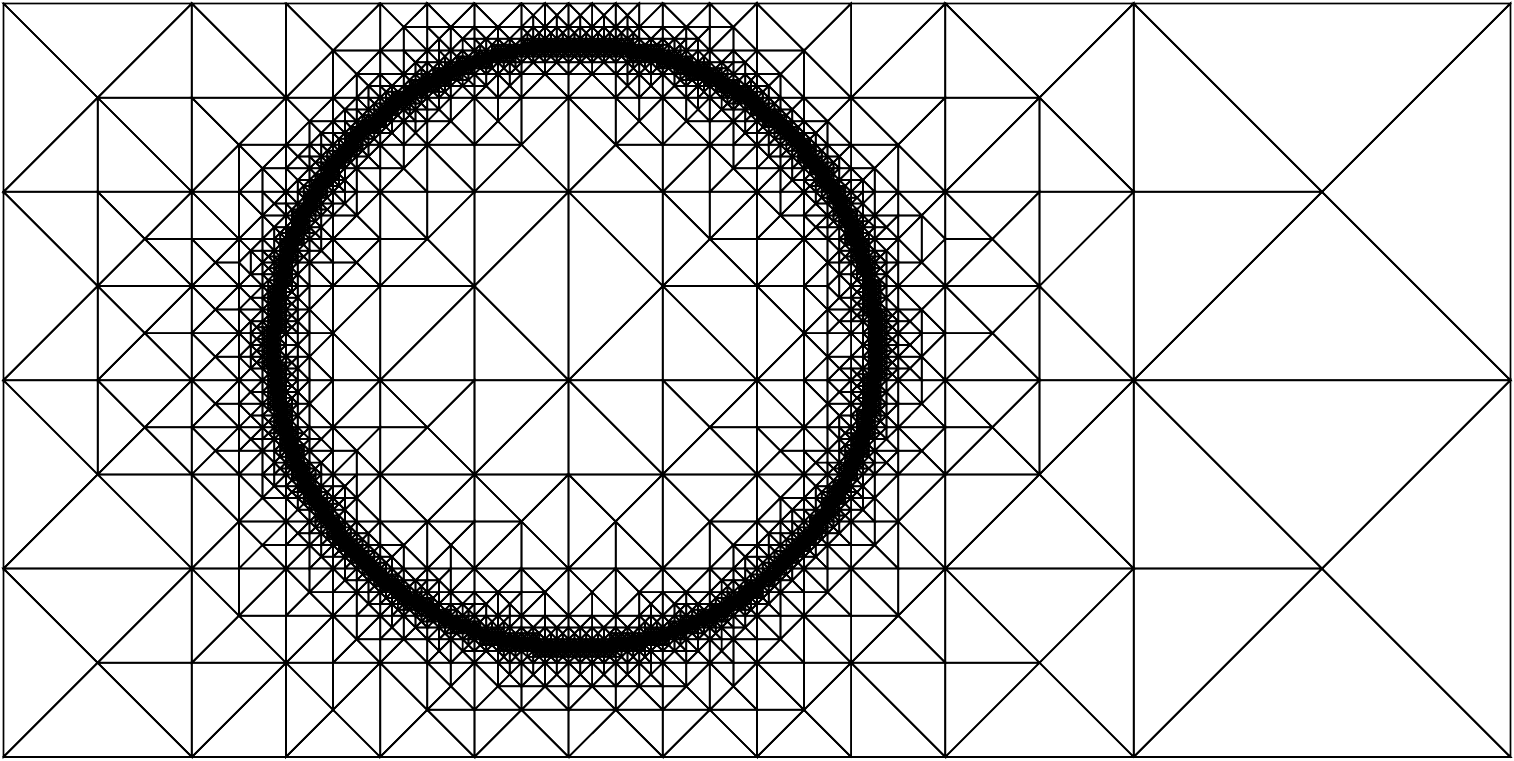}
\end{minipage}\hfill
\begin{minipage}{0.45\textwidth}
\includegraphics[width=\textwidth]{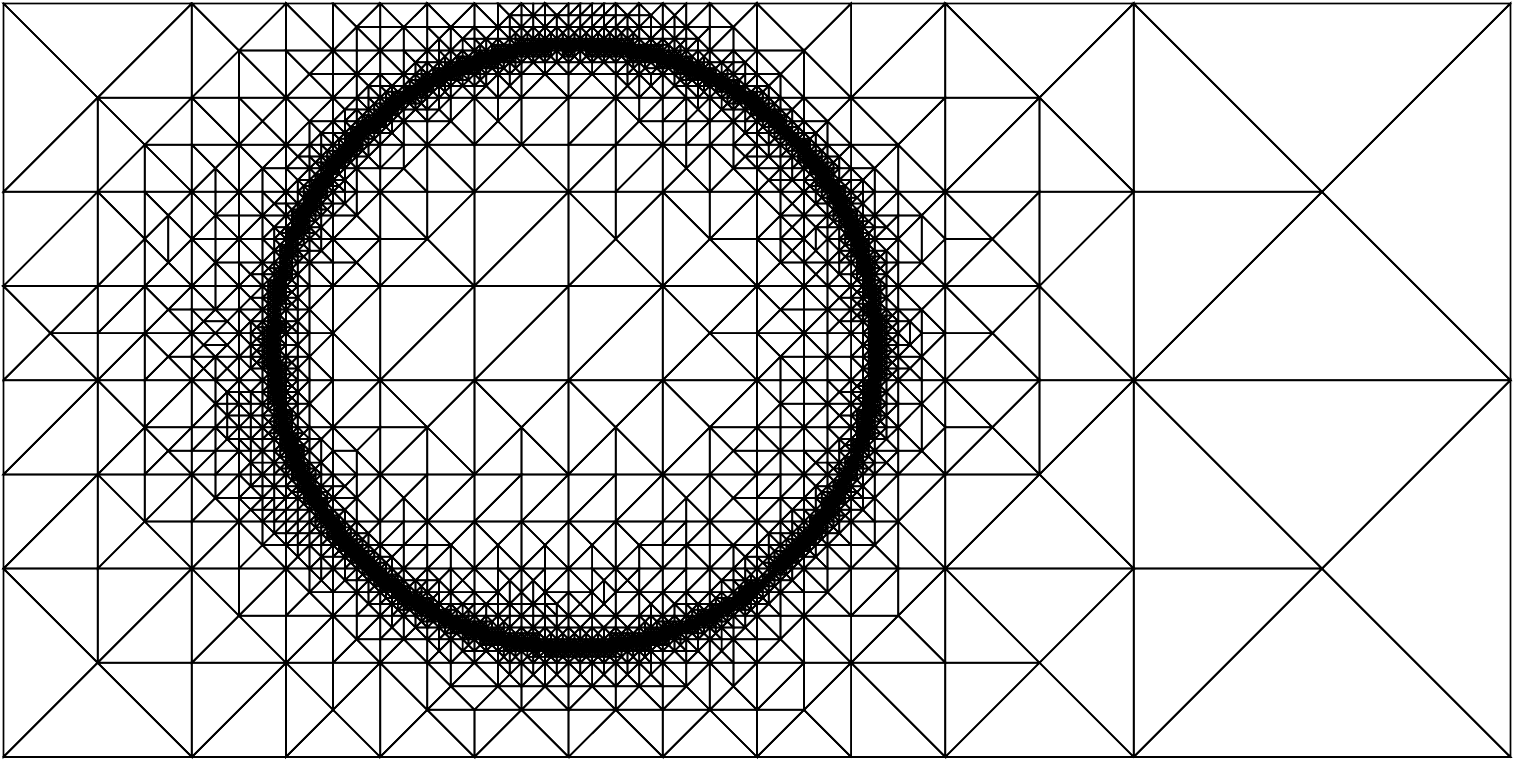}
\end{minipage}\\
\begin{minipage}{0.45\textwidth}
\vspace*{2mm}
\centering $t_1$: center at $(0.759\,,0.545)$
\vspace*{2mm}
\end{minipage}\hfill
\begin{minipage}{0.45\textwidth}
\vspace*{2mm}
\centering $t_1$: center at $(0.759\,,0.545)$
\vspace*{2mm}
\end{minipage}\\
\begin{minipage}{0.45\textwidth}
\includegraphics[width=\textwidth]{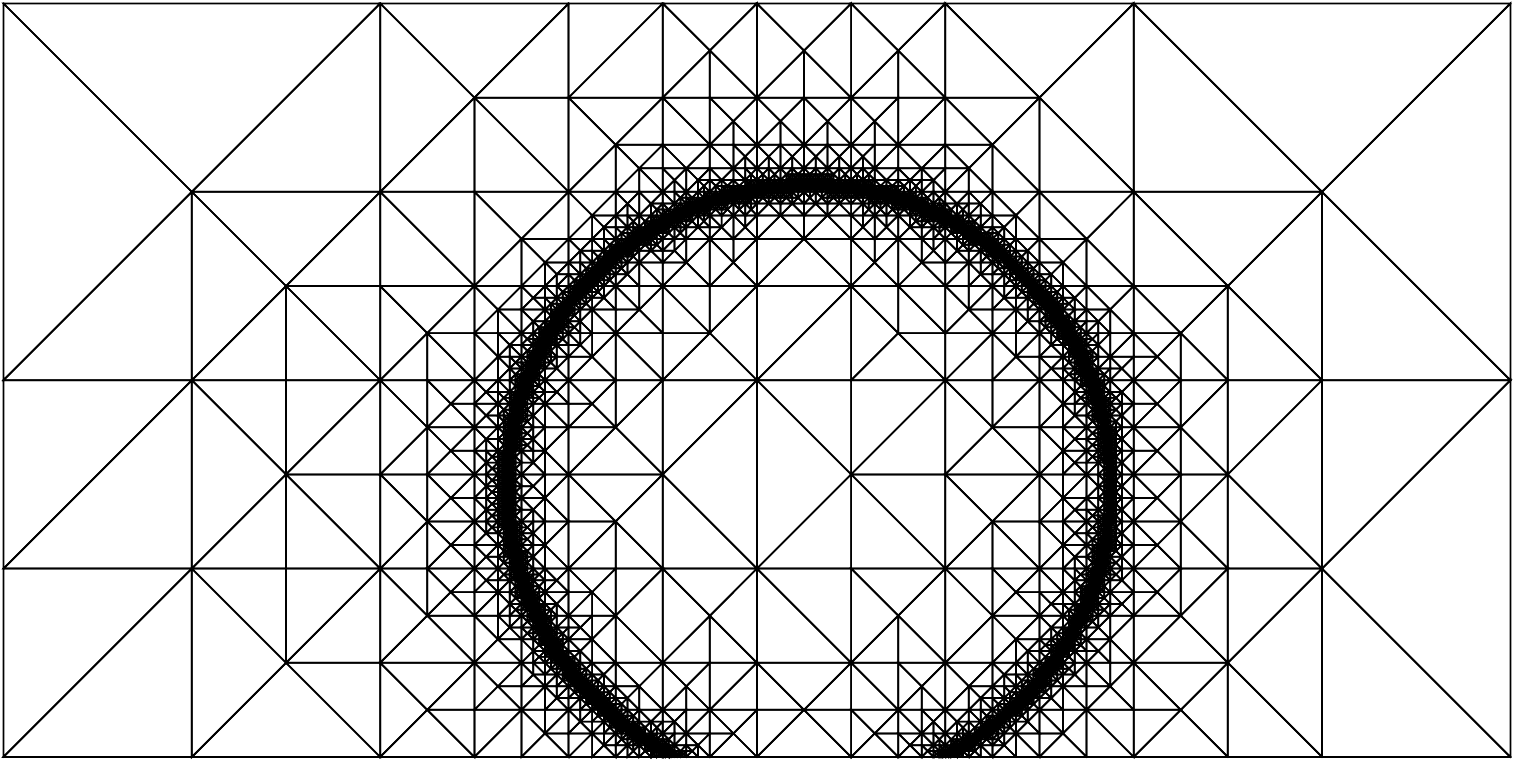}
\end{minipage}\hfill
\begin{minipage}{0.45\textwidth}
\includegraphics[width=\textwidth]{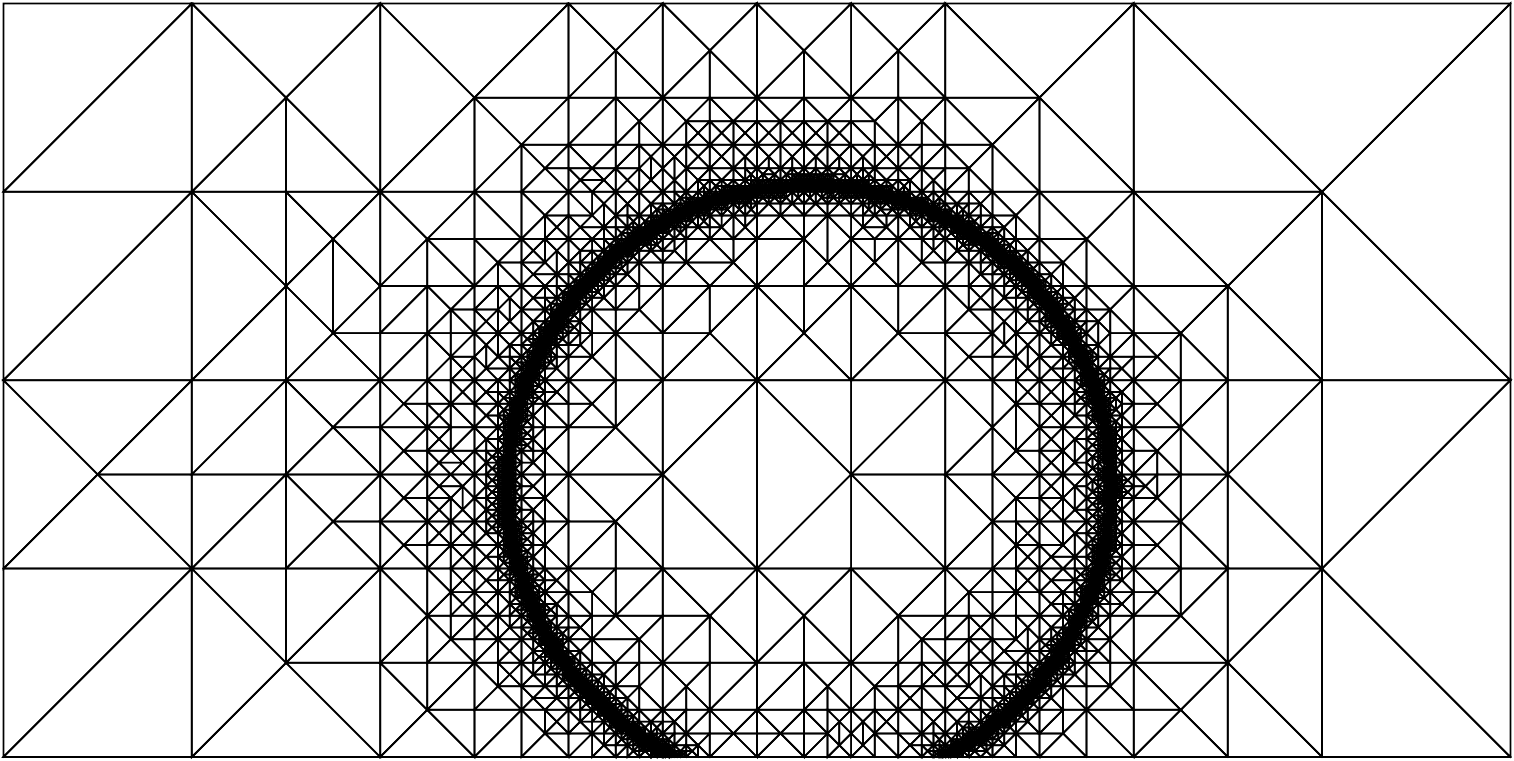}
\end{minipage}\\
\begin{minipage}{0.45\textwidth}
\vspace*{2mm}
\centering $t_2$: center at $(1.069\,,0.359)$
\vspace*{2mm}
\end{minipage}\hfill
\begin{minipage}{0.45\textwidth}
\vspace*{2mm}
\centering $t_2$: center at $(1.069\,,0.359)$
\vspace*{2mm}
\end{minipage}\\
\begin{minipage}{0.45\textwidth}
\includegraphics[width=\textwidth]{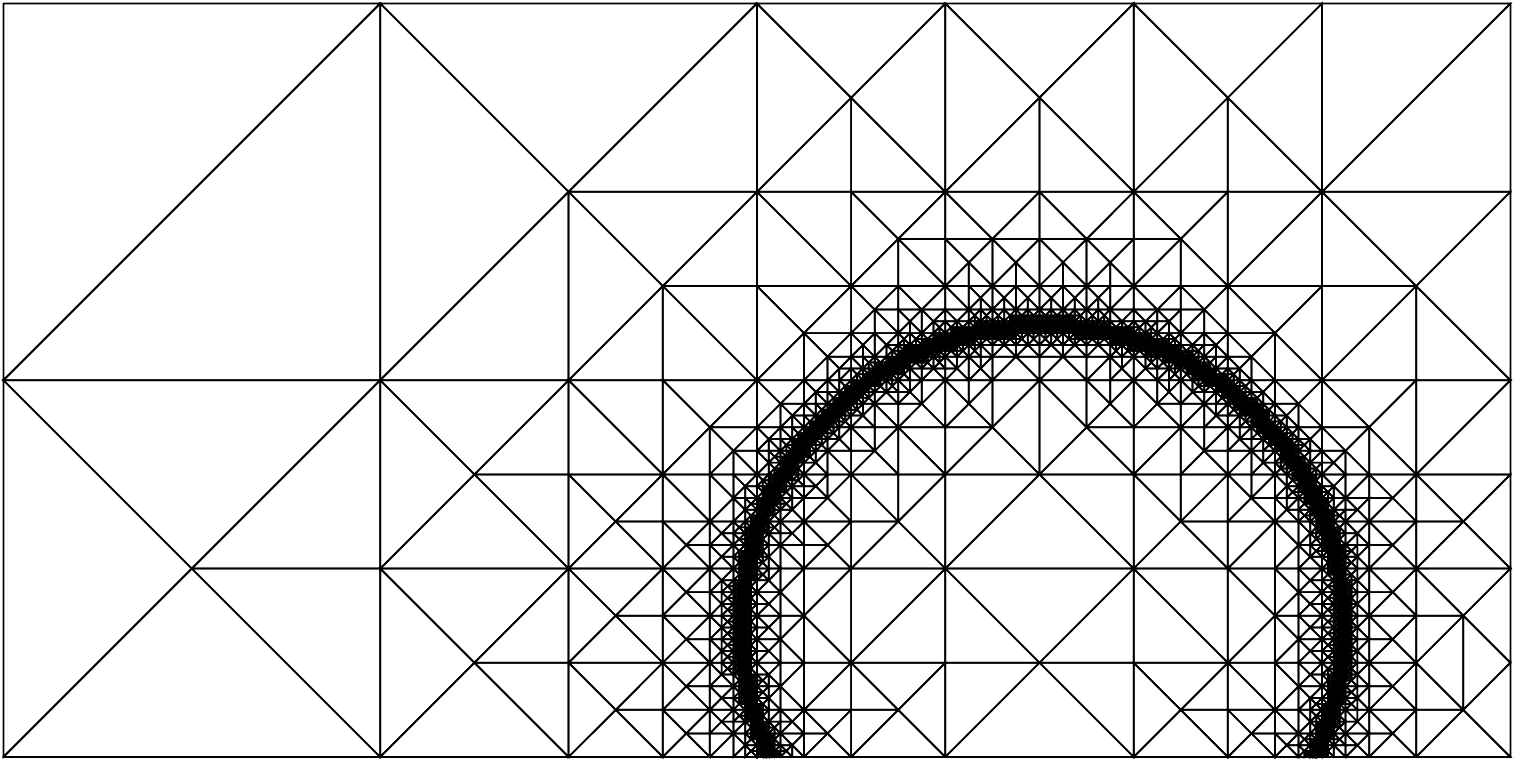}
\end{minipage}\hfill
\begin{minipage}{0.45\textwidth}
\includegraphics[width=\textwidth]{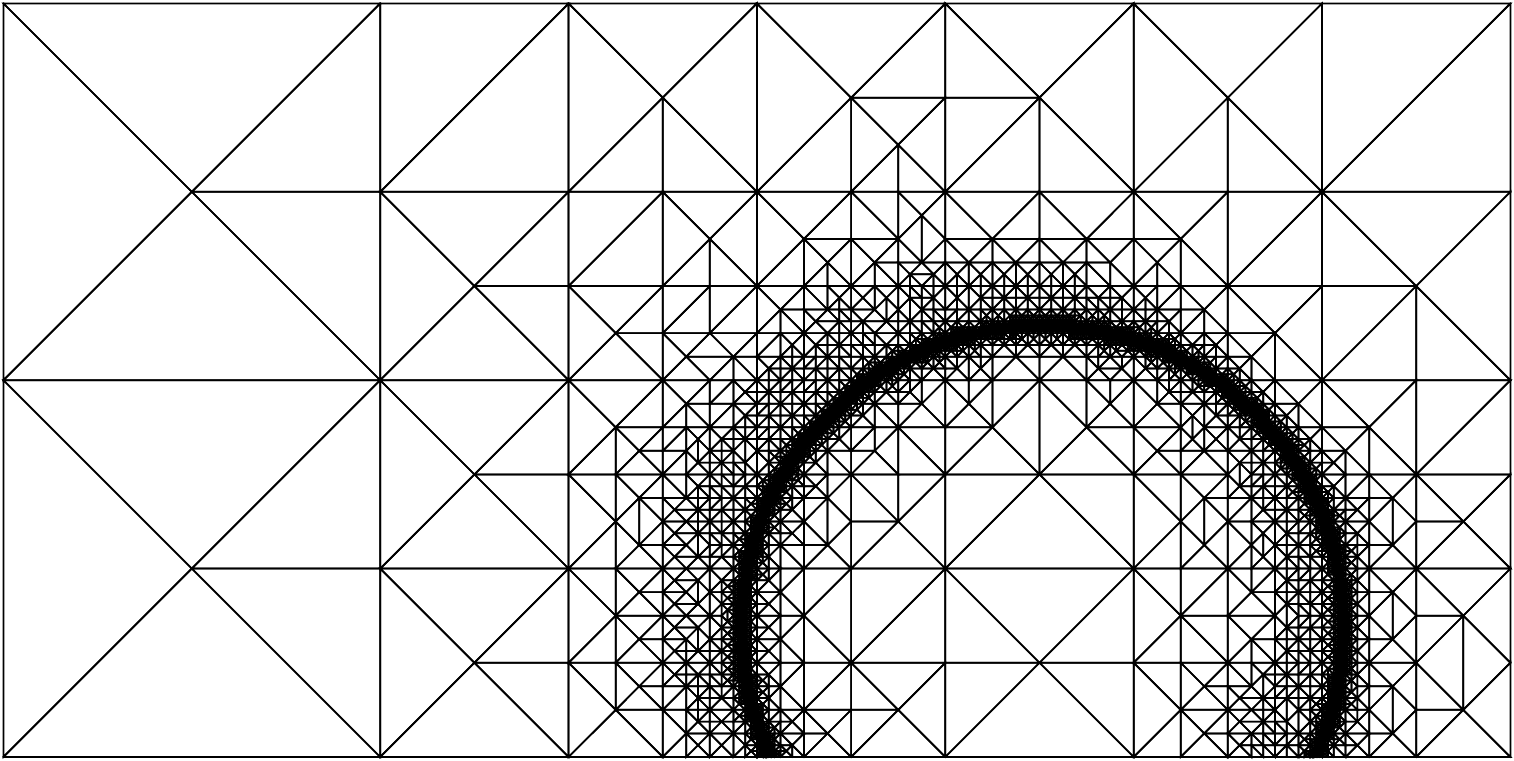}
\end{minipage}\\
\begin{minipage}{0.45\textwidth}
\vspace*{2mm}
\centering $t_3$: center at $(1.379\,,0.172)$
\vspace*{2mm}
\end{minipage}\hfill
\begin{minipage}{0.45\textwidth}
\vspace*{2mm}
\centering $t_3$: center at $(1.379\,,0.172)$
\vspace*{2mm}
\end{minipage}\\
\caption{Moving circle: triangulations at different time frames $t_1,t_2,t_3$ with $N_{\mathrm{max}} = 10^4$. Left: $N_{\mathrm{min}}=10^2$. Right: $N_{\mathrm{min}}=10^3$. On the right, we observe a pollution effect; the circle draws a tail. }
\label{fig:movingcircle}
\end{figure}

We first explain the procedure in this example. We define a circle with center and radius. If this circle intersects an element we mark this element for refinement. This is done by the function \texttt{markcircle} and together with the refinement routine this ensures that we capture the shape of the circle. We do this until a maximal number $N_{\mathrm{max}}$ of coordinates is reached, or the element size becomes too small. We then mark all elements for coarsening, coarsen and repeat this until we fall below a given minimal number $N_{\mathrm{min}}$ of nodes.  Subsequently, in the next time step $t_1$, we again refine along the circle (now moved to another position!), but do not start from the initial mesh but from the coarsened mesh from time step $t_0$. This is done consecutively, and we get a sequence of triangulations capturing the moving circle, see left column of Figure~\ref{fig:movingcircle} and cf.~\texttt{example1/} in \cite{ameshcoars}.

Note, that the triangulations we get are highly sensitive to the choice of the parameter $N_{\mathrm{min}}$. We observe a pollution effect if $N_{\mathrm{min}}$ is chosen too big, see right column of Figure~\ref{fig:movingcircle}. However, in general, a pollution effect does not falsify the computation. It only means that the shape of the circle is not captured in the best possible way and nodes exist where they are not necessarily needed. In principal, it is not just the parameter $ N_{\mathrm{min}}$ that is responsible for a pollution effect. It also depends on how fast the front for refinement advances, how many time steps are considered, how $N_{\mathrm{max}}$ is chosen, etc.\,. However, in most applications the pollution effect is controlled as error estimators are often used to regulate the error and thus also adapt the mesh, cf. \texttt{example2/} and \texttt{example3/} in \cite{ameshref}.

\subsection{Efficiency of the Coarsening Routine}\label{sect:MN}

We want to determine how efficient our coarsening routine is in the sense of element- and coordinate-ratios in between two coarsening steps in comparison to two refinement steps. Let therefore $\mathcal{T}_i$ be the triangulation after $i$ refinement steps and $\mathcal{N}_i$ denotes the set of nodes of the triangulation $\mathcal{T}_i$. We determine the element-ratio \[\rho_{\mathrm{elem}}^i = \frac{\# \mathcal{T}_{i+1}}{\# \mathcal{T}_i}\] and the coordinate-ratio \[\rho_{\mathrm{coord}}^i= \frac{\#\mathcal{N}_{i+1}}{\#\mathcal{N}_i}\] in each refinement step. Let $N$ be the number refinement steps. We expect $1 \leq \rho_{\mathrm{elem}}^i\leq 4$ for all $i = 1,\ldots, N$ due to the refinement patterns. We compute the actual ratios for an adaptive refinement along a circle. They  are presented in Table~\ref{tab:ratios_ref}. We see that the geometric means are $\overline{\rho}_{\mathrm{elem}} = 2.40$ and $\overline{\rho}_{\mathrm{coord}} = 2.15$. 

Analogously, we determine the ratios for coarsening steps. To this end, let $\hat{\mathcal{T}}_j$ be the triangulation received after $j$ coarsening steps, $j\in \left\{1,\ldots,M\right\}$, and $\hat{\mathcal{N}}_j$ is the corresponding set of nodes. Note, that $\mathcal{T}_N = \hat{\mathcal{T}}_0$, i.\,e.\,, we start our coarsening routine from the finest mesh and mark all elements for coarsening. As our coarsening routine is not inverse we need more coarsening steps than refinement steps to recover the initial triangulation. In other words, $M\geq N$. A blue refinement is coarsened in a two-step procedure. Thus, we expect a refinement--coarsening step ratio of about $1:2$ and consequently $M\approx 2\cdot N$. For a better quantification, we compute the element-ratio
\[ \hat{\rho}_{\mathrm{elem}}^j = \frac{\# \hat{\mathcal{T}}_{j}}{\# \hat{\mathcal{T}}_{j+1}} \] and the coordinate-ratio \[\hat{\rho}_{\mathrm{coord}}^j= \frac{\#\hat{\mathcal{N}}_{j}}{\#\hat{\mathcal{N}}_{j+1}}\] in each coarsening step $j$. The results are shown in Table~\ref{tab:ratios_coars}. The geometric means are given by $\overline{\hat{\rho}}_{\mathrm{elem}} =1.55$ and $\overline{\hat{\rho}}_{\mathrm{coord}} =1.47$. In this example we get $M=2\cdot N$ and in other experiments we also observed $M\approx 2\cdot N$. In terms of efficiency, this means that our coarsening strategy is not as efficient as its refinement counterpart. We need to expect twice as many coarsening steps to undo the refinement as is needed for refining. To get a better feeling about time efficiency, we measured the time for the refinement and coarsening part for this example. The computational time for the refinement part is 0.0665 seconds while the coarsening part takes 0.2171 seconds. Since twice as many coarsening steps are required than refinement steps, we conclude that one coarsening step is slightly more time-consuming than one refinement step. However, coarsening from $10^4$ degrees of freedom to $6$ can be done in about a fifth of a second, which is still very efficient.  \\

\begin{minipage}{0.48\textwidth}
\centering
\captionof{table}{Element- and coordinate-ratios in between two refinement steps.}
\label{tab:ratios_ref}
\begin{tabular}{rrrrr}
 \multicolumn{1}{c}{\boldmath $i$} &\multicolumn{1}{c}{ \boldmath  $\#\mathcal{T}_i$} &\multicolumn{1}{c}{ \boldmath $\#\mathcal{N}_i$} & \multicolumn{1}{c}{ \boldmath$\rho_{\mathrm{elem}}^i$} & \multicolumn{1}{c}{ \boldmath$\rho_{\mathrm{coord}}^i$}\\ \hline
0 & 4 & 6 & - & -\\
1 & 13 & 12 & 3.25 & 2.00\\
2 & 39 & 28& 2.00 & 2.33\\
3 & 123 & 74 & 3.15 & 2.64\\
4 & 297 & 164 & 2.41 & 2.22\\
5 & 693 & 365 & 2.33 & 2.23\\
6 & 1482 & 762 & 2.14 & 2.09\\
7 & 3085 & 1568 & 2.08 & 2.06\\
8 & 6239 & 3147 & 2.02 & 2.01\\
9 & 12597 & 6328 & 2.02 & 2.01 \\
10 & 25221 & 12642 & 2.00 & 2.00\\
&&&&\\
&&&&\\
&&&&\\
&&&&\\
&&&&\\
&&&&\\
&&&&\\
&&&&\\
&&&&\\
&&&&\\
\end{tabular}
\end{minipage}\hfill
\begin{minipage}{0.48\textwidth}
\centering
\captionof{table}{Element- and coordinate-ratios in between two coarsening steps.}
\label{tab:ratios_coars}
\begin{tabular}{rrrrr}
\multicolumn{1}{c}{ \boldmath $j$} & \multicolumn{1}{c}{ \boldmath $\#\hat{\mathcal{T}}_j$ }& \multicolumn{1}{c}{ \boldmath$\#\hat{\mathcal{N}}_j$} & \multicolumn{1}{c}{ \boldmath$\hat{\rho}_{\mathrm{elem}}^j$} &\multicolumn{1}{c}{ \boldmath $\hat{\rho}_{\mathrm{coord}}^j$}\\ \hline
0 & 25221 & 12642 & - & - \\
1 & 16610 & 8335 & 1.52 & 1.52\\
2 & 13454 & 6756 & 1.23 & 1.23\\
3 & 8851 & 4453 & 1.52 & 1.52 \\
4 & 6956 & 3505 & 1.27 & 1.27 \\
5 & 4484 & 2268 & 1.55 & 1.55 \\
6 & 3485 & 1768 & 1.29 & 1.28 \\
7 & 2199 & 1123 & 1.58 & 1.57 \\
8 & 1684 & 865 & 1.31 & 1.30 \\
9 & 1052 & 547 & 1.60 & 1.58 \\
10 & 800 & 421 & 1.32 & 1.30 \\
11 & 486 & 261 & 1.65 & 1.61 \\
12 & 360 & 198 & 1.35 & 1.32 \\
13 & 203 & 115 & 1.77 & 1.72 \\
14 & 143 & 85 & 1.42 & 1.35\\
15 & 70 & 45 & 2.04 & 1.89 \\
16 & 48 & 34 & 1.46 & 1.32 \\
17 & 19 & 16 & 2.53 & 2.13 \\
18 & 12 & 11 & 1.58 & 1.45 \\
19 & 6 & 7 & 2.00 & 1.57 \\
20 & 4 & 6 & 1.50 & 1.17\\
\end{tabular}
\end{minipage}

\revision{\subsection{Scalability of the Coarsening Routine}
This leads us to the question of scalability. To this end, we measure 20 times the computational time of the coarsening routine by use of MATLAB's builtin \matlab{tic/toc} and plot the mean of the measured times above the number of nodes  for newest vertex bisection and red-green-blue refinement. The numerical experiments on NVB are based on the implementation in \cite{p1afem}. This is diplayed in Figure~\ref{fig:scale}. The plot shows an almost linear behavior between the number of elements and the computational time in seconds. We see that RGB has some offset which is explained by a more involved determination of red middle elements and the fact that a CLOSURE step is carried out within a while-loop. A priori, it is unclear how often the while-loop iterates. This uncertainty in the CLOSURE step was already examined for the refinement routine in \cite{funkenschmidt}. Due to the structure of the mesh, for an adaptive RGB refined mesh of an initial weak BDD triangulation, at most two iterations are needed. The CLOSURE step is thus predictable and can not cause a huge increase in computational time.
%
%
\definecolor{mycolor1}{rgb}{0.00000,0.44700,0.74100}%
\definecolor{mycolor2}{rgb}{0.85000,0.32500,0.09800}%
\begin{figure}
\centering
\begin{tikzpicture}

\begin{axis}[%
width=0.75\textwidth,
height=0.25\textheight,
at={(2.401in,1.28in)},
scale only axis,
xmode=log,
xmin=7,
xmax=1000000,
xminorticks=true,
xlabel style={font=\color{white!15!black}},
xlabel={number of nodes},
ymode=log,
ymin=0.0001,
ymax=10,
yminorticks=true,
ylabel style={font=\color{white!15!black}},
ylabel={computational time [s]},
axis background/.style={fill=white},
legend style={at={(0.1,0.704)},anchor=south west, legend cell align=left, align=left, draw=white!15!black}
]
\addplot [color=mycolor1, mark=triangle, mark options={solid, mycolor1}]
  table[row sep=crcr]{%
811167	4.9065489708\\
535624	3.14202917735\\
434375	2.6166304213\\
286376	1.6039741625\\
225324	1.29531917935\\
146476	0.78578596405\\
114649	0.6264443756\\
74282	0.38352574025\\
57773	0.2922663673\\
37338	0.17332798055\\
28982	0.13593288475\\
18716	0.0813490177\\
14543	0.0625928281\\
9393	0.03675581855\\
7276	0.02885393165\\
4694	0.01735467345\\
3629	0.01350037815\\
2327	0.00980127685\\
1796	0.0078617738\\
1137	0.00501747135\\
875	0.0042500157\\
553	0.0029633933\\
424	0.0024501759\\
263	0.00175175995\\
198	0.0015417947\\
115	0.0014166153\\
85	0.0011681161\\
45	0.00125064595\\
34	0.0009112769\\
16	0.0009161894\\
11	0.00100413925\\
7	0.00107439805\\
};
\addlegendentry{RGB}

\addplot [color=mycolor2, mark=o, mark options={solid, mycolor2}]
  table[row sep=crcr]{%
835135	2.9148303216\\
640556	2.13493937855\\
470138	1.5734318517\\
332364	1.0846874101\\
235076	0.7593663075\\
166190	0.5195997892\\
117548	0.36315741705\\
83104	0.2438660294\\
58780	0.1695051417\\
41554	0.114864975\\
29390	0.08203458315\\
20772	0.05435086115\\
14690	0.0371924118\\
10380	0.0252625229\\
7338	0.01748132315\\
5182	0.01176272185\\
3658	0.00814031715\\
2576	0.0053987895\\
1812	0.00377355695\\
1266	0.0025010359\\
884	0.00175936695\\
614	0.0012445871\\
426	0.0009229767\\
292	0.00068318685\\
198	0.0005145303\\
130	0.0003868203\\
87	0.00031890985\\
53	0.0002714325\\
34	0.0002344383\\
20	0.00020336695\\
14	0.00019078055\\
8	0.00018039365\\
};
\addlegendentry{NVB}

\addplot [color=black]
  table[row sep=crcr]{%
50000	0.1\\
500000	1\\
};

\addplot [color=black]
  table[row sep=crcr]{%
50000	0.1\\
500000	0.1\\
};

\addplot [color=black]
  table[row sep=crcr]{%
500000	0.1\\
500000	1\\
};

\addplot[area legend, draw=black, fill=white!80!black, fill opacity=0.4]
table[row sep=crcr] {%
x	y\\
50000	0.1\\
500000	0.1\\
500000	1\\
50000	0.1\\
}--cycle;

\node[below, align=center]
at (axis cs:158113.91,0.093) {1};
\node[right, align=left]
at (axis cs:550644.147,0.316) {-1};
\end{axis}
\end{tikzpicture}%
\caption{\revision{Scalability of our RGB coarsening algorithm in comparison to the NVB coarsening algorithm implemented in \cite{p1afem}. A nearly linear behavior between the number of nodes and the computational time can be observed. }}
\label{fig:scale}
\end{figure}
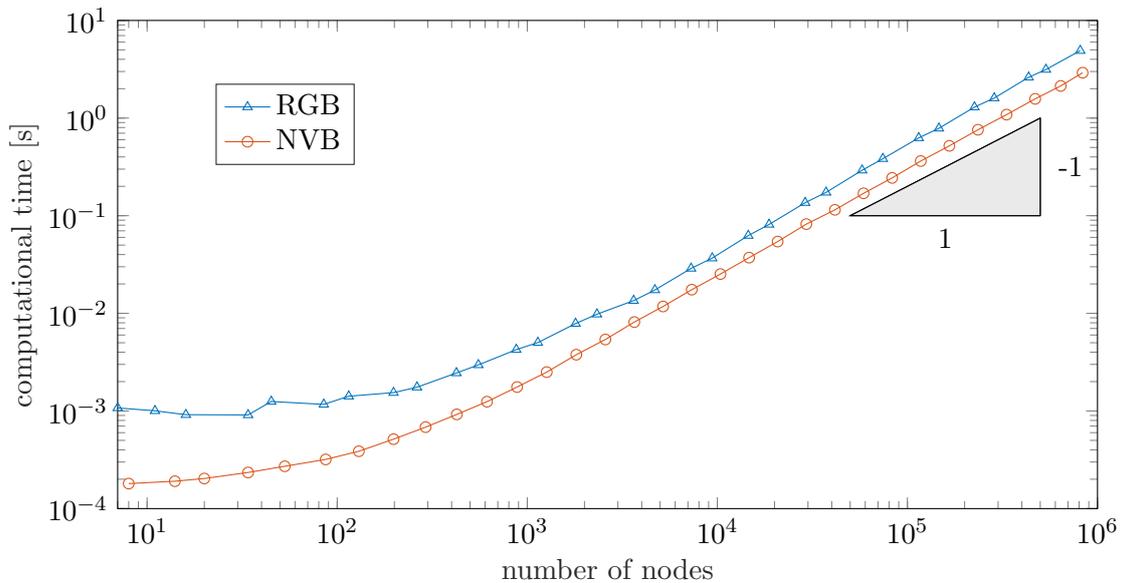}

\subsection{Local Coarsening}
\begin{figure}[!htb]
\begin{minipage}{0.45\textwidth}
\includegraphics[width=\textwidth]{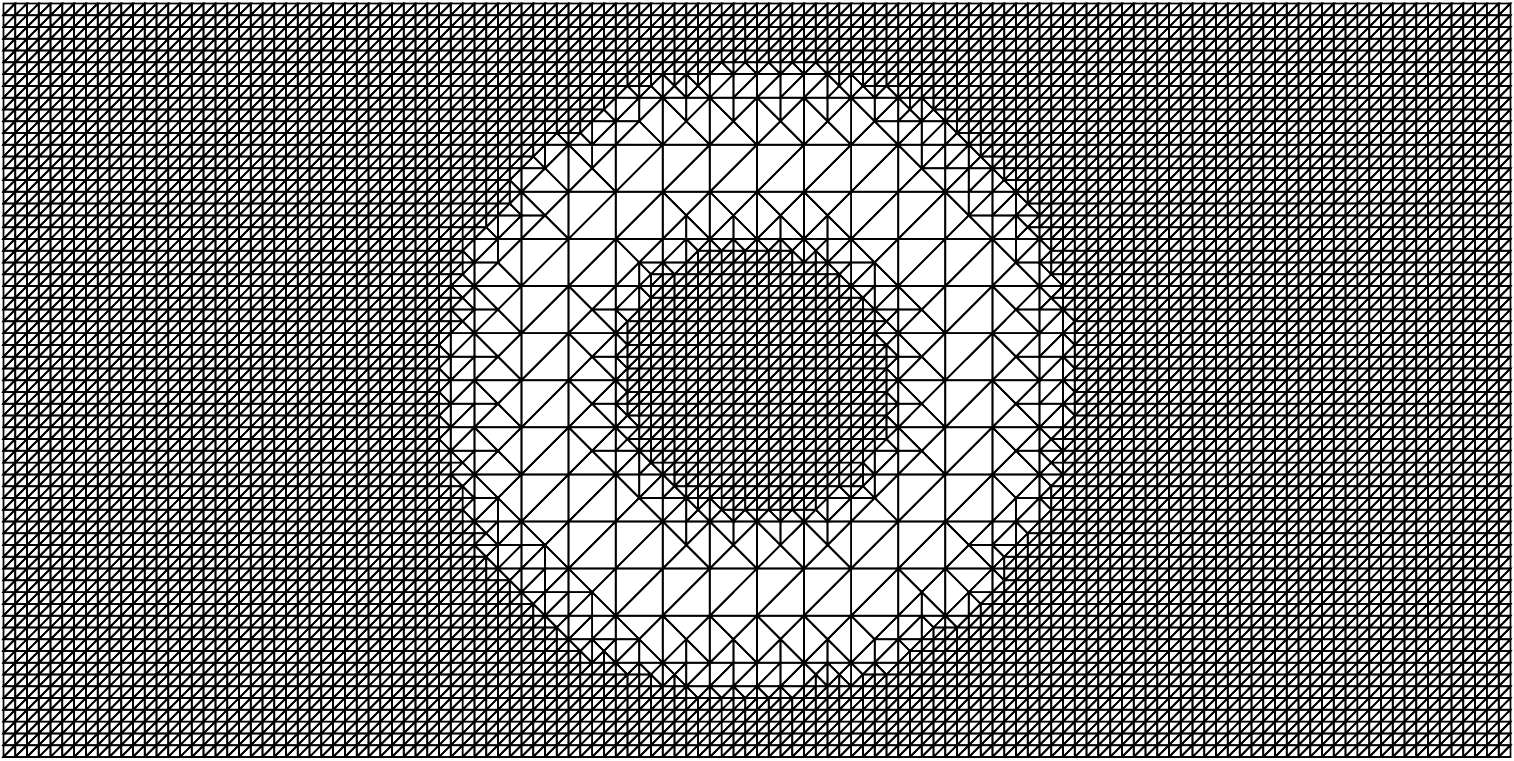}
\end{minipage}\hfill
\begin{minipage}{0.45\textwidth}
\includegraphics[width=\textwidth]{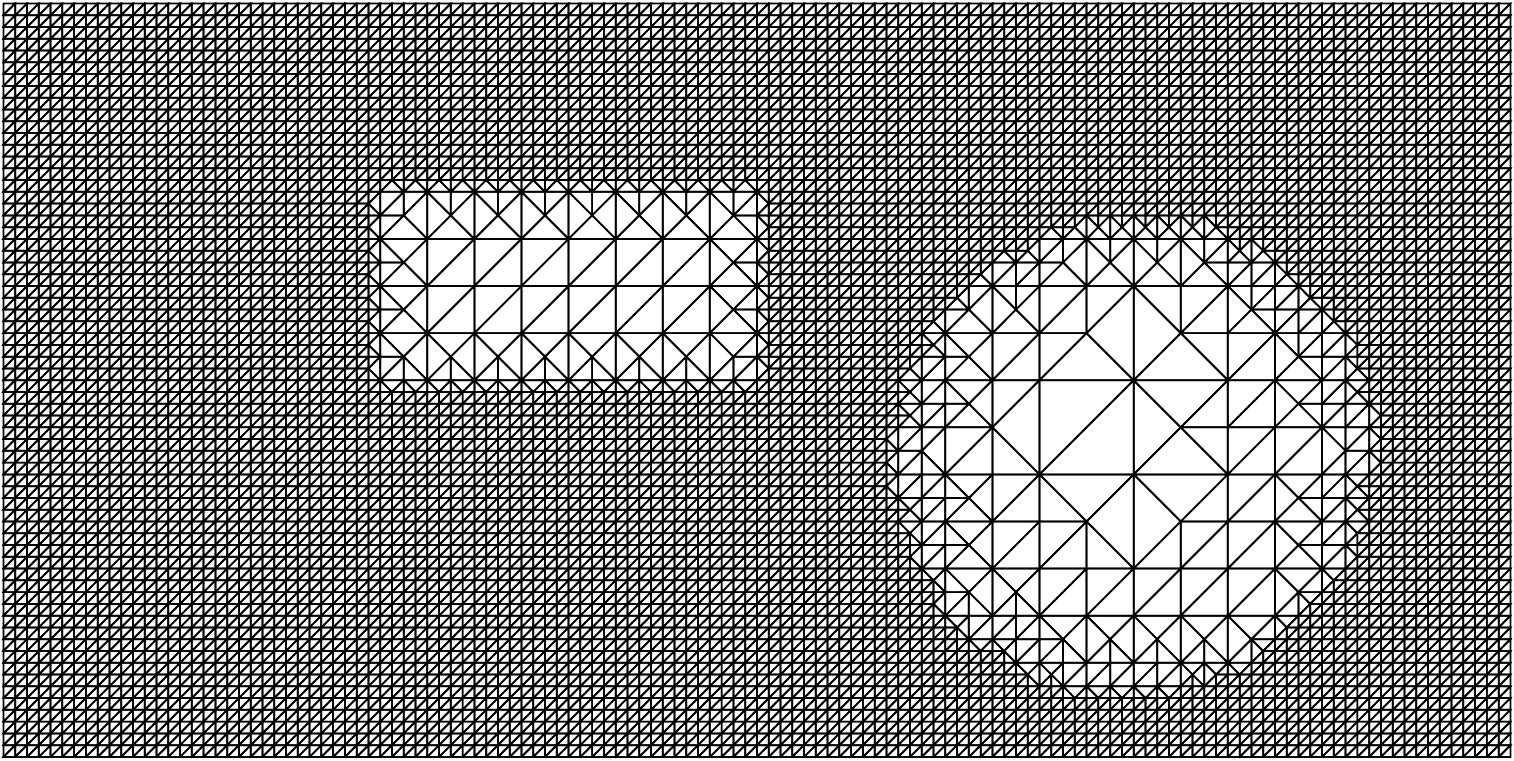}
\end{minipage}
\caption{Local coarsening.}
\label{fig:local}
\end{figure}

So far, we have used coarsening by marking all elements for coarsening. This time, we want to show that our algorithm can be used in an adaptive setting, cf. \texttt{example2/} or \texttt{example3/} or more generally, for local coarsening. To this end, we define a discrete point set and mark elements for coarsening that include this point, see Section~\ref{sect:minimalexample}. We start with a fine triangulation and proceed with this marking strategy and our coarsening algorithm. Figure~\ref{fig:local} shows possible local coarsened meshes. 

In summary, our proposed coarsening algorithm can be used in an interplay of refinement and coarsening, can coarsen locally and recover the initial triangulation - although not quite as efficiently. However, the latter point does usually not play a major role, as only a few coarsening steps are integrated in an adaptive routine.

\bibliographystyle{acm}
\bibliography{references}
\end{document}